\documentclass[11pt]{article}

\usepackage{amssymb}
\usepackage{esvect}
\usepackage{latexsym}
\usepackage[tbtags]{amsmath}
\usepackage{amscd}
\usepackage{amsthm}
\usepackage{color}
\usepackage{enumerate}
\usepackage{bm}
\usepackage{bbm}
\usepackage{placeins}
\usepackage{setspace}
\usepackage{adjustbox} 
\usepackage{subfigure} 
\usepackage{lipsum}
\usepackage{mathtools}
\usepackage[titletoc,title]{appendix}
\usepackage[T1]{fontenc}
\usepackage{dsfont}
\usepackage{natbib}
\usepackage{url}
\usepackage[font=small,labelfont=bf]{caption}
\usepackage{authblk} 
\usepackage[hyperfootnotes=false]{hyperref} 
\hypersetup{
colorlinks = true,
linkcolor = blue, 
filecolor = magenta,
urlcolor = blue, 
citecolor = blue 
}
\usepackage{cleveref}
\usepackage{sectsty}
\usepackage{fullpage}
\usepackage[symbol,multiple]{footmisc}
\usepackage{tikz}
\usepackage{booktabs}
\usepackage{multirow}
\usepackage[textsize=tiny]{todonotes}

\usepackage{colortbl,xcolor}

\definecolor{vlightorange}{rgb}{1,0.97,0.9}       
\definecolor{lightorange}{rgb}{1,0.9,0.75}
\definecolor{midorange}{rgb}{1,0.8,0.6}      
\definecolor{strongorange}{rgb}{1,0.65,0.4}  
\definecolor{deeporange}{rgb}{1,0.5,0.2} 

\definecolor{vlightgreen}{rgb}{0.95, 1.0, 0.95}
\definecolor{lightgreen}{rgb}{0.85, 1.0, 0.85}
\definecolor{midgreen}{rgb}{0.75, 0.95, 0.75}
\definecolor{stronggreen}{rgb}{0.65, 0.90, 0.65}
\definecolor{deepgreen}{rgb}{0.55, 0.85, 0.55}

\definecolor{vlightblue}{rgb}{0.93, 0.97, 1}
\definecolor{lightblue}{rgb}{0.80, 0.90, 1}
\definecolor{midblue}{rgb}{0.60, 0.80, 1}
\definecolor{strongblue}{rgb}{0.40, 0.65, 0.95}
\definecolor{deepblue}{rgb}{0.25, 0.50, 0.85}

\newcommand{\mE}{\mathbb{E}}

\newcommand{\mV}{\mathrm{Var}}

\newcommand{\bX}{\bm{X}}

\newcommand{\bZ}{\bm{Z}}

\newcommand{\given}{\,|\,}

\newtheorem{theorem}{Theorem}[section]   
\newtheorem{lemma}[theorem]{Lemma}

\newtheorem{corollary}[theorem]{Corollary}

\DeclarePairedDelimiter\ceil{\lceil}{\rceil}
\DeclarePairedDelimiter\floor{\lfloor}{\rfloor}
\DeclareMathOperator*{\argmax}{arg\,max}
\DeclareMathOperator*{\argmin}{arg\,min}
\DeclareMathOperator*{\sargmin}{sargmin}

\allowdisplaybreaks

\title{Locally minimax optimal confidence sets for the best model}

\author[1]{Ilmun Kim}
\author[2,3]{Aaditya Ramdas}
\affil[1]{Department of Mathematical Sciences, KAIST\\}
\affil[2]{Department of Statistics and Data Science, Carnegie Mellon University\\}
\affil[3]{Machine Learning Department, Carnegie Mellon University}
\affil[ ]{\texttt{ilmunk@kaist.ac.kr \ aramdas@cmu.edu  }}

\begin{document}
\maketitle

\begin{abstract}
	This paper tackles a fundamental inference problem: given $n$ observations from a distribution $P$ over $\mathbb{R}^d$ with unknown mean $\bm{\mu}$, we must form a confidence set for the index (or indices) corresponding to the smallest component of $\bm{\mu}$. By duality, we reduce this to testing, for each $r$ in $1,\ldots,d$, whether $\mu_r$ is the smallest. Based on the sample splitting and self-normalization approach of~\cite{kim2024dimension}, we propose ``dimension-agnostic'' tests that maintain validity regardless of how $d$ scales with $n$, and regardless of arbitrary ties in $\bm{\mu}$. Notably, our validity holds under mild moment conditions, requiring little more than finiteness of a second moment, and permitting possibly strong dependence between coordinates. In addition, we establish the \emph{local} minimax separation rate for this problem, which adapts to the cardinality of a confusion set, and show that the proposed tests attain this rate. Furthermore, we develop robust variants that continue to achieve the same minimax rate under heavy-tailed distributions with only finite second moments. While these results highlight the theoretical strength of our method, a practical concern is that sample splitting can reduce finite-sample power. We show that this drawback can be substantially alleviated by the multi-split aggregation method of \citet{guo2025rank}. Finally, empirical results on simulated and real data illustrate the strong performance of our approach in terms of type I error control and power compared to existing methods. 
\end{abstract}

\section{Introduction}
\label{sec:intro}
Suppose that $\bX_1, \ldots, \bX_{2n}$ are i.i.d.~random vectors in $\mathbb{R}^d$, $d \geq 2$, with unknown distribution $P$ and mean $\bm{\mu} \coloneqq (\mu_1, \ldots, \mu_d)^\top$. Denoting $[d] \coloneqq \{1, \ldots, d\}$,  the goal of discrete argmin inference is to form a confidence set for 
\[\Theta= \Theta(P) \coloneqq \argmin_{k \in [d]} \mu_k,\]
which is the set of all coordinates whose mean equals the smallest in $\bm{\mu}$—this problem is conceptually equivalent to \emph{discrete argmax inference} since a confidence set for $\argmax_{k \in [d]} \mu_k$ can be obtained by simply negating the samples. Apart from being a fundamental and easy-to-state problem, discrete argmin inference has modern applications. For example, suppose that we have $d$ pre-trained black-box machine learning models (like large language models released by different companies), and we would like to choose the best one(s) for some particular task. To this end, we can evaluate these models on $2n$ unseen i.i.d.\ test data points (from the task distribution) using some task-appropriate loss function, and let $X_{k,i}$ denote the loss of the $k$-th model ($k\in[d]$) on the $i$-th data point ($i \in [2n]$). Then, discrete argmin inference corresponds to identifying the model(s) with minimum risk (expected loss). Because the same test data points are being used for all models, and because the models may be similar, it is important that we allow for strong correlations between the coordinates. Another natural application includes identifying the best treatment(s) amongst many, in a multi-armed randomized clinical trial. Since we may be interested in comparing a large number of models or treatments,
we tackle this problem under high-dimensional settings where the ambient dimension $d$ may vary with the sample size\footnote{We refer to $n$ as the sample size for notational convenience, although the total number of observations is $2n$. This choice simplifies the presentation in later sections involving sample splitting.} $n$; to reflect this dependence explicitly, we denote it by $d_n$, though we omit the subscript when the distinction is not essential.

Noting the duality between confidence sets and hypothesis tests, a large part of the paper will focus on solving the following dual testing problem: given some fixed $r \in [d]$, we test the null and alternative hypotheses given by
\begin{align} \label{Eq: original hypotheses}
	H_0: r \in \Theta  \quad \text{versus} \quad H_1: r \notin \Theta.
\end{align}
Let $\psi_r: \{\bX_1, \ldots, \bX_{2n}\} \to \{0,1\}$ denote a test function that rejects the null hypothesis $H_0: r \in \Theta$ when $\psi_r = 1$. Our objective is then two fold: (a) to construct a test that controls the type I error rate at a nominal level $\alpha \in (0,1)$, and (b) to achieve high (and potentially optimal) power over a broad class of distributions; higher test power will yield a smaller confidence set. 
For (a), we will develop a test that remains asymptotically valid (as $n\to\infty$) regardless of the relationship between the dimension $d$ and the sample size $n$. Such a test is referred to as \emph{dimension-agnostic} (DA), as formalized by \citet{kim2024dimension}. While the DA property can be trivially satisfied without regard to power, the real challenge lies in achieving both DA validity and minimax-optimal power under the alternative. 
We achieve this by adapting the versatile ``sample splitting plus self-normalization'' approach of~\cite{kim2024dimension}, that has been adapted to many other problems since its first preprint appeared in 2020.

Although sample splitting is crucial for DA validity, it entails finite-sample costs in power and stability. To mitigate these issues, we follow the multi-split aggregation method of \citet{guo2025rank}, which repeats the procedure across random splits and aggregates the results; see \Cref{Sec: simulation}.

With these considerations in place, we now formalize our objectives in terms of both dual and primal goals.

\paragraph{Formal (dual) goal.} Let $\mathcal{P}_n$ denote a generic class of distributions with dimension $d_n$. 
Let $\mathcal{P}_{0,r} \subseteq \mathcal{P}_n$ denote those distributions under which $\mu_r$ is the smallest, meaning that the null hypothesis is true. We seek to ensure the following DA control of the type I error:
\begin{align} \label{Eq: DA control}
\limsup_{n \to \infty} \sup_{P \in \mathcal{P}_{0,r}} P(\psi_r = 1) \leq \alpha, \quad \text{regardless of the sequence } (d_n)_{n=1}^\infty.
\end{align}
We also want to ensure that the test has high power under the alternative. That is, when $\mu_r$ is not among the smallest and the gap between $\mu_r$ and the smallest mean is sufficiently large, the test should be able to detect this with high probability. We refer to \Cref{Sec: Power Analysis} for a technical formulation of this power requirement, especially our nuanced goal of local minimax optimality.

\paragraph{Formal (primal) goal.}
We seek a set 
$\widehat{\Theta}$ which is an asymptotically valid DA confidence set for the argmin, satisfying
\begin{align}
	\label{eq:DA-confidence-set}
	\liminf_{n \to \infty} \inf_{P \in \mathcal{P}_n} \inf_{r \in \Theta(P)} P\bigl(r \in \widehat{\Theta}\bigr) \geq 1 - \alpha, \quad \text{regardless of the sequence } (d_n)_{n=1}^\infty.
	\end{align}
We will also later give tight conditions under which $P(r \notin \widehat \Theta) \to 1$ for all $r \notin \Theta$.

As mentioned earlier, our testing results can be inverted to yield such confidence sets for the argmin. We simply run $d$ such DA tests $\psi_1,\dots,\psi_d$ and then let \[
\widehat{\Theta} \coloneqq \{k \in [d] : \psi_k = 0\}\] 
denote the set of indices not rejected by the corresponding tests. This duality between testing and confidence set construction provides a principled alternative to classical methods for constructing confidence sets for the argmin index.

We highlight that the above notion of coverage in \eqref{eq:DA-confidence-set} is uniform in the distributional sense (since the $\inf_P$ follows, rather than precedes, the $\liminf_n$) but pointwise in $\Theta$. To elaborate the latter point, it is helpful to consider the following three types of coverage discussed in the literature:
\begin{enumerate}
	\item \emph{Weak coverage}: At least one element of $\Theta$ is covered, i.e., $P(\Theta \cap \widehat{\Theta} \neq \emptyset) \geq 1 - \alpha$;
	\item \emph{Pointwise coverage}: Every element of $\Theta$ is covered, i.e., $\inf_{r \in \Theta} P(r \in \widehat{\Theta}) \geq 1 - \alpha$;
	\item \emph{Uniform coverage}: The entire set $\Theta$ is covered, i.e., $P(\Theta \subseteq \widehat{\Theta}) \geq 1 - \alpha$.
\end{enumerate}
Our main focus is on achieving pointwise coverage, which lies between weak and uniform coverage in strength—it is more stringent than weak coverage but less demanding than uniform coverage. Notably, all three notions coincide when $\Theta$ is a singleton. In practice, the choice among these guarantees reflects a trade-off between statistical validity and inferential power, and the most appropriate criterion may vary depending on application. 

While our proposed method is designed to attain pointwise coverage, we also show that a simple yet non-trivial modification leads to confidence sets with uniform coverage detailed in \Cref{Sec: MCS}.

\paragraph{Related work.} 
The most directly related works to ours—which we will compare to empirically and theoretically—are the model confidence set of \cite{hansen2011model}, the bootstrap approach of \cite{mogstad2024inference} and the cross-validation plus privacy approach of~\cite{zhang2024winners}. The first of these targets uniform coverage, and it tends to yield very wide sets in practice (and is also extremely slow to run). The second paper targets a slightly different rank inference problem which can be tweaked to yield both a pointwise and uniform coverage solution for our problem, while the third approach also targets pointwise coverage like us. Empirically, both papers tend to  perform worse than our method  across a range of settings. Theoretically, we prove that our approach is locally (and thus globally) minimax optimal; in contrast, the second paper did not study efficiency, while existing theoretical results for the third fall short of minimax optimality (both globally and locally). We expand on these works, and many more, in the broader literature survey below.

The argmin inference problem has a long-standing history in statistics and related fields, dating back at least to the work of \citet{bechhofer1954single,gupta1956decision}, with further developments documented in classical texts such as \citet{gibbons1977selecting,gupta1979multiple}. Although a comprehensive review would take up too much space, we highlight several representative contributions that help situate our work in the broader literature. Early work in this area primarily focused on constructing confidence intervals for the argmin index under parametric assumptions, such as normality or known error distributions \citep[e.g.,][]{gupta1965some,dudewicz1970confidence,nelson2001comparisons,boesel2003using}. In particular, \citet{gupta1965some} proposed early solutions to argmin inference by developing multiple decision procedures for selecting the index with the smallest mean among several normal populations. Focusing on  pointwise coverage, \citet{futschik1995confidence} proposed a two-stage selection procedure that improves upon the subset selection method of \citet{gupta1965some}, though their approach still relies on certain conditions for error distributions and independence among the coordinates.

More recent developments have adopted nonparametric or model-agnostic approaches. \citet{hall2009using} proposed bootstrap-based methods to quantify uncertainty in empirical rankings, including the $m$-out-of-$n$ and independent-component bootstrap to address issues of inconsistency and dependence. While their approach provides valuable insights into ranking variability, it does not directly target argmin inference or provide formal confidence sets for the best-performing index. 

\citet{xie2009confidence} addressed inference in the presence of ties and near-ties by constructing marginal confidence intervals for population ranks using smooth rank estimators and nonstandard bootstrap procedures. Their method improves upon conventional bootstrap intervals, offering better coverage properties under ties and near ties. However, their framework is designed primarily for a fixed number of groups and relies on a smoothing parameter that must be carefully chosen. Pursuing similar goals, \citet{mogstad2024inference} proposed procedures for constructing marginal and simultaneous confidence sets for ranks using valid pairwise comparisons under weak assumptions.  While their method accommodates heteroskedasticity and ties, it does not provide a detailed analysis of power (equivalently, the expected length of the confidence set), and its performance in high-dimensional settings remains unexplored. A related strand of work is the model confidence set (MCS) framework proposed by \citet{hansen2011model}, which constructs a confidence set for the best-performing model under a user-specified loss function. This framework aims to achieve uniform coverage guarantees, but doing so incurs high computational costs (see e.g., \Cref{tab:execution_time_summary}) and often yields procedures with limited power in practice. Additionally, the MCS approach lacks a formal power analysis and does not pursue optimality in distinguishing small differences between competing models. Building on \citet{hansen2011model}, \citet{arnold2024sequential} developed a sequential MCS procedure with time-uniform coverage guarantees, but their method is limited to bounded score functions.

In contrast to these nonparametric approaches, \citet{fan2024ranking} developed a parametric framework for rank inference in multiway comparison designs based on a generalized Plackett--Luce model. Their method focuses on estimating latent ranking parameters from observed choices and achieves optimal convergence rates for individual ranks. However, it relies on a specific model assumption and is designed for a fixed number of groups.

A separate line of work has focused on post-selection inference, which aims to provide valid inference after a data-driven selection step. In this context, \citet{hung2019} introduced a selective inference framework for verifying top ranks in exponential family models via pairwise testing, though their method is restricted to a specific model class and requires tie-breaking to enforce a unique top rank. \citet{sood2024selective} proposed a conceptually unifying framework for selective inference via p-values, which is demonstrated in the context of inference on winners and rank verification. However, the application of their framework is limited to exponential family models or independent p-values, and focuses on validity over efficiency. More recently, \citet{goldwasser2025gaussian} introduced selective inference procedures for verifying the winner and top-$K$ ranks under independent but heteroskedastic Gaussian data, and \citet{sood2025powerful} extended this line of work to settings with arbitrary Gaussian covariance structures. Nevertheless, both approaches rely on the assumption of Gaussianity with a known covariance matrix, which is a strong and often unrealistic requirement in practice. Finally, \citet{painsky2025near} analyze multinomial benchmark rankings under a \emph{fixed} category size, which is an orthogonal setting to our high-dimensional mean-comparison framework in which \(d=d_n\) may grow with \(n\).

Recent work of \citet{zhang2024winners} proposed a general framework for argmin inference in high-dimensional settings with an emphasis on pointwise coverage. Their approach combines cross-validation with exponentially weighted comparisons to construct valid confidence sets for the argmin index. It is model-agnostic and accommodates ties, near ties, and complex dependence structures, making it broadly applicable across diverse data settings. However, the procedure requires careful tuning—such as the choice of weighting parameters and cross-validation strategy—which may influence its practical performance. While the method performs well in many settings, our empirical results in \Cref{fig:type1} suggest that its validity may be sensitive to the problem context, particularly in maintaining type I error control. Moreover, as shown in \Cref{fig:high-dim}, their method exhibits significant power loss in certain regimes, indicating that the test may not achieve a minimax separation rate and highlighting the need for further research to improve its performance. Finally, their theoretical guarantees are established under the assumptions of uniformly bounded data, which is very light-tailed, whereas our results extend to heavy-tailed data, requiring slightly more than existence of a second moment.

The first of our methods is related to a proposal in the latest version of \citet[][Section 4.5]{takatsu2025bridging}, which was done in parallel to our work. Both works utilize the sample-splitting and self-normalization techniques of~\cite{kim2024dimension} to establish DA validity with the pointwise coverage guarantee. While their work establishes the validity of the confidence set, it offers only a brief discussion without a comprehensive theoretical or empirical investigation. In contrast, we provide a thorough analysis including establishing local minimax optimality and empirical evaluations to other methods. Further, we also propose a novel noise-adjusted method that can substantially improve the power under heteroskedasticity, along with additional variants that are robust to heavy-tailed data (see \Cref{Sec: robust}),  both of which are also proven to be locally minimax optimal (the first against light-tailed data, the second against heavy-tailed data).

A related connection arises from the best-arm identification problem in the multi-armed bandit literature~\citep[e.g.,][Chapter 33]{lattimore2020bandit}, where the goal is to identify the most favorable arm based on sample data. However, most bandit methods emphasize sequential decision-making (sampling different coordinates adaptively) rather than fixed-sample inference or confidence set construction. Nonetheless, insights from this literature may inform future developments in rank and argmin inference.

\paragraph{Our contributions.} With the prior work in view, we develop a novel method for argmin inference that satisfies the following key desiderata:
\begin{itemize}
	\item[(i)] \emph{Dimension-agnostic performance}: valid in both low- and high-dimensional settings, without relying on dimension-specific assumptions, and requiring only mild moment conditions;
	\item[(ii)] \emph{Powerful inference}: power that adapts to the cardinality of the confusion set in \eqref{Eq: confusion set} that determines the difficulty of the problem and attains local minimax rates across different regimes;
	\item[(iii)] \emph{Robustness to data characteristics}: accommodating ties and near ties in the mean vector, and remaining valid under strong dependence among components of $\bX$;
	\item[(iv)] \emph{Model-agnostic and tuning-free implementation}: applicable without parametric model assumptions and requiring no (non-trivial or difficult to set) tuning parameters.
\end{itemize}
To the best of our knowledge, no existing method simultaneously satisfies all of these arguably natural desiderata. Our proposed framework is designed to fill this gap.

While our approach builds on the fundamental principle of sample splitting and self-normalization formalized by \citet{kim2024dimension}, it goes beyond a straightforward extension. The discrete argmin inference problem poses unique challenges, particularly the efficient and robust estimation of the runner-up index so as to attain local minimax optimality under minimal assumptions. We address this with explicit selectors, including a noise-adjusted rule that improves power under heteroskedasticity, and robust variants that retain optimality under heavy-tailed distributions. Finally, we develop a new two-step procedure for constructing dimension-agnostic confidence sets, which, for the first time, achieves uniform coverage in high-dimensional settings; see \Cref{Sec: MCS}.

\paragraph{Organization.} The remainder of this paper is organized as follows. In \Cref{Sec: DA argmin test}, we present the proposed DA method, which ensures asymptotic validity under minimal conditions. In \Cref{Sec: Power Analysis}, we derive the minimax separation rate for argmin inference and show that our proposed tests achieve this rate. In \Cref{Sec: robust}, we introduce a robust variant of the initial proposal that achieves the same separation rate under weaker moment conditions. \Cref{Sec: MCS} proposes and analyzes DA model confidence sets with uniform coverage. \Cref{Sec: simulation} presents empirical results demonstrating the competitive performance of the proposed method compared to existing approaches. We conclude in \Cref{Sec: conclusion} by summarizing the paper and discussing potential directions for future research. The omitted proofs and technical results are provided in \Cref{App: Proof}, and additional simulation restuls are presented in \Cref{Sec: additional-simulation}.

\paragraph{Notation.} We use boldface letters (e.g., $\bm{\mu}$, $\bm{\Sigma}$) te denote vectors and matrices, and regular (non-bold) letters for scalars. The operators $\vee$ and $\wedge$ denote the maximum and minimum, respectively, and the symbol $\bm{e}_k$ denotes the $k$-th standard basis vector in $\mathbb{R}^d$. Following convention, the standard normal cumulative distribution function is denoted by $\Phi(\cdot)$, and $N(\bm{\mu}, \bm{\Sigma})$ refers to a multivariate normal distribution with mean vector $\bm{\mu}$ and covariance matrix $\bm{\Sigma}$. The symbol $\bm{I}_d$ denotes the $d \times d$ identity matrix.  $|\mathbb{S}|$ denotes the  cardinality of a set $\mathbb{S}$. We use $o(1)$ to denote a sequence that tends to zero as $n \to \infty$.

\section{Dimension-agnostic argmin test} \label{Sec: DA argmin test}
In this section, we introduce our proposed testing procedure for the hypotheses in \eqref{Eq: original hypotheses}, and establish its asymptotic validity. To this end, we adopt the DA approach introduced by \citet{kim2024dimension} to construct a test that remains valid regardless of the behavior of the dimension $d$. Let \[
s \coloneqq \sargmin_{k \in [d] \setminus \{r\}} \mu_k,
\] 
where `$\sargmin$' denotes the \emph{smallest} index attaining the minimum value (that is, the smallest index in the set $\argmin_{k \in [d] \setminus \{r\}} \mu_k$). This allows us to
reformulate the original hypotheses in \eqref{Eq: original hypotheses} as
\begin{align*}
	H_0: \mu_r - \mu_s \leq 0 \quad \text{versus} \quad H_1: \mu_r - \mu_s > 0,
\end{align*}
which simply determines the positivity of $\mu_r - \mu_s$. When $s$ is known, this problem can be tackled using a standard one-sided $t$-test. However, the complexity arises when $s$ is unknown and needs to be estimated from the data. To handle this, we use a sample splitting strategy where one subset is used to estimate $s$ (\emph{model selection}), and another is used to construct a test (\emph{inference}), typically using some form of self-normalization. 

This ``sample splitting plus self normalization'' is a fundamental principle of the DA approach. After its introduction in~\cite{kim2024dimension}, this technique for DA inference (as opposed to just inference) has been successfully applied to various high-dimensional inference problems~\citep[e.g.,][]{liu2022multiple,shekhar2022permutation,shekhar2023permutation,gao2023dimension,martinez2023efficient,zhang2024another,lundborg2024projected,liu2024projection,zhang2024novel,takatsu2025bridging}. By extending this framework to the discrete argmin inference problem, our work ensures asymptotic validity under mild moment conditions and achieves minimax-optimal power across both low- and high-dimensional regimes, even for heavy-tailed data. 

The next subsection describes the proposed DA argmin test in detail. 

\subsection{Our procedure} \label{Sec: Procedure}
Before presenting our procedure, we first describe a natural approach to the argmin inference, which uses the full sample mean vector $\overline{\bX} \coloneqq (\overline{X}_1,\ldots,\overline{X}_d)^\top = \frac{1}{2n} \sum_{i=1}^{2n} \bX_i$ without sample splitting. Specifically, this method computes the maximum of the $d-1$ one-sided $t$-statistics given by
\begin{align*}
\max_{k \in [d]\setminus {r}} \frac{\overline{X}_r - \overline{X}_k}{\widehat{\sigma}_{r,k}},
\end{align*}
where $\widehat{\sigma}^2_{r,k}$ denotes an estimator of the variance for the difference $\overline{X}_r - \overline{X}_k$. While intuitive, this approach involves \emph{double dipping} as the same data are employed for both identifying the most significant component and performing inference, complicating calibration particularly in high-dimensional scenarios. Bootstrap-based calibration methods, such as those employed by \citet{mogstad2024inference}, are a viable option to address this issue. However, their methodology is limited to fixed-dimensional settings and computationally expensive due to the need for repeated resampling. Moreover, \citet{mogstad2024inference} do not provide theoretical guarantees regarding statistical efficiency or adaptivity to the intrinsic difficulty of argmin inference.

Notably, our proposed approach is also based upon the same underlying statistic. However, we circumvent the calibration difficulty by splitting the dataset into two independent halves: the first half is used for selecting the component that maximizes the statistic, and the second half is dedicated to conducting the inference through a one-sided $t$-statistic evaluation. As we illustrate later, this two-step procedure is not only simple to implement but also leads to a test that is both dimension-agnostic and locally minimax optimal. 

We now proceed to a detailed presentation of our procedure.

\paragraph{DA approach.} Denote the sample means as $\overline{\bX}^{(1)} = (\overline{X}^{(1)}_1,\ldots,\overline{X}^{(1)}_d)^\top = \frac{1}{n} \sum_{i=1}^{n} \bX_i$ and $\overline{\bX}^{(2)} = (\overline{X}^{(2)}_1,\ldots,\overline{X}^{(2)}_d)^\top = \frac{1}{n} \sum_{i=n+1}^{2n} \bX_i$, which are constructed from the first and second halves of the samples, respectively. To address the argmin inference problem, we propose a simple two-step procedure that separates the selection and inference stages:
\begin{enumerate}
	\item \textbf{Selection.} Estimate the argmin $s$ using the second half of samples. We propose two different approaches for this purpose. The first estimator is the plug-in estimator, defined as 
	\begin{align*}
		\widehat{s}_{\mathrm{plug}} \coloneqq \sargmin_{k \in [d] \setminus \{r\}} \overline{X}^{(2)}_k,
	\end{align*}
	which directly selects the index corresponding to the smallest sample mean in the second half of the data. Alternatively, we propose a noise-adjusted estimator that accounts for the potentially differing noise level associated with each component. Denote 
    \[
    \bm{\gamma}_k \coloneqq \boldsymbol{e}_r - \boldsymbol{e}_k,
    \]
    where we recall that $\boldsymbol{e}_r$ and $\boldsymbol{e}_k$ are the $r$-th and $k$-th standard basis vectors in $\mathbb{R}^d$, respectively. The noise-adjusted estimator is defined as
	\begin{align*}
		\widehat{s}_{\mathrm{adj}} \coloneqq \sargmin_{k \in [d] \setminus \{r\}} \frac{\overline{X}^{(2)}_k - \overline{X}^{(2)}_r}{\sqrt{\bm{\gamma}_k^\top \widehat{\bm{\Sigma}}^{(2)} \bm{\gamma}_k}\vee \kappa},
	\end{align*}
	where $\kappa > 0$ is a small constant (set to $10^{-8}$ in our experiments) included to prevent instability in variance estimation. The matrix $\widehat{\bm{\Sigma}}^{(2)}$ above is the sample covariance matrix computed from $\bX_{n+1},\ldots,\bX_{2n}$. This noise-adjusted estimator essentially finds an index that maximizes a signal-to-noise ratio, defined as the mean difference divided by the standard deviation, rather than the mean difference itself. 
	\item \textbf{Inference.} Given $\widehat{s} = \widehat{s}_{\mathrm{plug}}$ or $\widehat{s}=\widehat{s}_{\mathrm{adj}}$, we determine whether the mean difference 
	\begin{align*}
		\overline{X}^{(1)}_r - \overline{X}^{(1)}_{\widehat{s}} = \bm{\gamma}_{\widehat{s}}^\top \overline{\bX}^{(1)} 
	\end{align*}
	is significantly positive. Specifically, we reject the null hypothesis if 
	\begin{align*}
        \sqrt{n}\bm{\gamma}_{\widehat{s}}^\top \overline{\bX}^{(1)} > z_{1-\alpha} \sqrt{\bm{\gamma}_{\widehat{s}}^\top \widehat{\bm{\Sigma}}^{(1)} \bm{\gamma}_{\widehat{s}}},	
	\end{align*}
	where $\widehat{\bm{\Sigma}}^{(1)}$ is the sample covariance matrix based on $\bX_1,\ldots,\bX_n$ and $z_{1-\alpha}$ is the $1-\alpha$ quantile of $N(0,1)$ with $\alpha \in (0,1)$. 
\end{enumerate}
We refer to the test derived from this two-step procedure as the DA argmin test. A few remarks are in order about the procedure:
\begin{itemize}
	\item In essence, the proposed argmin test is a standard one-sided $t$-test to determine whether $\mu_r - \mu_{\widehat{s}}$ is positive. Under the null, $\mu_r - \mu_{\widehat{s}} \leq 0$ holds for any choice of $\widehat{s} \in [d]\!\setminus\!\{r\}$, and thus the test maintains its asymptotic validity even if $\widehat{s}$ is incorrectly selected. On the other hand, under the alternative, $\widehat{s}$ is expected to satisfy $\mu_r - \mu_{\widehat{s}} > 0$ with high probability. This positive gap leads to significant power in detecting deviations from the null. 
	\item Sample splitting plays a crucial role in this framework. Without sample splitting, the samples are reused for both selection and inference, which results in strongly dependent summands in the test statistic. This strong dependency structure breaks the conditions for central limit theorem and leads to invalid inference.
	\item Despite its central role, sample splitting has drawbacks mentioned earlier: (i) the reduced sample size for both selection and inference can lower (practical) power and (ii) the results may vary with different random splits. To mitigate these issues, we adopt the multi-split aggregation method of \citet{guo2025rank} as described in \Cref{Sec: simulation}. This strategy reduces randomness and improves power by using the samples more efficiently.
    \item Recall that our test can be easily inverted (by repeating it for each coordinate) to produce a DA confidence set as outlined in~\eqref{eq:DA-confidence-set}.
\end{itemize}
In the following sections, we examine the theoretical properties of the DA argmin test, focusing on its asymptotic validity and power analysis. These theoretical results apply to both selection procedures, namely $\widehat{s}_{\mathrm{plug}}$ and $\widehat{s}_{\mathrm{adj}}$, and thus we denote either estimator simply by $\widehat{s}$ whenever the distinction is not necessary.

\subsection{Asymptotic validity} \label{Sec: Asymptotic Validity}
To establish the asymptotic validity of the proposed argmin test, we impose a mild moment condition on the contrasts $W_1,\ldots,W_d$ where each 
\begin{align} \label{Eq: definition of Wk}
	W_{k} \coloneqq \bm{\gamma}_k^\top (\bX - \bm{\mu})
\end{align}
represents the difference between the $r$-th and the $k$-th centered coordinates. To motivate the form of our condition, consider a class of null distributions $\mathcal{P}_{0,r}$ for $H_0:r \in \Theta$ and note that a standard Berry--Esseen bound for normalized sums (of i.i.d.~copies of the random variable $W_k$) typically involves the third absolute moment. In particular,  for asymptotic normality to hold uniformly over $\mathcal P_{0,r}$, one commonly encountered condition is that
\begin{align*}
	\max_{k \in [d] \setminus \{r\}} \sup_{P \in \mathcal{P}_{0,r}} \mE_P \left[ \frac{|W_k|^3}{n^{1/2} \{\mE_P[W_k^2]\}^{3/2}} \right] = o(1),
\end{align*}
where the maximum over $k$ ensures uniform convergence across all coordinates excluding the target coordinate $r$. Rather than requiring this third-moment condition, we impose a strictly weaker moment condition. Specifically, we assume that 
\begin{align} \label{Eq: moment condition}
	\max_{k \in [d] \setminus \{r\}} M_k \coloneqq \max_{k \in [d] \setminus \{r\}} \sup_{P \in \mathcal{P}_{0,r}} \mE_P \Biggl[ \frac{W_k^2}{\mE_P[W_k^2]}\min \Biggl\{ 1, \frac{|W_{k}|}{n^{1/2} (\mE_P[W_k^2])^{1/2}} \Biggr\} \Biggr] = o(1).
\end{align}
A few remarks on this condition are in order.
\begin{itemize}
	\item The moment condition \eqref{Eq: moment condition} serves a similar role to the remainder term in a Berry--Esseen bound, but with a lighter tail requirement that allows for a broader class of distributions. For example, the $t$-distribution with $3$ degrees of freedom lacks a finite third moment, yet it satisfies the truncated second moment condition in \eqref{Eq: moment condition}. 
	Interestingly, the truncated moment condition is in fact equivalent to Lindeberg's condition for the central limit theorem, which characterizes necessary and sufficient conditions for convergence of general triangular arrays. We establish this equivalence in \Cref{lem: UCLT}. In light of this lemma, our imposition of such a condition should be viewed as both natural and minimal.
	\item We also emphasize that the moment condition~\eqref{Eq: moment condition} is a one-dimensional requirement on the contrasts $W_k=\bm{\gamma}_k^\top(\bX-\bm{\mu})$ and, by itself, places essentially no restriction on the joint dependence among the coordinates of $\bX$. Consequently, the moment condition~\eqref{Eq: moment condition} allows strong dependence between the coordinates of $\bX$. For instance, if $\bX$ follows a multivariate normal distribution, the condition holds for any positive semi-definite covariance matrix, provided that the variance of each $W_k$ is positive.  In particular, it allows the correlations between $\boldsymbol{e}_r^\top \bX$ and $\boldsymbol{e}_k^\top \bX$ to approach one at an arbitrary rate, and the remaining components of $\bX$ (excluding the $r$-th coordinate) to be arbitrarily dependent. More broadly, the same conclusion extends beyond the multivariate normal, for example, to families in which the kurtosis of $W_k$ is uniformly bounded across $k$.
	\item In the moment condition~\eqref{Eq: moment condition}, the maximum over $k$ is taken outside the expectation, which is a considerably weaker requirement than placing the maximum inside. In typical scenarios, the truncated moments $M_k$ are of comparable order across different $k$, so uniform convergence is expected to hold regardless of how $d$ grows with $n$. This observation underlies the dimension‑agnostic property of our proposed test, as stated in \Cref{The: DA validity} below.
	\item We highlight that imposing condition~\eqref{Eq: moment condition} on the vector $\bX$ itself rather than their componentwise difference does not guarantee the asymptotic normality.
\end{itemize}

The asymptotic validity of the DA argmin test over $\mathcal{P}_{0,r}$ follows directly from the theorem below.

\begin{theorem} \label{The: DA validity}
	There exists a constant $C > 0$ such that the following inequality holds
	\begin{align*}
		& \sup_{P \in \mathcal{P}_{0,r}} \sup_{t \in \mathbb{R}} \Bigg| P \Biggl(\frac{\sqrt{n}\bm{\gamma}_{\widehat{s}}^\top \bigl(\overline{\bX}^{(1)} -  \bm{\mu} \bigr) }{\sqrt{\bm{\gamma}_{\widehat{s}}^\top \widehat{\bm{\Sigma}}^{(1)} \bm{\gamma}_{\widehat{s}}}} \leq t\Biggr) - \Phi(t) \Bigg| \leq  \min \biggl\{1, C \!\! \max_{k \in [d] \setminus \{r\}} M_k\biggr\}.
	\end{align*}
    By Assumption~\eqref{Eq: moment condition}, we conclude that the DA argmin test is asymptotically valid uniformly over $\mathcal{P}_{0,r}$ in the sense of \eqref{Eq: DA control}.
\end{theorem}

We highlight again that the validity of the DA argmin test requires only a mild moment condition, slightly stronger than the existence of a finite second moment. This flexibility enables the DA argmin test to remain reliable even in high-dimensional and heavy-tailed settings. In contrast, existing methods often impose more stringent assumptions, such as the uniformly bounded random variable condition in \citet{zhang2024winners} and the parametric assumption in classical approaches~\citep[e.g.,][]{gupta1965some}. More importantly, the proposed procedure attains the minimax separation rate for the argmin inference problem, as detailed in the following section. Finally, \Cref{The: DA validity} applies to \emph{any} data-driven selection rule $\widehat{s}$ computed exclusively from the second half of the data, making the validity guarantee robust to the choice of selection procedure.

\section{Power analysis} \label{Sec: Power Analysis}
\begin{figure}[t]
	\centering
	\includegraphics[width=1\textwidth]{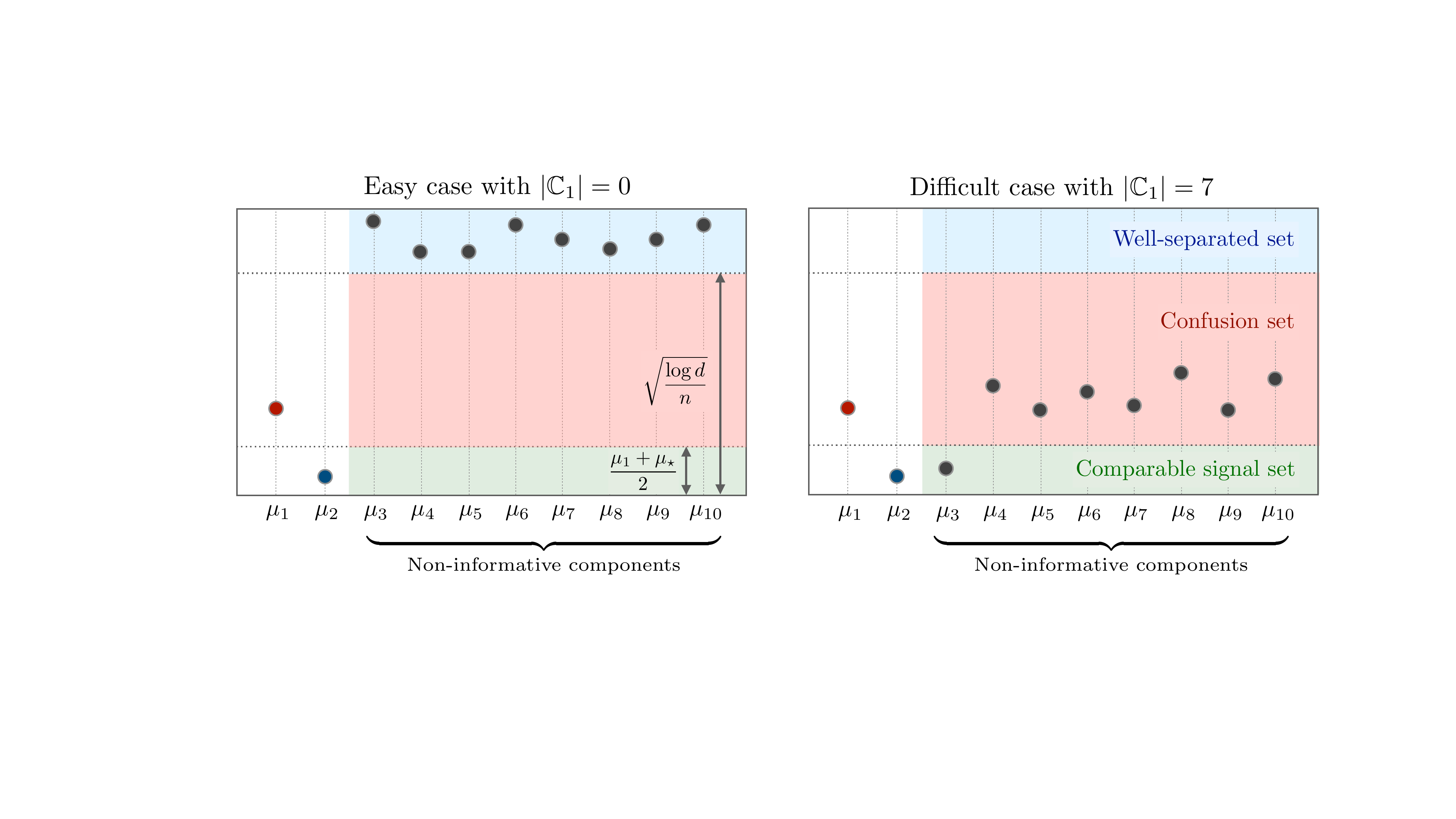}
	\caption{Illustration of the confusion set $\mathbb{C}_r$ with $r=1$. The left panel depicts a scenario with $|\mathbb{C}_1| = 0$ where $\mu_3,\ldots,\mu_{10}$ are sufficiently larger than the minimum $\mu_{\star} = \mu_2$ relative to $\mu_1$, allowing the argmin to be easily identified based on samples. In contrast, the right panel illustrates a scenario with $|\mathbb{C}_1| = 7$ where $\mu_4,\ldots,\mu_{10}$ are closer to $\mu_{\star}$, making it more difficult to distinguish the argmin from the other components. Note that $\mu_3$ on the right is excluded from the confusion set because it violates the lower-bound condition in \eqref{Eq: confusion set}. In this case, $\mu_1 -\mu_3$ is comparable in size to $\mu_1 -\mu_\star$.}
	\label{fig:confusion}
\end{figure}

We next analyze the power of the DA argmin test under the alternative hypothesis. As a first step in our analysis, we introduce the notion of a confusion set, which characterizes the difficulty of the argmin inference problem. 
Denote  \[
\mu_{\star} \coloneqq \min_{k \in [d]} \mu_k,
\] 
and define the set 
\[
\Theta_{-r} \coloneqq \argmin_{k \in [d] \setminus \{r\}} \mu_k.
\]
Under the alternative, $\mu_r$ is not in the argmin set, so $\mu_{\star}=\min_{k \in [d] \setminus \{r\}} \mu_k$ and is attained by every element of $\Theta_{-r}$, implying that $\mu_r > \mu_{\star} = \mu_s$ for all $s \in \Theta_{-r}$. The confusion set for the index $r$ is defined as:
\begin{align} \label{Eq: confusion set}
	\mathbb{C}_r \coloneqq \biggl\{k \in [d] \! \setminus \! (\{r\} \cup \Theta_{-r}) : \frac{\mu_r - \mu_{\star}}{2} \leq \mu_k - \mu_{\star} \leq C_n \sqrt{\frac{\log (d)}{n}} \biggr\}, 
\end{align}
where $C_n$ is any positive sequence such that $C_n \to \infty$ as $n \to \infty$. Here the constant $1/2$ in the lower bound is arbitrary and can be replaced by any constant in $(0,1)$. Note that by construction, $\mathbb{C}_r$ excludes the index $r$, but under the alternative, it also excludes every index $s \in \Theta_{-r}$ because $\mu_s - \mu_{\star}$ equals 0 but the lower bound in \eqref{Eq: confusion set} is positive. See Figure~\ref{fig:confusion} for an illustration.

\paragraph{Remarks on the confusion set $\mathbb{C}_r$:}
\begin{itemize}
	\item To better understand the role of the confusion set, first consider the case where $\mu_k - \mu_{\star} > C_n \sqrt{\log(d)/n}$. In this scenario, $\mu_k$ is sufficiently far from $\mu_{\star}$, making it unlikely for index $k$ to be selected as the sample argmin. Such indices are therefore not problematic for inference and can be effectively disregarded when assessing the difficulty of the argmin inference problem.
	Next, consider the case where $\mu_k - \mu_{\star} < (\mu_r - \mu_{\star})/2$, under which it holds that
	\begin{align*}
	\mu_r - \mu_k > \frac{1}{2}(\mu_r - \mu_{\star}).
	\end{align*}
	In the event that $\widehat{s} = k$, the resulting signal $\mu_r - \mu_{\widehat{s}}$ remains sufficiently large, comparable in magnitude to $\mu_r - \mu_{\star}$ up to a constant factor, thereby allowing reliable detection of the difference between $\mu_r$ and $\mu_{\star}$.
	
	Taken together, these observations suggest that the confusion set $\mathbb{C}_r$ comprises indices for which the signal $\mu_r - \mu_k$ is not large enough to ensure reliable discrimination between $\mu_r$ and $\mu_{\star}$. In other words, the confusion set captures the subset of indices that truly contribute to the difficulty of the argmin inference problem.
	\item The confusion set appearing in \citet{zhang2024winners} is given by
	\begin{align} \label{Eq: confusion set zhang}
		\widetilde{\mathbb{C}}_r \coloneqq \biggl\{k \in [d] \! \setminus \! (\{r\} \cup \Theta_{-r}) : \frac{\mu_r - \mu_{\star}}{2} \leq \mu_k - \mu_{\star} \leq \frac{1}{\lambda} \bigl( \log d + 3 \sqrt{\log V}\bigr) \biggr\},
	\end{align}
	where $\lambda = o(\sqrt{n})$ and $V$ denotes the number of folds in cross-validation. The main difference from ours lies in their upper bound, which is less restrictive than the one in \eqref{Eq: confusion set}. Thus their confusion set is larger than ours, leading to a worse rate.
\end{itemize}

\begin{figure}[t]
	\centering
	\includegraphics[width=.85\textwidth]{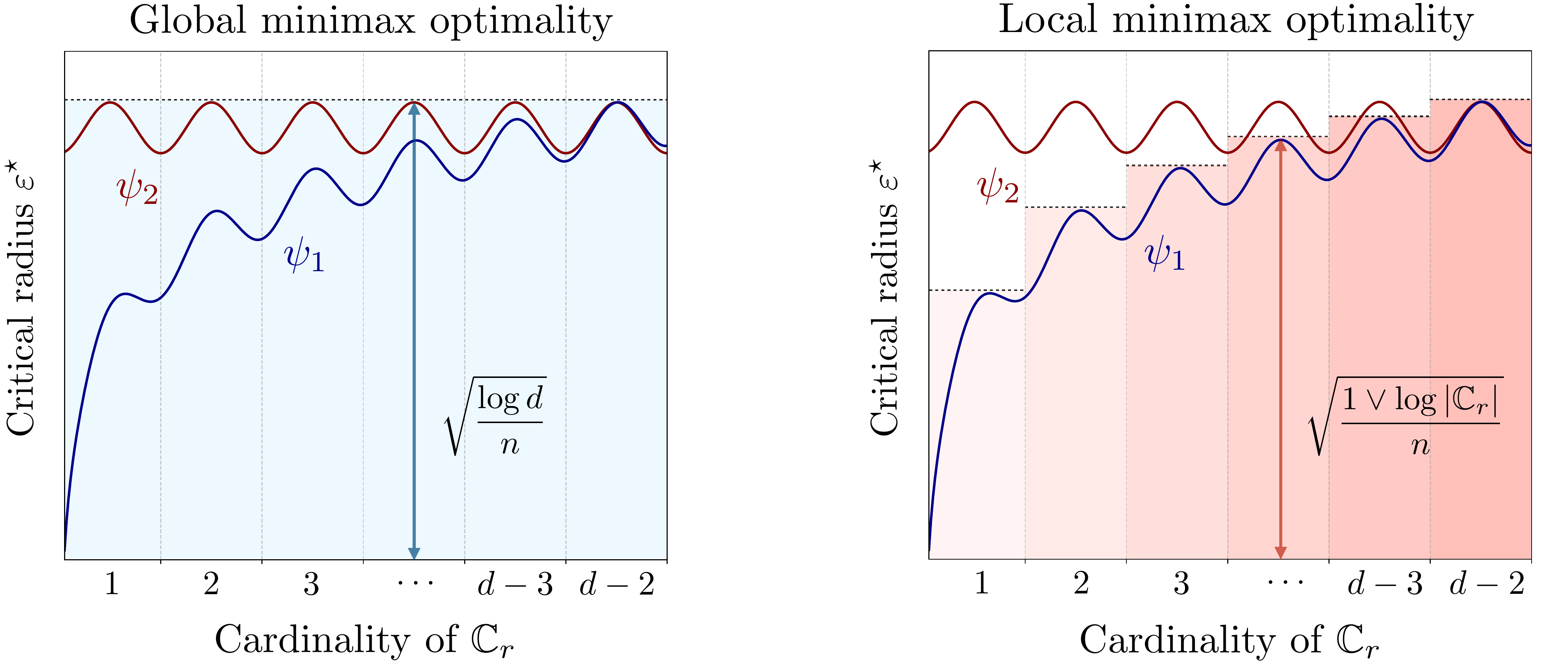}
	\caption{Illustration of global vs.~local minimax optimality.
	\emph{Left}: both tests $\psi_1$ and $\psi_2$ achieve global minimax optimality, with uniform separation rates below the global critical radius $\varepsilon^\star \asymp \sqrt{\log(d)/n}$ (indicated by the dotted line), which is independent of the confusion set size $|\mathbb{C}_r|$. \emph{Right}: only $\psi_1$ achieves local minimax optimality by adapting to the confusion set size through the $|\mathbb{C}_r|$-dependent critical radius $\varepsilon^\star = \sqrt{(1 \vee \log |\mathbb{C}_r|)/n}$ as defined in \eqref{Eq: critical radius} (indicated by the dotted lines).}
	\label{fig:local minimax}
\end{figure}

Having defined the confusion set, we now explain the main objective of this section. Let $\mathcal{P}$ be a collection of distributions where $\bX \sim P \in \mathcal{P}$ is a sub-Gaussian random vector in $\mathbb{R}^d$ with a fixed variance proxy $\sigma^2$. That is, we assume that for every unit vector $v \in \mathbb{R}^d$, the one-dimensional projection $\langle \bm{v}, \bX \rangle$ is a sub-Gaussian random variable with parameter $\sigma^2$; i.e., 
\[\mE \bigl[\exp\bigl(\lambda \langle \bm{v}, \bX \rangle\bigr)\bigr] \leq \exp\bigl(\lambda^2 \sigma^2/2\bigr) \text{ for all } \lambda \in \mathbb{R}.
\]
Note that, in particular, the variance of $\langle \bm{v}, \bX \rangle$ is at most $\sigma^2$ for every unit norm $\bm{v} \in \mathbb{R}^d$.
Now define a class of local alternatives that share the same cardinality of the confusion set as
\begin{align*}
	\mathcal{P}_{1,r}(\varepsilon; \tau) \coloneqq \bigl\{P \in \mathcal{P} : \mu_{r} - \mu_{\star} \geq \varepsilon \ \text{and} \ |\mathbb{C}_r| = \tau \bigr\},
\end{align*}
where $\varepsilon > 0$ is a positive constant and $\tau \in \{0,1,\ldots,d-2\}$. We aim to characterize the condition on $\varepsilon$ under which the asymptotic uniform power of the DA test is one for distributions in $\mathcal{P}_{1,r}(\varepsilon;\tau)$. In particular, we claim that if $\varepsilon$ is sufficiently larger than the critical radius $\varepsilon^\star$ defined as
\begin{align} \label{Eq: critical radius}
	\varepsilon^\star = \varepsilon^\star(\tau)  \coloneqq \sqrt{\frac{1 \vee  \log(\tau)}{n}},
\end{align}
then the DA test has asymptotic power one. Moreover, we show in \Cref{Thm: lower bound} that if $\varepsilon$ is sufficiently smaller than $\varepsilon^\star$, then no asymptotic level-$\alpha$ test can achieve nontrivial uniform power over distributions in $\mathcal{P}_{1,r}(\varepsilon;\tau)$. This implies that the DA argmin test is \emph{locally} minimax optimal: it achieves the best possible separation rate for each fixed confusion set size $\tau$, adapting to the intrinsic difficulty of the problem instance. \Cref{fig:local minimax} illustrates the distinction between global and local minimax optimality.

We formalize and prove these claims in the subsections that follow.

\subsection{Upper bound} \label{Sec: upper bound}
We start with a positive result that characterizes the condition under which the DA argmin test has asymptotic power one. The following result holds for both selection procedures, $\widehat{s}_{\mathrm{plug}}$ and $\widehat{s}_{\mathrm{adj}}$.
\begin{theorem} \label{Thm: upper bound}
	For any $\tau \in \{0,1,\dots,d-2\}$, suppose that $\varepsilon \geq C_n' \varepsilon^\star$ where $C_n'$ is any positive sequence diverging to infinity as $n \to \infty$. Then the asymptotic uniform power of the DA argmin test over $\mathcal{P}_{1,r}(\varepsilon;\tau)$ equals one:
	\begin{align*}
		\lim_{n \to \infty} \inf_{P \in \mathcal{P}_{1,r}(\varepsilon;\tau)} P \Big( \sqrt{n}\bm{\gamma}_{\widehat{s}}^\top \overline{\bX}^{(1)} > z_{1-\alpha} \sqrt{\bm{\gamma}_{\widehat{s}}^\top \widehat{\bm{\Sigma}}^{(1)} \bm{\gamma}_{\widehat{s}}} \Big) = 1.
	\end{align*}
    Since the DA argmin test does not depend on knowledge of $\tau$, Theorem~\ref{Thm: lower bound} implies that it is locally minimax optimal.
\end{theorem}
\Cref{Thm: upper bound} shows that the DA argmin test achieves a uniform separation rate that adapts to the unknown cardinality of the confusion set. In particular, when the cardinality $|\mathbb{C}_r|$ is constant, the test attains the parametric $1/\sqrt{n}$-rate. More generally, the separation rate depends logarithmically on $|\mathbb{C}_r|$, with the worst-case rate being $\sqrt{\log(d)/n}$. A related result by \citet[][Theorem 4.1]{zhang2024winners} shows that their test is powerful when $\mu_r - \mu_{\star}$ is sufficiently larger than $\lambda^{-1}\bigl(\log |\widetilde{\mathbb{C}}_r| + \log\log(d) + \log\log V \bigr)$. Under the assumption $\lambda = o(\sqrt{n})$, which is required for the validity of their procedure, this comparison highlights that our test achieves a sharper (and indeed optimal, as shown in \Cref{Thm: lower bound}) separation rate than that of \citet{zhang2024winners}. We refer to empirical results in \Cref{fig:high-dim} that support this claim. 

\Cref{Thm: upper bound} yields a direct implication for the DA confidence set for $\Theta$, which is constructed by inverting the DA argmin test. Specifically, it implies that any index $r \notin \Theta$ is asymptotically excluded from the DA confidence set, provided that the mean gap $\mu_r - \mu_{\star}$ is sufficiently larger than $\sqrt{(1 \vee \log |\mathbb{C}_r|)/n}$. We formalize this implication in the following corollary.
\begin{corollary} \label{Cor: upper bound}
	For any $\tau \in \{0,1,\dots,d-2\}$, suppose that the $r$-th mean gap satisfies $\mu_r - \mu_\star \geq \varepsilon \geq C_n' \varepsilon^\star$ where $C_n'$ is any positive sequence diverging to infinity as $n \to \infty$. Let $\widehat{\Theta}_{\mathrm{DA}}$ denote the confidence set constructed by inverting the DA argmin test. Then the index $r$ is excluded from $\widehat{\Theta}_{\mathrm{DA}}$ with probability tending to one:
	\begin{align*}
		\lim_{n \to \infty} \inf_{P \in \mathcal{P}_{1,r}(\varepsilon;\tau)} P \big( r \notin \widehat{\Theta}_{\mathrm{DA}} \big) = 1.
	\end{align*}
\end{corollary}

As formally established later in \Cref{Cor: lower bound}, the above result is minimax optimal in the sense that no asymptotically valid confidence set can reliably exclude the index $r \notin \Theta$ when the mean gap $\mu_r - \mu_{\star}$ is sufficiently smaller than $\varepsilon^\star$. 

\subsection{Lower bound}
In this subsection, we establish a lower bound for the separation rate $\varepsilon$ and show that the DA argmin test is minimax rate optimal. Building on this, we further show that the DA confidence set also achieves minimax rate optimality. Recall that $\mathcal{P}_{0,r}$ represents the collection of null distributions satisfying $r \in \Theta$ and the moment condition specified in \eqref{Eq: moment condition}. Let $\Psi_{\alpha}$ be the set of all asymptotic level-$\alpha$ tests over $\mathcal{P}_{0,r}$ defined as
\begin{align*}
	\Psi_{\alpha} = \Psi(\alpha,r) \coloneqq \biggl\{ \psi: \limsup_{n \to \infty} \sup_{P \in \mathcal{P}_{0,r}} P(\psi = 1) \leq \alpha \biggr\}.
\end{align*} 
The following result illustrates that any test in $\Psi_{\alpha}$ cannot achieve a separation rate smaller than $\varepsilon^\star$.

\begin{theorem} \label{Thm: lower bound}
	Let $\alpha \in (0,1/2)$ and $\beta \in (0, 1-2\alpha)$. 
	There exists a constant $c > 0$ that only depends on $\alpha, \beta$ and $\sigma$ such that if $\varepsilon \leq c \varepsilon^\star$, then the asymptotic minimax type II error is at least $\beta$:
	\begin{align*}
		\liminf_{n \to \infty} \inf_{\psi \in \Psi_{\alpha}} \sup_{P \in \mathcal{P}_{1,r}(\varepsilon;\tau)} P(\psi = 0) \geq \beta.
	\end{align*}
\end{theorem}
We emphasize an adaptive nature of this lower bound, which ranges from a parametric $1/\sqrt{n}$-rate to a $\sqrt{\log(d)/n}$-rate depending on the cardinality of the confusion set. Intuitively, when the confusion set is small, the search cost for the argmin index is negligible, allowing the rate to remain parametric. However, as the confusion set grows, the search cost increases, and in the worst-case scenario, the rate degrades to $\sqrt{\log(d)/n}$. The proof of \Cref{Thm: lower bound} builds on this intuition by carefully designing $\bm{\mu}$ to accommodate confusion sets of varying cardinalities.

Let $\mathcal{A}_{\alpha}$ denote the collection of all asymptotic $1-\alpha$ confidence sets for $\Theta$ defined as
\begin{align*}
	\mathcal{A}_{\alpha} = \Bigl\{ \widehat{\Theta}: \liminf_{n \to \infty} \inf_{P \in \mathcal{P}} \inf_{r \in \Theta(P)} P(r \in \widehat{\Theta}) \geq 1 - \alpha \Bigr\}.
\end{align*}
By the duality between confidence sets and tests, \Cref{Thm: lower bound} reveals a fundamental limitation in constructing confidence sets for $\Theta$. For some $r \notin \Theta$, if the mean gap $\mu_r - \mu_{\star}$ is substantially smaller than $\sqrt{(1 \vee \log|\mathbb{C}_r|)/n}$, then no asymptotically valid confidence set can  ensure the exclusion of $r$. This limitation is formalized in the following corollary.

\begin{corollary} \label{Cor: lower bound}
	Let $\alpha \in (0,1/2)$ and $\beta \in (0, 1-2\alpha)$. There exists a constant $c > 0$ that only depends on $\alpha, \beta$ and $\sigma$ such that if $\varepsilon \leq c \varepsilon^\star$, then the worst-case probability of correctly excluding $r \notin \Theta$ across all asymptotically valid confidence sets is at most $1-\beta$:
	\begin{align*}
		\limsup_{n \to \infty} \sup_{\widehat{\Theta} \in \mathcal{A}_{\alpha}} \inf_{P \in \mathcal{P}_{1,r}(\varepsilon;\tau)} P(r \notin \widehat{\Theta}) \leq 1-\beta \quad \text{for any $r \in [d]$.}
	\end{align*}
\end{corollary}

We reiterate that \Cref{Cor: upper bound} and \Cref{Cor: lower bound} taken together establish the local minimax optimality of the DA confidence set at the separation rate $\varepsilon^\star$. We next introduce a robust variant of the DA argmin test that is designed to attain the minimax separation rate under heavy-tailed distributions.

\section{Robust DA argmin test} \label{Sec: robust}

In the previous section, we established that the proposed DA argmin tests (and the DA confidence sets) attain the minimax separation rate under sub-Gaussian assumptions. As a natural next step, we extend these tests to handle heavy-tailed distributions by developing a robust variant. This robust version is specifically designed to retain desirable power guarantees even when the data exhibit outliers or lack sub-Gaussian tails. The central idea is to replace $\widehat{s}$ with a robust alternative that is less sensitive to outliers. 

To this end, we employ the median-of-means (MoM) estimator for estimating the argmin $s$. The MoM estimator, which traces back to \citet{nemirovsky1983problem,jerrum1986random}, has been extensively studied in the literature \citep[e.g.,][]{alon1996space, lerasle2011robust, hsu2016loss, bubeck2013bandits, lugosi2019}. It is defined as the median of the sample means over $V$ disjoint subsets of the data. Formally, let $B_1,\ldots,B_V$ be a partition of $[n]$ into equally sized blocks, each of size $|B_v| = n/V$, and assume $V \leq n/2$. The MoM estimator of $\mu_k$ for $k \in [d]$ is then defined as
\begin{align*}
	\widehat{\mu}_{\mathrm{MoM},k} \coloneqq \mathrm{median} \biggl\{ \frac{1}{|B_v|} \sum_{i \in B_v} X_{i,k} : v \in [V] \biggr\},
\end{align*}
where $X_{i,k}$ denotes the $k$-th component of $\bX_i \in \mathbb{R}^d$. Unlike the empirical mean, the MoM estimator achieves sub-Gaussian concentration under only finite second moments and mitigates the influence of extreme values. Building on this property, we propose a robust DA argmin test that achieves the minimax separation rate under finite variance assumptions. 

Let $\mathcal{P}^{\leq 2}$ denote the class of distributions on $\mathbb{R}^d$ whose marginal variances are uniformly bounded by $\sigma^2$, i.e., $\sup_{P \in \mathcal{P}^{\leq 2}} \max_{k \in [d]} \mathrm{Var}_P(X_k) \leq \sigma^2$, where $X_k$ is the $k$-th component of $\bX \sim P$. In particular, every $\sigma^2$-sub-Gaussian distribution belongs to $\mathcal{P}^{\leq 2}$. We then define the alternative hypothesis class as
\begin{align*}
	\mathcal{P}_{1,r}^{\leq 2}(\varepsilon;\tau) \coloneqq \bigl\{ P \in \mathcal{P}^{\leq 2}: \mu_r - \mu_{\star} \geq \varepsilon \ \text{and} \ |\mathbb{C}_r| = \tau \bigr\}.
\end{align*}
We first define the plug-in estimator $\widetilde{s}_{\mathrm{plug}}$ by replacing the sample means in $\widehat{s}_{\mathrm{plug}}$ with the MoM estimates:
\begin{align*}
	\widetilde{s}_{\mathrm{plug}} \coloneqq \sargmin_{k \in [d]\setminus \{r\}} \widehat{\mu}_{\mathrm{MoM},k}.
\end{align*}
Similarly, we define the noise-adjusted estimator $\widetilde{s}_{\mathrm{adj}}$ based on a noise-adjusted difference of MoM estimates:
\begin{align*}
	\widetilde{s}_{\mathrm{adj}} \coloneqq \sargmin_{k \in [d] \setminus \{r\}} \frac{\widehat{\mu}_{\mathrm{MoM},k} - \widehat{\mu}_{\mathrm{MoM},r}}{\sqrt{\bm{\gamma}_k^\top \widehat{\bm{\Sigma}}^{(2)} \bm{\gamma}_k} \vee \kappa},
\end{align*}
where $\kappa > 0$ is a small constant considered in $\widehat{s}_{\mathrm{adj}}$. We refer to the DA argmin test using either $\widetilde{s}_{\mathrm{plug}}$ or $\widetilde{s}_{\mathrm{adj}}$ as the \emph{robust DA argmin test}. Since the validity result in \Cref{The: DA validity} holds for any random variable $\widehat{s} \in [d]\!\setminus\!\{r\}$ independent of the first half of the sample, the robust variant inherits the same asymptotic validity guarantees as the original test. We now examine the asymptotic power of the robust DA argmin test under heavy-tailed distributions, which holds for both $\widetilde{s} = \widetilde{s}_{\mathrm{plug}}$ and $\widetilde{s} =\widetilde{s}_{\mathrm{adj}}$.

\begin{theorem} \label{Thm: robust}
	For any $\tau \in \{0,1,\dots,d-2\}$, suppose that $\varepsilon \geq C_n' \varepsilon^\star$ where $C_n'$ is any positive sequence diverging to infinity as $n \to \infty$ and $\varepsilon^\star$ was defined in~\eqref{Eq: critical radius}. Set $\eta = 1/2 \wedge (C_n'^{-1} \vee e^{-C_n} \vee e^{-n/18})$. Then the asymptotic uniform power of the robust DA argmin test with $V = 4.5\ceil{\log(1/\eta)}$ over $\mathcal{P}_{1,r}^{\leq 2}(\varepsilon;\tau)$ equals one:
	\begin{align*}
		\lim_{n \to \infty} \inf_{P \in \mathcal{P}_{1,r}^{\leq 2}(\varepsilon;\tau)} P \Big( \sqrt{n}\bm{\gamma}_{\widetilde{s}}^\top \overline{\bX}^{(1)} > z_{1-\alpha} \sqrt{\bm{\gamma}_{\widetilde{s}}^\top \widehat{\bm{\Sigma}}^{(1)} \bm{\gamma}_{\widetilde{s}}} \Big) = 1.
	\end{align*}
\end{theorem}
The above theorem establishes that the robust DA argmin test achieves the same minimax separation rate $\varepsilon^\star$  under heavy-tailed distributions with finite variance. The proof follows that of \Cref{Thm: upper bound} almost verbatim, with only minor variations outlined in \Cref{subsec:proof-robust}. While \Cref{Thm: robust} represents a clear improvement over \Cref{Thm: upper bound}, the MoM-based approach comes with several practical drawbacks. Most notably, its optimal performance depends on the choice of the partition parameter $\eta$, which itself depends on the sequences $C_n$ and $C_n'$. This limitation stems from the inherent dependence of the MoM estimator on the user-specified confidence level.

Although we focus on the MoM estimator as a concrete example, it is important to note that the proof of Theorem~\ref{Thm: robust} is more broadly applicable. In particular, the same power guarantee can be established for any robust estimator that exhibits sub-Gaussian tails under finite second moment conditions—such as Catoni's M-estimator~\citep{catoni2012challenging} and the truncated empirical mean~\citep{bubeck2013bandits}, with only minor changes to the proof in order to incorporate minor differences between the formal guarantees of these estimators. Moreover, the corresponding robust DA argmin confidence set can be constructed by inverting the robust DA argmin test.

Preliminary numerical results in \Cref{Sec: additional-robust} indicate that the MoM variant does not consistently outperform the original DA argmin test in heavy-tailed settings, possibly due to the loss of efficiency induced by data splitting. We also observe that the robust DA argmin test based on Catoni's M-estimator exhibits a similar performance to the original DA argmin test against heavy-tailed alternatives. These findings suggest that while current robust approaches offer theoretical advantages, developing a practically effective and robust alternative remains an important challenge for future research.

\section{Dimension-agnostic model confidence sets} \label{Sec: MCS}

While not the main point of our paper, we point out that our techniques can be used to derive confidence sets for the argmin that have uniform coverage guarantees (recall the definition of uniform coverage in Section~\ref{sec:intro}). We believe that what we present below is the first nontrivial approach to guarantee uniform coverage in \emph{high-dimensional} settings.

\paragraph{A revisit to the MCS.} Before introducing our approach, we first briefly review the model confidence set (MCS) approach of \citet{hansen2011model}, which is one of the most influential approaches in the literature for achieving uniform coverage. The MCS algorithm \citep[][Definition 2]{hansen2011model} is an iterative procedure that starts with a set of candidate models $\mathcal{M}^0$ and sequentially tests whether subsets $\mathcal{M} \subseteq \mathcal{M}_0$ are optimal. When the test fails to reject the null, the procedure accepts the subset as optimal; otherwise, it eliminates an worst-performing model in the subset. They show that the resulting confidence set has asymptotic uniform coverage under some conditions. However, as we elaborate in \Cref{App: MCS validity},\footnote{We thank Jing Lei for helpful discussions and sharing an explicit proof of \citet[][Theorem 1]{hansen2011model}.} the asymptotic validity of their procedure implicitly assumes that the size of the initial model set $\mathcal{M}^0$ is fixed. Without additional stronger assumptions, we believe that the MCS procedure should be viewed as a fixed-size model selection method. In contrast, the DA-MCS procedures we introduce below remain valid even when the number of candidate models grows with the sample size.

\paragraph{One-step construction of a DA-MCS with uniform coverage.}

We first present a very simple approach that guarantees uniform coverage. We simply run the DA argmin test on the full sample at level $\alpha/d$ (instead of $\alpha$). The uniform coverage guarantee is ensured by the union bound, and it is the first direct method we are aware of with valid uniform coverage in high-dimensional settings. Further, this method is still globally minimax optimal. However, this method is not locally minimax optimal in that it does not adapt to the cardinality of the confusion set, so we propose the following two-step construction. 

\paragraph{Two-step construction of a DA-MCS with uniform coverage.}
To provide a better benchmark with uniform coverage, we now introduce a modified version of our own confidence set construction that attains a \emph{uniform} coverage guarantee. Let $\psi_k(S,c)$ denote the application of our DA argmin test for $H_0: k \in \Theta$ to the dataset $S$ at level $c$.
\begin{enumerate}
\item[\textbf{1.}]\textbf{Pre‑screening.}
For each $k \in [d]$, apply the DA argmin test $\psi_k$ to the \emph{second} half of the data $D_2 \coloneqq \{\bX_{n+1},\ldots,\bX_{2n}\}$ at a nominal level
  $n^{-1/2}$, and define  the pre-screened set as
  \[
    \widehat\Theta^{(2)} \coloneqq
    \bigl\{k\in[d] : \text{the null for }k\text{ is \emph{not} rejected by $\psi_k(D_2,n^{-1/2})$}\bigr\}.
  \]
\item[\textbf{2.}]\textbf{Final inference.}
  For each $k \in [d]$, run the DA argmin test $\psi_k$ on the \emph{full} sample $D \coloneqq \{\bX_1,\ldots,\bX_{2n}\}$ at level
  $\alpha' := \alpha/(1 \vee |\widehat\Theta^{(2)}|)$. The final confidence set is given as
  \[
    \widehat\Theta_{\mathrm{DA}}^{\mathrm{uni}} \coloneqq
	\bigl\{k\in[d] : \text{the null for }k\text{ is \emph{not} rejected by $\psi_k (D, \alpha')$}\bigr\}.
  \]
\end{enumerate}
In essence, we apply the (pointwise) DA argmin test at the adjusted level $\alpha/(1 \vee |\widehat\Theta^{(2)}|)$ to guarantee uniform coverage where $|\widehat{\Theta}^{(2)}|$ serves as a data-driven proxy for the cardinality of the true argmin set $|\Theta|$.

Let $\mathcal{P}^{\le 3}$ be the collection of distributions $P$ satisfying the third moment condition
\[
\max_{k\in[d]\setminus\{r\}}
\mathbb{E}_{P}\!\left[\frac{|W_k|^{3}}
      {\bigl\{\mathbb{E}_{P}[W_k^{2}]\bigr\}^{3/2}}\right]
      \le C
\quad\text{for every }r\in\Theta(P),
\]
where $W_k$ is defined as in \eqref{Eq: definition of Wk} and $C>0$ is a universal constant. This third‑moment bound is slightly stronger than the truncated
second‑moment requirement in~\eqref{Eq: moment condition} as discussed in the main text. Additionally, we will assume that  
\[
\sup_{P\in\mathcal{P}^{\le 3}}\!|\Theta(P)| = o\bigl(n^{1/2}\bigr),
\]
which ensures that the remainder term in the Berry--Esseen bound is uniformly negligible over all $r \in \Theta$. Notably, this assumption imposes no restriction on the ambient dimension $d$, but only on the cardinality of the argmin set. Under these conditions, the following proposition shows that the above two-step construction yields a confidence set with uniform coverage.
\begin{theorem} \label{Thm: uniform coverage}
	As long as $\sup_{P\in\mathcal{P}^{\le 3}}\!|\Theta(P)| = o\bigl(n^{1/2}\bigr)$, the confidence set $\widehat\Theta_{\mathrm{DA}}^{\mathrm{uni}}$ from the two-step construction satisfies
	\begin{align*}
		\liminf_{n \to \infty}\inf_{P \in \mathcal{P}^{\leq 3} }P(\Theta \subseteq \widehat{\Theta}_{\mathrm{DA}}^{\mathrm{uni}}) \geq 1 - \alpha.
	\end{align*}
\end{theorem}
The proof of \Cref{Thm: uniform coverage} is provided in \Cref{App: uniform coverage}. Simulation results in \Cref{Sec: Sim MCS,Sec: additional-simulation-mcs} compare the empirical performance of the one‑step and two‑step DA‑MCS procedures and demonstrate that the two‑step variant achieves uniform coverage much closer to the nominal level \(1-\alpha\).

\paragraph{Confidence set for the smallest mean.} A closely related task to constructing a confidence set for the argmin set is developing a confidence set for the smallest mean $\mu_\star = \min_{k \in [d]} \mu_k$. Let $\overline{X}_k$ and $\widehat{\sigma}_k$ denote the sample mean and sample standard deviation for the $k$-th population, respectively. A natural starting point is to determine a critical threshold $t_{1-\alpha}$ satisfying
\begin{align*}
	P\biggl( \max_{k \in [d]} \frac{\sqrt{2n}|\overline{X}_k - \mu_k|}{\widehat{\sigma}_k} \leq t_{1-\alpha} \biggr) \geq 1 - \alpha + o(1),
\end{align*}
which in turn yields an asymptotically valid confidence set for $\mu_{\star}$ as 
\begin{align*} 
\mathcal{C}_1 =	\left[ \min_{k \in [d]} \biggl\{ \overline{X}_k - t_{1-\alpha} \frac{\widehat{\sigma}_k}{\sqrt{2n}} \biggr\}, ~ \min_{k \in [d]} \biggl\{\overline{X}_k + t_{1-\alpha} \frac{\widehat{\sigma}_k}{\sqrt{2n}} \biggr\} \right].
\end{align*}
The simplest choice $t_{1-\alpha}=z_{1-\frac{\alpha}{2d}}$ follows from a Bonferroni correction, but this strategy is notoriously conservative in high‑dimensional settings. A refined approach works as follows: 
\begin{enumerate}
    \item First run the DA-MCS procedure on $D_2$ at a nominal level $\gamma_n$ tending to zero, thereby obtaining a uniformly valid screening set $\widehat{\Theta}_{\mathrm{DA}}^{\mathrm{uni}}$.
    \item Then apply the above construction together with the Bonferroni correction to $D_1$ but only for the indices in this data‑adaptive subset $\widehat{\Theta}_{\mathrm{DA}}^{\mathrm{uni}}$ of size $\widehat{d} \coloneqq |\widehat{\Theta}_{\mathrm{DA}}^{\mathrm{uni}}|$, to obtain the interval
\begin{align*} 
	\mathcal{C}_2 = \left[ \min_{k \in \widehat{\Theta}_{\mathrm{DA}}^{\mathrm{uni}}} \biggl\{ \overline{X}_k^{(1)} - z_{1-\frac{\alpha}{2\widehat{d}}} \frac{\widehat{\sigma}_k^{(1)}}{\sqrt{n}} \biggr\}   , ~ \min_{k \in \widehat{\Theta}_{\mathrm{DA}}^{\mathrm{uni}}} \biggl\{\overline{X}_k^{(1)} + z_{1-\frac{\alpha}{2\widehat{d}}} \frac{\widehat{\sigma}_k^{(1)}}{\sqrt{n}} \biggr\}\right],
\end{align*}
	where $\overline{X}_k^{(1)}$ and $\widehat{\sigma}_k^{(1)}$ are the sample mean and sample standard deviation for the $k$-th population based on $D_1$.
\end{enumerate}
This data-adaptive approach has the following uniform coverage guarantee. 
\begin{theorem} \label{Thm: CS for the smallest mean}
If $\sup_{P\in\mathcal{P}^{\le 3}}\!|\Theta(P)| = o\bigl(n^{1/2}\bigr)$ and $\gamma_n \to 0$, then $\inf_{P\in\mathcal{P}^{\le 3}}P(\Theta \subseteq \widehat{\Theta}_{\mathrm{DA}}^{\mathrm{uni}}) = 1 - o(1)$. If we further have $\sup_{P \in \mathcal{P}^{\leq 3}}\mE_P[\widehat{d}] = o(n^{1/2})$, then 
 \begin{align*}
		\liminf_{n \to \infty}\inf_{P \in \mathcal{P}^{\leq 3} }P(\mu_\star \in \mathcal{C}_2) \geq 1 - \alpha.
\end{align*}
\end{theorem}

When compared with $\mathcal{C}_1$, the confidence set $\mathcal{C}_2$ involves a less favorable factor of $1/\sqrt{n}$ than $1/\sqrt{2n}$, which is a consequence of sample splitting. Despite this efficiency loss by a factor of $\sqrt{2}$, $\mathcal{C}_2$ is nevertheless expected to yield narrower intervals than the confidence set $\mathcal{C}_1$ whenever $\widehat d \ll d$. See \Cref{Sec: simulation-smallest-mean} for supporting numerical evidence. While one could improve efficiency further via repeated sample splitting and aggregation, we do not pursue this extension here as they fall outside the main focus of our paper.

\section{Simulations} \label{Sec: simulation}
In this section, we conduct simulation studies to evaluate the finite-sample performance of the DA argmin test and other existing methods under the setting $r=1$. Specifically, we compare the following methods in terms of size and power:
\begin{itemize}
    \item\texttt{LOO}: The method proposed by \citet{zhang2024winners}, using the data-driven parameter selection procedure recommended by the authors.
    \item \texttt{Bonferroni}: The one-sided $t$-test with Bonferroni correction. Specifically, it performs one-sided $t$-tests for $H_0: \mu_1 \leq \mu_k$ versus $H_1: \mu_1 > \mu_k$ for each $k \in \{2, 3, \ldots, d\}$ at the adjusted level $\alpha/(d-1)$, and rejects the null if any of the tests is significant.
    \item \texttt{csranks}: The method based on rank confidence intervals by \citet{mogstad2024inference}. It constructs confidence intervals for ranks by approximating the distribution of the maximum of pairwise mean differences via a bootstrap method. The null hypothesis is then rejected whenever the (marginal) lower bound of the confidence interval for the rank of the first population includes rank one. The procedure is implemented using the \texttt{csranks} package available on CRAN.
    \item \texttt{MCS}: The method introduced by \citet{hansen2011model}, implemented via the \texttt{MCS} package on CRAN. We adopt the default settings provided by the package, except that the number of bootstrap replications is reduced from $B = 5,\!000$ to $B = 100$ to mitigate computational overhead. 
    \item \texttt{DA-plug}: Our proposed DA argmin test using the plug-in selection method $\widehat{s} = \widehat{s}_{\mathrm{plug}}$.
    \item \texttt{DA-plug}$^{\times 10}$: This variant averages the \texttt{DA-plug} test statistics over 10 random data splits. The null is rejected if the averaged statistic exceeds a threshold determined via the subsampling method of \citet{guo2025rank}.
    \item \texttt{DA-adj}: Our proposed DA argmin test using the noise-adjusted selection method $\widehat{s} = \widehat{s}_{\mathrm{adj}}$.
    \item \texttt{DA-adj}$^{\times 10}$: This variant averages the \texttt{DA-adj} test statistics over 10 random data splits. The null is rejected if the averaged statistic exceeds a threshold determined via the subsampling method of \citet{guo2025rank}.
\end{itemize}
We examine the type I error rates of these methods across various significance levels $\alpha$ in \Cref{Sec: type1}, and investigate their power and validity under homoskedasticity in \Cref{Sec: power} and under heteroskedasticity in \Cref{Sec: power-unequal}. Additional empirical results in high-dimensional settings and on real-world data are presented in \Cref{Sec: highdim} and in \Cref{Sec: real data} anc \Cref{Sec: real-data-argmax}, respectively. Experiments on uniform coverage are in \Cref{Sec: Sim MCS} and \Cref{Sec: additional-simulation-mcs} and on heavy-tailed settings in \Cref{Sec: additional-robust}. We also refer to \Cref{tab:execution_time_summary} for a summary of execution times of all methods.

\begin{table}[htbp]
	\begin{minipage}{0.58\textwidth}
	\textbf{Computational efficiency.} \Cref{tab:execution_time_summary} summarizes the execution times (in seconds) for each method evaluated in our simulations, with dimension $d = 100$ and total sample size $2n = 1000$. All procedures were implemented in \texttt{R} and executed on a single core. Among them, the \texttt{MCS}$^{\times 5000}$ method—using the default setting of $B = 5,\!000$ bootstrap replications—is by far the most computationally demanding, followed by the reduced version \texttt{MCS}$^{\times 100}$ with $B = 100$. In contrast, the proposed DA methods are highly efficient, with the base versions completing in under 0.01 seconds. While the aggregated variants (\texttt{DA-plug}$^{\times 10}$ and \texttt{DA-adj}$^{\times 10}$) incur additional computational cost due to repeated data splits, they remain practical for moderate-scale applications.
	\end{minipage}%
	\hfill
	\begin{minipage}{0.4\textwidth}
	\centering
	\small
	\captionof{table}{Elapsed Time in Seconds}
	\begin{tabular}{l r}
	  \toprule
	  \textbf{Method} & \textbf{Elapsed Time} \\
	  \midrule
	  \texttt{LOO}         & 0.090 \\
	  \texttt{Bonferroni}  & 0.008 \\
	  \texttt{MCS}$^{\times 5000}$         & 633.124 \\
	  \texttt{MCS}$^{\times 100}$         & 28.86 \\
	  \texttt{csranks}     & 0.012 \\
	  \texttt{DA-plug}     & 0.007 \\
	  \texttt{DA-adj}      & 0.010 \\
	  \texttt{DA-plug}$^{\times 10}$   & 5.141 \\
	  \texttt{DA-adj}$^{\times 10}$    & 10.806 \\
	  \bottomrule
	\end{tabular}
	\label{tab:execution_time_summary}
	\end{minipage}
	\end{table}

\subsection{Type I error rate across nominal levels} \label{Sec: type1}
\citet{zhang2024winners} establish asymptotic validity of \texttt{LOO} under relatively strong conditions, including bounded random variables and a lower bound on the smallest eigenvalue of the covariance matrix. Although these assumptions might be relaxed through more refined theoretical developments, it remains unclear whether the practical implementation of \texttt{LOO}, especially when data-driven tuning is used, ensures reliable type I error control in finite samples. In this subsection, we examine this aspect through an empirical investigation along with the empirical size of the other methods.

To this end, we evaluate the empirical type I error rates of \texttt{LOO}, \texttt{Bonferroni}, \texttt{csranks}, \texttt{DA-plug}, and \texttt{DA-adj} under a simple yet informative setting. Specifically, we consider $\bX \sim N(\bm{\mu}, \bm{I}_d)$, with $\bm{\mu} = (0,0,0,0)^\top$ for $d = 4$ and $\bm{\mu} = (0,0,0,0,10,\ldots,10)^\top$ for $d = 100$, and generate $2n \in \{500, 2000, 5000\}$ samples. We compute the empirical rejection rates across various significance levels $\alpha \in \{0.01, 0.05, \ldots, 0.45, 0.50\}$.

\begin{figure}[h]
	\centering
	\includegraphics[width=1\textwidth]{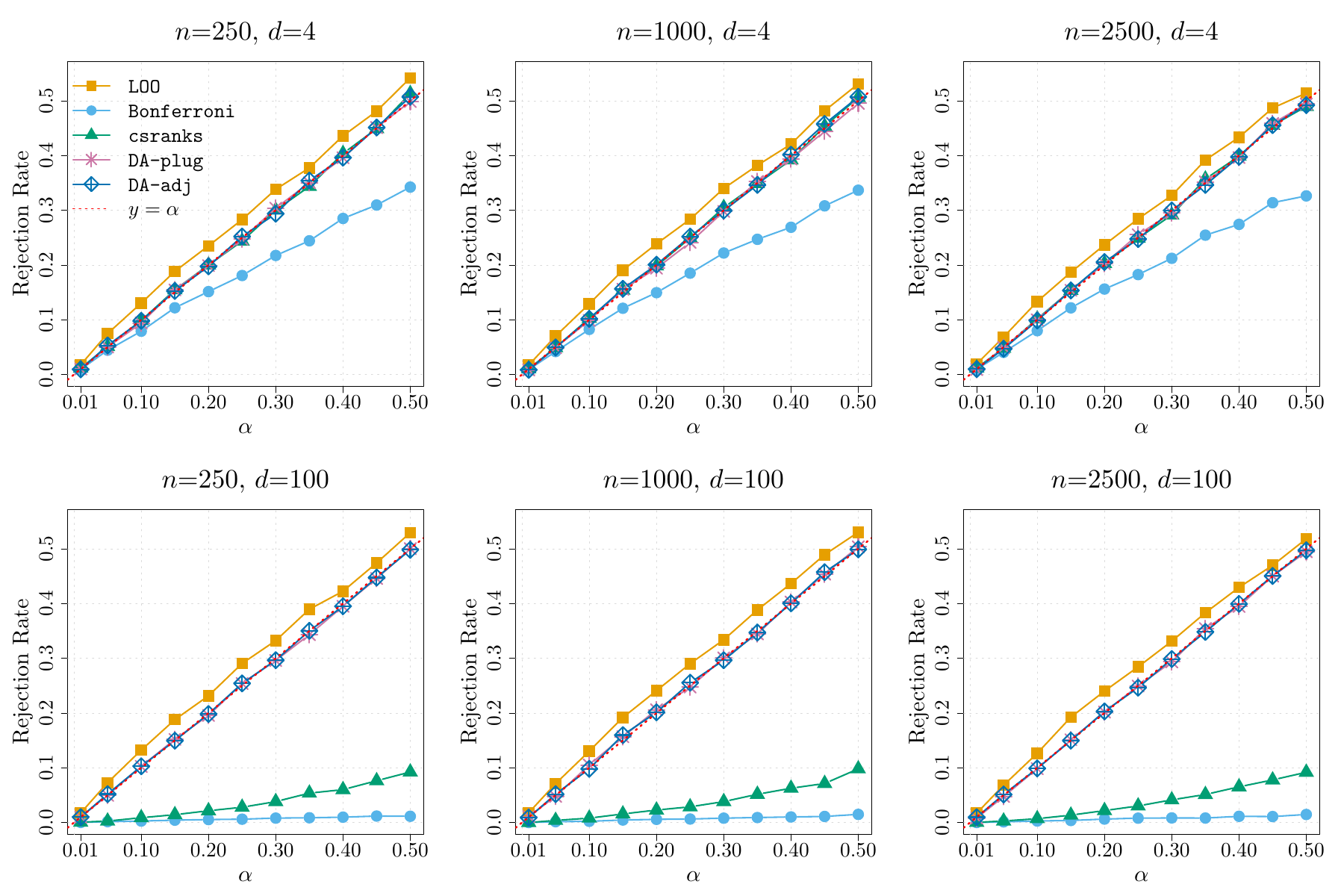}
	\caption{Empirical type I error rates for \texttt{LOO}, \texttt{Bonferroni}, \texttt{csranks}, \texttt{DA-plug}, and \texttt{DA-adj} are presented across various sample sizes and dimensions. The results consistently indicate that \texttt{LOO} tends to be liberal in controlling the type I error rate, even as the sample size increases, whereas \texttt{Bonferroni} generally exhibits conservative behavior. The \texttt{csranks} method performs well when $d=4$ but becomes increasingly conservative at $d=100$. In contrast, both \texttt{DA-plug} and \texttt{DA-adj} reliably maintain the nominal error level across different significance levels $\alpha$ and combinations of $n$ and $d$.}
	\label{fig:type1}
\end{figure}

The results, summarized in \Cref{fig:type1}, are based on $10,\!000$ repetitions. As shown in the figure, the \texttt{LOO} method tends to be liberal in its type I error, and the gap between the empirical and nominal levels (ranging from $0$ to $0.05$) does not diminish as the sample size increases. This observation suggests that the violation—albeit relatively mild—is not merely a finite-sample artifact. While our empirical settings are limited, these findings underscore that the theoretical guarantees of \texttt{LOO} may not fully translate into reliable practical performance, particularly regarding type I error control. In contrast, both \texttt{DA-plug} and \texttt{DA-adj} consistently exhibit accurate type I error control across all considered settings. The \texttt{Bonferroni} method tends to be conservative, with its conservativeness becoming more pronounced in higher dimensions. The \texttt{csranks} method, on the other hand, performs reliably when $d=4$ but becomes increasingly conservative when $d=100$. Consequently, these results support the use of more stable alternatives such as \texttt{DA-plug} and \texttt{DA-adj} in applications where rigorous and tight control of the type I error is essential. The \texttt{DA-plug}$^{\times 10}$ and \texttt{DA-adj}$^{\times 10}$ methods as well as \texttt{MCS} are excluded from this analysis due to their computational demands. Their performance is evaluated separately in the subsequent sections.

\subsection{Power and validity under homoskedasticity} \label{Sec: power}
We next explore the empirical power and size of the considered tests under various signal structures and homoskedastic covariance settings. Specifically, we consider a simulation setup where $\bX \sim N(\bm{\mu},\bm{\Sigma})$ with $2n = 1,\!000$ and $d = 100$. Each scenario is repeated $5,\!000$ times to approximate the power and the size except for \texttt{MCS}, which is repeated $500$ times due to its higher computational demands (see \Cref{tab:execution_time_summary}). To represent different signal structures, we examine three distinct mean vectors under the alternative:
\begin{enumerate}
	\item[(i)] $\bm{\mu}^{(a)} = (0.1, 0, 0.1, 0.1, \ldots, 0.1)^\top$, representing small values for  non-informative components;
	\item[(ii)] $\bm{\mu}^{(b)} = (\mu^{(b)}_1, \ldots, \mu^{(b)}_d)^\top$, where $\mu^{(b)}_1 = 0.2$ and $\mu^{(b)}_k = 0.1 + \frac{k - 2}{d - 2} \times 0.9$ for $k \in \{2, \ldots, d\}$, representing gradually increasing values for the non-informative components;
	\item[(iii)] $\bm{\mu}^{(c)} = (0.05, 0, 0, 0, 10, 10, \ldots, 10)^\top$, representing large values for  non-informative components.
\end{enumerate}
The covariance structure of the features follows a Toeplitz form where the covariance matrix is defined as $\Sigma_{k_1k_2} = \rho^{|k_1 - k_2|}$ for $k_1,k_2 \in [d]$, with $\rho \in \{0, 0.4, 0.8\}$ representing no correlation, moderate correlation, and strong correlation, respectively. 

To assess the empirical size, we construct the mean vectors $\bm{\mu}^{(a,0)}$, $\bm{\mu}^{(b,0)}$, and $\bm{\mu}^{(c,0)}$ by replacing the first component of $\bm{\mu}^{(a)}$, $\bm{\mu}^{(b)}$, and $\bm{\mu}^{(c)}$ with their respective minimum values, while keeping the remaining components and the covariance structure unchanged. This modification ensures that the null hypothesis is satisfied.

\begin{table}[h]
	\centering
	
	\begin{minipage}{\textwidth}
		\small	
		\centering
		\caption{Empirical power at the significance level $\alpha = 0.05$ for different mean structures and correlation levels under equal variance. The highest power in each scenario is highlighted in bold, and deeper color intensity indicates higher power. Our DA methods are the most powerful, except in the final case where LOO dominates it, but Table~\ref{tab:type1-equal} shows that LOO does not control type I error in this setting. Among methods that do control type I error, DA methods are the most powerful across the board.} \label{tab:rejection-rates}
		\renewcommand{\arraystretch}{1.0}
		\setlength{\tabcolsep}{5pt}
		\begin{tabular}{lccccccccc}
		\toprule
		\multirow{2}{*}{\textbf{Method}} & \multicolumn{3}{c}{$\bm{\mu}^{(a)}$ + equal variance} & \multicolumn{3}{c}{$\bm{\mu}^{(b)}$ + equal variance} & \multicolumn{3}{c}{$\bm{\mu}^{(c)}$ + equal variance} \\
		& $\rho=0$ & $\rho=0.4$ & $\rho=0.8$ & $\rho=0$ & $\rho=0.4$ & $\rho=0.8$ & $\rho=0$ & $\rho=0.4$ & $\rho=0.8$ \\
		\cmidrule(lr){2-4} \cmidrule(lr){5-7} \cmidrule(lr){8-10}
		\texttt{LOO} & \cellcolor{vlightorange}0.098 & \cellcolor{lightorange}0.157 & \cellcolor{lightorange}0.181
					& \cellcolor{midorange}0.430 & \cellcolor{midorange}0.437 & \cellcolor{strongorange}0.686
					& \cellcolor{midorange}\textbf{0.427} & \cellcolor{midorange}\textbf{0.480} & \cellcolor{deeporange}\textbf{0.786} \\
		\texttt{Bonferroni} & \cellcolor{lightorange}0.154 & \cellcolor{vlightorange}0.106 & \cellcolor{vlightorange}0.025
					& \cellcolor{lightorange}0.285 & \cellcolor{lightorange}0.205 & \cellcolor{vlightorange}0.051
					& \cellcolor{vlightorange}0.040 & \cellcolor{vlightorange}0.025 & \cellcolor{vlightorange}0.001 \\
		\texttt{csranks} & \cellcolor{lightorange}0.231 & \cellcolor{midorange}0.456 & \cellcolor{deeporange}0.980
					& \cellcolor{lightorange}0.363 & \cellcolor{midorange}0.543 & \cellcolor{deeporange}0.980
					& \cellcolor{vlightorange}0.073 & \cellcolor{vlightorange}0.099 & \cellcolor{midorange}0.380 \\
		\texttt{MCS} & \cellcolor{vlightorange}0.000 & \cellcolor{vlightorange}0.004 & \cellcolor{vlightorange}0.008
					& \cellcolor{vlightorange}0.048 & \cellcolor{vlightorange}0.054 & \cellcolor{lightorange}0.140
					& \cellcolor{midorange}0.352 & \cellcolor{midorange}0.354 & \cellcolor{strongorange}0.646 \\
		\texttt{DA-plug} & \cellcolor{lightorange}0.219 & \cellcolor{lightorange}0.305 & \cellcolor{midorange}0.501
					& \cellcolor{lightorange}0.371 & \cellcolor{midorange}0.424 & \cellcolor{strongorange}0.679
					& \cellcolor{lightorange}0.205 & \cellcolor{lightorange}0.238 & \cellcolor{midorange}0.426 \\
		\texttt{DA-plug}$^{\times 10}$ & \cellcolor{midorange}\textbf{0.307} & \cellcolor{midorange}0.401 & \cellcolor{strongorange}0.727
					& \cellcolor{strongorange}\textbf{0.593} & \cellcolor{strongorange}0.674 & \cellcolor{deeporange}0.957
					& \cellcolor{midorange}0.310 & \cellcolor{midorange}0.359 & \cellcolor{strongorange}0.655 \\
		\texttt{DA-adj} & \cellcolor{lightorange}0.232 & \cellcolor{midorange}0.448 & \cellcolor{deeporange}0.931
					& \cellcolor{lightorange}0.365 & \cellcolor{midorange}0.506 & \cellcolor{deeporange}0.932
					& \cellcolor{lightorange}0.207 & \cellcolor{lightorange}0.250 & \cellcolor{midorange}0.477 \\
		\texttt{DA-adj}$^{\times 10}$ & \cellcolor{midorange}\textbf{0.307} & \cellcolor{strongorange}\textbf{0.589} & \cellcolor{deeporange}\textbf{0.988}
					& \cellcolor{strongorange}0.585 & \cellcolor{strongorange}\textbf{0.728} & \cellcolor{deeporange}\textbf{0.994}
					& \cellcolor{midorange}0.300 & \cellcolor{midorange}0.370 & \cellcolor{strongorange}0.697 \\
		\bottomrule
		\end{tabular}
	\end{minipage}
	
	\vspace{1.5em}
	
	\begin{minipage}{\textwidth}
		\small
		\centering
		\caption{Empirical type I error at the significance level $\alpha = 0.05$ for different mean structures and correlation levels under equal variance. Blue shading indicates over-rejection (liberal tests), green indicates under-rejection (conservative tests), and white indicates appropriate rejection rates (correct coverage). Our DA methods maintain the right coverage throughout; others are either too conservative or anti-conservative.}
		\renewcommand{\arraystretch}{1.0}
		\setlength{\tabcolsep}{5pt}
		\begin{tabular}{lccc ccc ccc}
		\toprule
		\multirow{2}{*}{\textbf{Method}} 
		& \multicolumn{3}{c}{$\bm{\mu}^{(a,0)}$ + equal variance} 
		& \multicolumn{3}{c}{$\bm{\mu}^{(b,0)}$ + equal variance} 
		& \multicolumn{3}{c}{$\bm{\mu}^{(c,0)}$ + equal variance} \\
		& $\rho=0$ & $\rho=0.4$ & $\rho=0.8$ 
		& $\rho=0$ & $\rho=0.4$ & $\rho=0.8$ 
		& $\rho=0$ & $\rho=0.4$ & $\rho=0.8$ \\
		\cmidrule(lr){2-4} \cmidrule(lr){5-7} \cmidrule(lr){8-10}
		\texttt{LOO} & \cellcolor{stronggreen}0.011 & \cellcolor{deepgreen}0.006 & \cellcolor{deepgreen}0.000 & \cellcolor{stronggreen}0.014 & \cellcolor{stronggreen}0.012 & \cellcolor{deepgreen}0.007 & \cellcolor{lightblue}0.071 & \cellcolor{lightblue}0.073 & \cellcolor{lightblue}0.067 \\
		\texttt{Bonferroni} & \cellcolor{deepgreen}0.001 & \cellcolor{deepgreen}0.000 & \cellcolor{deepgreen}0.000 & \cellcolor{deepgreen}0.001 & \cellcolor{deepgreen}0.000 & \cellcolor{deepgreen}0.000 & \cellcolor{deepgreen}0.001 & \cellcolor{deepgreen}0.001 & \cellcolor{deepgreen}0.000 \\
		\texttt{csranks} & \cellcolor{deepgreen}0.003 & \cellcolor{deepgreen}0.001 & \cellcolor{deepgreen}0.002 & \cellcolor{deepgreen}0.001 & \cellcolor{deepgreen}0.002 & \cellcolor{deepgreen}0.003 & \cellcolor{deepgreen}0.004 & \cellcolor{deepgreen}0.003 & \cellcolor{deepgreen}0.006 \\
		\texttt{MCS} & \cellcolor{deepgreen}0.000 & \cellcolor{deepgreen}0.000 & \cellcolor{deepgreen}0.000 & \cellcolor{deepgreen}0.000 & \cellcolor{deepgreen}0.000 & \cellcolor{deepgreen}0.000 & \cellcolor{lightgreen}0.030 & \cellcolor{midgreen}0.024 & 0.040 \\
		\texttt{DA-plug} & \cellcolor{stronggreen}0.019 & \cellcolor{midgreen}0.021 & \cellcolor{midgreen}0.024 & \cellcolor{lightgreen}0.033 & \cellcolor{midgreen}0.026 & \cellcolor{midgreen}0.027 & 0.053 & 0.049 & 0.047 \\
		\texttt{DA-plug}$^{\times 10}$ & \cellcolor{midgreen}0.023 & \cellcolor{stronggreen}0.019 & \cellcolor{midgreen}0.025 & \cellcolor{lightgreen}0.030 & \cellcolor{lightgreen}0.031 & \cellcolor{midgreen}0.028 & 0.053 & 0.056 & 0.051 \\
		\texttt{DA-adj} & \cellcolor{midgreen}0.021 & \cellcolor{midgreen}0.020 & \cellcolor{lightgreen}0.030 & \cellcolor{midgreen}0.028 & \cellcolor{midgreen}0.029 & \cellcolor{midgreen}0.028 & 0.051 & 0.049 & 0.046 \\
		\texttt{DA-adj}$^{\times 10}$ & \cellcolor{midgreen}0.025 & \cellcolor{midgreen}0.023 & \cellcolor{lightgreen}0.034 & \cellcolor{lightgreen}0.032 & \cellcolor{lightgreen}0.031 & \cellcolor{midgreen}0.029 & 0.051 & 0.054 & 0.052 \\

		\bottomrule
		\end{tabular}
		\label{tab:type1-equal}
		\end{minipage}
\end{table}



The simulation results in \Cref{tab:rejection-rates,tab:type1-equal} summarize the empirical power and size of the considered methods under homoskedastic settings. When comparing the proposed methods, the two DA tests (\texttt{DA-plug} and \texttt{DA-adj}) exhibit similar power when $\rho=0$. In all other scenarios, however, \texttt{DA-adj} consistently outperforms \texttt{DA-plug}, by accounting more effectively for the correlation structure. The aggregated versions (\texttt{DA-plug}$^{\times 10}$ and \texttt{DA-adj}$^{\times 10}$) further improve power, albeit at the cost of increased computation. Among all methods, \texttt{DA-adj}$^{\times 10}$ generally demonstrates strong power across most settings and frequently achieves the highest power. One notable exception is the scenario with $\bm{\mu}^{(c)}$, where \texttt{LOO} shows slightly higher power. However, in this case, \texttt{LOO} also exhibits inflated type I error rates, as demonstrated in \Cref{tab:type1-equal}, which may compromise the validity of the power comparison. The \texttt{Bonferroni} procedure shows limited power in most settings due to its conservative nature. While \texttt{csranks} performs reasonably well under strong correlation $(\rho = 0.8)$, it generally yields lower power than \texttt{DA-adj}$^{\times 10}$ in other cases. The \texttt{MCS} method has limited power in the first two scenarios, whereas it performs well in the last scenario with $\bm{\mu}^{(c)}$. Although no single method dominates across all scenarios, the proposed DA-argmin test—particularly with noise-adjusted selection and aggregation—consistently demonstrates strong and robust power while maintaining correct size control across a range of signal structures and correlation levels.

\subsection{Power and validity under heteroskedasticity} \label{Sec: power-unequal}
To assess robustness under heteroskedasticity, we also consider an unequal variance setting, where the diagonal elements of the covariance matrix are modified such that $\Sigma_{kk} = 20$ for $k \in \{3, 4, \ldots, d\}$, while the remaining diagonal entries are set to $1$. All other simulation settings remain the same as in the homoskedastic case.

\begin{table}[h!]
	\centering
	
	\begin{minipage}{\textwidth}
		\centering
		\small	
		\caption{Empirical power at the significance level $\alpha = 0.05$ for different mean structures and correlation levels under unequal variance. The highest power in each scenario is highlighted in bold, and deeper color intensity indicates higher power. When comparing methods with valid type-I error (Table~\ref{tab:type1-unequal}), DA methods perform very favorably across settings.} \label{tab:rejection-rates-unequal-var}
		\renewcommand{\arraystretch}{1.0}
		\setlength{\tabcolsep}{5pt}
		\begin{tabular}{lccccccccc}
		\toprule
		\multirow{2}{*}{\textbf{Method}} & \multicolumn{3}{c}{$\bm{\mu}^{(a)}$ + unequal variance} & \multicolumn{3}{c}{$\bm{\mu}^{(b)}$ + unequal variance} & \multicolumn{3}{c}{$\bm{\mu}^{(c)}$ + unequal variance} \\
		& $\rho=0$ & $\rho=0.4$ & $\rho=0.8$ & $\rho=0$ & $\rho=0.4$ & $\rho=0.8$ & $\rho=0$ & $\rho=0.4$ & $\rho=0.8$ \\
		\cmidrule(lr){2-4} \cmidrule(lr){5-7} \cmidrule(lr){8-10}
		\texttt{LOO} & \cellcolor{vlightorange}0.084 & \cellcolor{vlightorange}0.115 & \cellcolor{midorange}0.380
					& \cellcolor{vlightorange}0.000 & \cellcolor{vlightorange}0.001 & \cellcolor{lightorange}0.181
					& \cellcolor{lightorange}\textbf{0.258} & \cellcolor{midorange}\textbf{0.351} & \cellcolor{strongorange}\textbf{0.703} \\
		\texttt{Bonferroni} & \cellcolor{lightorange}0.171 & \cellcolor{lightorange}0.130 & \cellcolor{vlightorange}0.055
					& \cellcolor{lightorange}0.166 & \cellcolor{lightorange}0.103 & \cellcolor{vlightorange}0.030
					& \cellcolor{vlightorange}0.017 & \cellcolor{vlightorange}0.006 & \cellcolor{vlightorange}0.003 \\
		\texttt{csranks} & \cellcolor{lightorange}\textbf{0.184} & \cellcolor{midorange}\textbf{0.381} & \cellcolor{deeporange}\textbf{0.962}
					& \cellcolor{lightorange}0.162 & \cellcolor{midorange}0.363 & \cellcolor{deeporange}0.961
					& \cellcolor{vlightorange}0.019 & \cellcolor{vlightorange}0.041 & \cellcolor{lightorange}0.223 \\
		\texttt{MCS} & \cellcolor{vlightorange}0.004 & \cellcolor{vlightorange}0.002 & \cellcolor{vlightorange}0.004
					& \cellcolor{vlightorange}0.000 & \cellcolor{vlightorange}0.000 & \cellcolor{vlightorange}0.000
					& \cellcolor{lightorange}0.140 & \cellcolor{lightorange}0.156 & \cellcolor{lightorange}0.166 \\
		\texttt{DA-plug} & \cellcolor{vlightorange}0.049 & \cellcolor{vlightorange}0.052 & \cellcolor{vlightorange}0.042
					& \cellcolor{vlightorange}0.062 & \cellcolor{vlightorange}0.067 & \cellcolor{vlightorange}0.059
					& \cellcolor{vlightorange}0.098 & \cellcolor{lightorange}0.128 & \cellcolor{lightorange}0.202 \\
		\texttt{DA-plug}$^{\times 10}$ & \cellcolor{vlightorange}0.050 & \cellcolor{vlightorange}0.052 & \cellcolor{vlightorange}0.050
					& \cellcolor{vlightorange}0.080 & \cellcolor{vlightorange}0.080 & \cellcolor{vlightorange}0.073
					& \cellcolor{lightorange}0.125 & \cellcolor{lightorange}0.145 & \cellcolor{lightorange}0.240 \\
		\texttt{DA-adj} & \cellcolor{lightorange}0.122 & \cellcolor{midorange}0.259 & \cellcolor{deeporange}0.841
					& \cellcolor{lightorange}0.217 & \cellcolor{midorange}0.384 & \cellcolor{deeporange}0.916
					& \cellcolor{lightorange}0.135 & \cellcolor{lightorange}0.188 & \cellcolor{midorange}0.462 \\
		\texttt{DA-adj}$^{\times 10}$ & \cellcolor{lightorange}0.160 & \cellcolor{midorange}0.343 & \cellcolor{deeporange}0.946
					& \cellcolor{midorange}\textbf{0.294} & \cellcolor{strongorange}\textbf{0.517} & \cellcolor{deeporange}\textbf{0.982}
					& \cellcolor{lightorange}0.164 & \cellcolor{midorange}0.251 & \cellcolor{strongorange}0.605 \\
		\bottomrule
		\end{tabular}	
	\end{minipage}
	
	\vspace{1.5em}
	
	\begin{minipage}{\textwidth}
	\small
	\centering
	\caption{Empirical type I error at the significance level $\alpha = 0.05$ for different mean structures and correlation levels under unequal variance. Blue shading indicates over-rejection (liberal tests), green indicates under-rejection (conservative tests), and white indicates appropriate rejection rates (correct coverage).}
	\renewcommand{\arraystretch}{1.0}
	\setlength{\tabcolsep}{5pt}
	\begin{tabular}{lccc ccc ccc}
	\toprule
	\multirow{2}{*}{\textbf{Method}} 
	& \multicolumn{3}{c}{$\bm{\mu}^{(a,0)}$ + unequal variance} 
	& \multicolumn{3}{c}{$\bm{\mu}^{(b,0)}$ + unequal variance} 
	& \multicolumn{3}{c}{$\bm{\mu}^{(c,0)}$ + unequal variance} \\
	& $\rho=0$ & $\rho=0.4$ & $\rho=0.8$ 
	& $\rho=0$ & $\rho=0.4$ & $\rho=0.8$ 
	& $\rho=0$ & $\rho=0.4$ & $\rho=0.8$ \\
	\cmidrule(lr){2-4} \cmidrule(lr){5-7} \cmidrule(lr){8-10}
	\texttt{LOO} & \cellcolor{deepgreen}0.000 & \cellcolor{deepgreen}0.000 & \cellcolor{deepgreen}0.000 & \cellcolor{deepgreen}0.000 & \cellcolor{deepgreen}0.000 & \cellcolor{deepgreen}0.000 & \cellcolor{lightblue}0.070 & \cellcolor{vlightblue}0.064 & \cellcolor{vlightblue}0.065 \\
\texttt{Bonferroni} & \cellcolor{stronggreen}0.005 & \cellcolor{deepgreen}0.003 & \cellcolor{deepgreen}0.003 & \cellcolor{deepgreen}0.002 & \cellcolor{deepgreen}0.001 & \cellcolor{deepgreen}0.001 & \cellcolor{deepgreen}0.002 & \cellcolor{deepgreen}0.001 & \cellcolor{deepgreen}0.001 \\
\texttt{csranks} & \cellcolor{stronggreen}0.005 & \cellcolor{stronggreen}0.006 & \cellcolor{deepgreen}0.004 & \cellcolor{deepgreen}0.003 & \cellcolor{deepgreen}0.001 & \cellcolor{deepgreen}0.003 & \cellcolor{deepgreen}0.002 & \cellcolor{deepgreen}0.001 & \cellcolor{deepgreen}0.002 \\
\texttt{MCS} & \cellcolor{deepgreen}0.000 & \cellcolor{deepgreen}0.000 & \cellcolor{deepgreen}0.000 & \cellcolor{deepgreen}0.000 & \cellcolor{deepgreen}0.000 & \cellcolor{deepgreen}0.000 & 0.042 & 0.042 & \cellcolor{lightgreen}0.034 \\
\texttt{DA-plug} & \cellcolor{stronggreen}0.016 & \cellcolor{stronggreen}0.014 & \cellcolor{stronggreen}0.018 & \cellcolor{midgreen}0.023 & \cellcolor{midgreen}0.021 & \cellcolor{midgreen}0.023 & 0.048 & 0.052 & 0.048 \\
\texttt{DA-plug}$^{\times 10}$ & \cellcolor{stronggreen}0.011 & \cellcolor{stronggreen}0.013 & \cellcolor{stronggreen}0.012 & \cellcolor{midgreen}0.024 & \cellcolor{midgreen}0.022 & \cellcolor{midgreen}0.021 & 0.053 & 0.046 & 0.047 \\
\texttt{DA-adj} & \cellcolor{stronggreen}0.019 & \cellcolor{stronggreen}0.019 & \cellcolor{stronggreen}0.018 & \cellcolor{midgreen}0.027 & \cellcolor{midgreen}0.024 & \cellcolor{midgreen}0.029 & 0.054 & 0.052 & 0.050 \\
\texttt{DA-adj}$^{\times 10}$ & \cellcolor{stronggreen}0.015 & \cellcolor{stronggreen}0.012 & \cellcolor{stronggreen}0.012 & \cellcolor{midgreen}0.024 & \cellcolor{midgreen}0.024 & \cellcolor{midgreen}0.024 & 0.053 & 0.049 & 0.050 \\
	\bottomrule
	\end{tabular}
	\label{tab:type1-unequal}
	\end{minipage}
\end{table}

\Cref{tab:rejection-rates-unequal-var,tab:type1-unequal} reports the empirical performance of the methods under heteroskedastic variance settings. The results closely mirror those observed in the homoskedastic case, with \texttt{DA-adj}$^{\times 10}$ generally exhibiting strong power across most configurations and often achieving the highest power. Notably, the performance gap between \texttt{DA-plug} and \texttt{DA-adj} becomes more pronounced under heteroskedasticity, highlighting the advantage of noise-adjusted selection in the presence of non-uniform variances. As in the homoskedastic case, the \texttt{LOO} method attains the highest power in the scenario with $\bm{\mu}^{(c)}$, but this comes at the cost of inflated type I error rates, as evident in \Cref{tab:type1-unequal}. The \texttt{Bonferroni} procedure remains conservative, with limited detection power except in the $\bm{\mu}^{(a)}$ scenario without correlation. While \texttt{csranks} performs well in that particular setting, it generally underperforms relative to \texttt{DA-adj}$^{\times 10}$ in other configurations. As in the homoskedastic cases, the \texttt{MCS} method exhibits limited power in the first two scenarios, while it performs reasonably well in the last scenario with $\bm{\mu}^{(c)}$.

\subsection{Power and validity in high-dimensional settings} \label{Sec: highdim}
\begin{figure}[h!]
	\centering
	\includegraphics[width=0.92\textwidth]{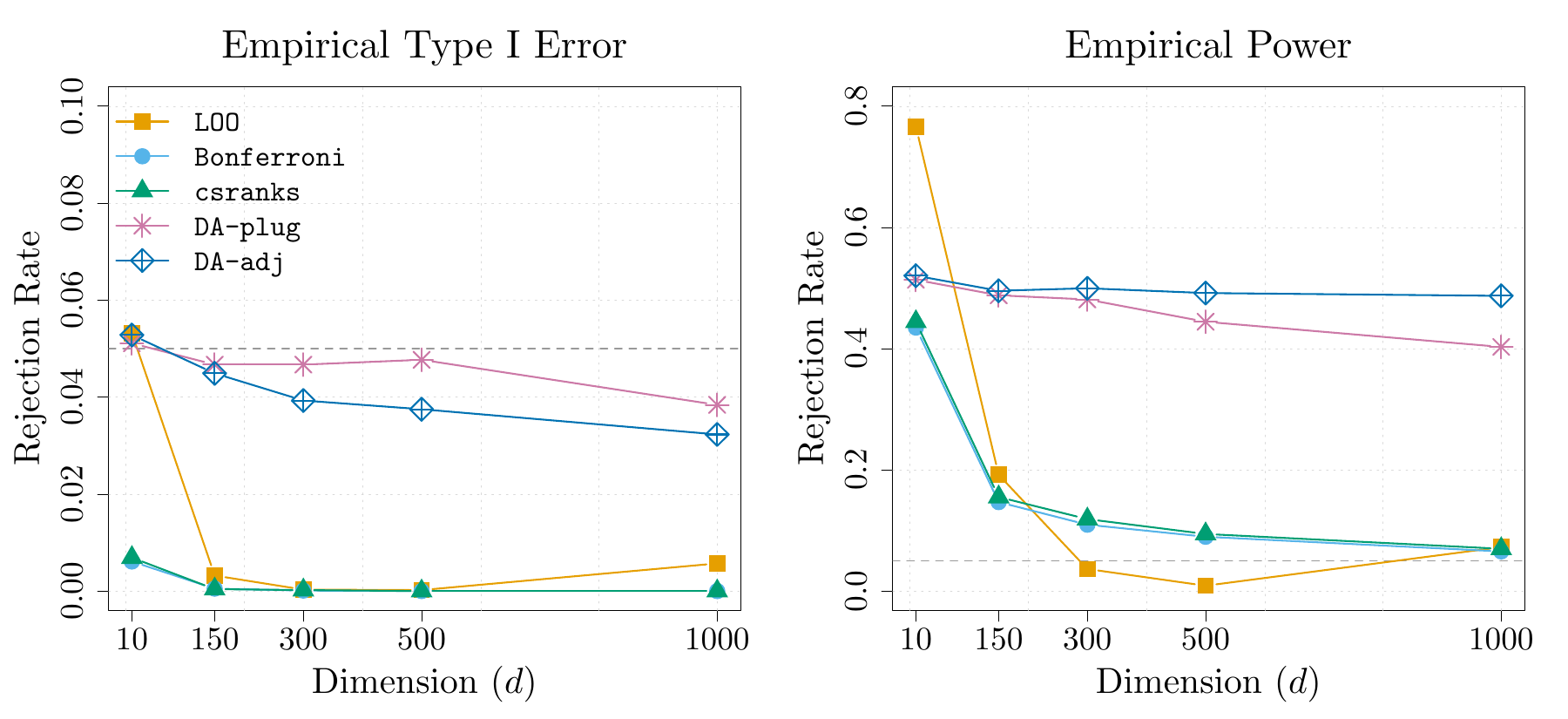}
	\caption{Empirical type I error rates (left) and power (right) of the considered methods across varying dimensions under the settings described in \Cref{Sec: highdim}. The dashed line represents the nominal level of 0.05. The results demonstrate the superior performance of the DA argmin tests in the considered high-dimensional settings, consistently maintaining strong power across all dimensions while controlling the type I error.}
	\label{fig:high-dim}
\end{figure}

We next investigate the performance of the considered methods across varying dimensional settings to assess their sensitivity to problem dimensionality. Specifically, we consider dimensions $d \in \{10, 150, 300, 500, 1000\}$ and evaluate the empirical rejection rates under the following configuration. The mean vector is set to $\bm{\mu} = (0, 0, 1, 1, \ldots, 1)^\top$ under the null and $\bm{\mu} = (0.15, 0, 1, 1, \ldots, 1)^\top$ under the alternative. The covariance matrix $\bm{\Sigma}$ is diagonal with entries $\Sigma_{kk} = 1$ for $k = \{1, 2\}$, and $\Sigma_{kk} = 20$ for $k \in \{3, \ldots, d\}$. The sample size is fixed at $n = 500$, and the significance level is set to $\alpha = 0.05$. The results shown in \Cref{fig:high-dim} are averaged over $10,\!000$ replications.

The left panel of \Cref{fig:high-dim} presents the empirical rejection rates under the null hypothesis. All methods adequately control the type I error rate below the nominal level of $0.05$ across all dimensions. Notably, the tests tend to become increasingly conservative as dimensionality grows, with the \texttt{DA-plug} and \texttt{DA-adj} methods exhibiting relatively less conservativeness compared to the others.

The right panel of \Cref{fig:high-dim} shows the empirical power of the methods under the alternative hypothesis. In the low-dimensional setting ($d = 10$), the \texttt{LOO} method achieves the highest power, followed by the proposed DA argmin tests (\texttt{DA-plug} and \texttt{DA-adj}). Interestingly, this trend reverses in higher dimensions, where the power of \texttt{LOO} deteriorates rapidly, becoming nearly close to the nomial level $\alpha$. This phenomenon may be attributed to the weighting nature of \texttt{LOO}, which assigns non-negligible weight to irrelevant components in high-dimensional settings. Similar patterns are observed for the \texttt{Bonferroni} and \texttt{csranks} methods, whose power also declines substantially as the dimension increases. While the power of \texttt{DA-plug} and \texttt{DA-adj} exhibits a mild decrease with dimensionality, these methods consistently outperform the others and remain competitive throughout. Between the two DA argmin tests, the noise-adjusted version (\texttt{DA-adj}) tends to have slightly higher power, particularly in higher dimensions.

The aggregated versions (\texttt{DA-plug}$^{\times 10}$ and \texttt{DA-adj}$^{\times 10}$) are not included in \Cref{fig:high-dim} due to their computational cost. However, given their strong performance in previous experiments, we expect them to yield even higher power in high dimensions while still controlling type I error. The \texttt{MCS} method is similarly excluded from this analysis due to its intensive computational demands.

\subsection{Real world data example} \label{Sec: real data}
We revisit the classification competition datasets analyzed by \citet{zhang2024winners} to illustrate the performance of the proposed DA argmin tests in a real-world setting. In these competitions, students trained classification models on a provided training dataset and subsequently predicted labels on a separate test dataset. The competitions took place in 2023 and 2024, attracting submissions of 44 and 39 prediction models, respectively. Model performance was evaluated based on binary classification errors, encoded as $0$ (correct) or $1$ (incorrect), on test datasets of size 183 in 2023 and 1236 in 2024. The primary objective was to identify the best-performing model, i.e., the one with the lowest classification error, and to construct a $95\%$ confidence set for this model. A detailed description of the datasets is available in \citet{zhang2024winners}.

We apply our proposed \texttt{DA-adj}$^{\times 50}$ method, a variant of the previously introduced \texttt{DA-adj}$^{\times 10}$ procedure but employing 50 random data splits to ensure stable inference. Since the performance of the simpler \texttt{DA-plug} variant was similar, we focus our presentation solely on the \texttt{DA-adj} method. Our method is compared against four established procedures: \texttt{LOO}, \texttt{Bonferroni}, \texttt{csranks}, and \texttt{MCS}. Given that several of these methods—including ours—depend on random data splits and thus can yield varying results, we follow \citet{zhang2024winners} and report averages computed over 100 replications.

Our results indicate that the \texttt{DA-adj} method consistently produces smaller inclusion sets across both competition years compared to the alternatives. This advantage is particularly pronounced in the 2023 competition, where the average size of the inclusion set produced by \texttt{DA-adj} is $17.35_{\pm 1.17}$ with the number after $\pm$ indicating the standard deviation. This value is substantially lower than the inclusion set sizes produced by \texttt{LOO} ($32.55_{\pm 1.14}$), \texttt{csranks} ($38.56_{\pm 0.61}$), \texttt{Bonferroni} ($41_{\pm 0}$), and \texttt{MCS} ($43_{\pm 0}$). Similarly, in the 2024 competition, our method continues to outperform its counterparts, achieving an average inclusion set size of $19.93_{\pm 0.57}$, compared to \texttt{LOO} ($25.48_{\pm 1.80}$), \texttt{csranks} ($28.69_{\pm 0.87}$), \texttt{Bonferroni} ($30_{\pm 0}$), and \texttt{MCS} ($37_{\pm 0}$).

\begin{figure}
	\centering
	\includegraphics[width=1\textwidth]{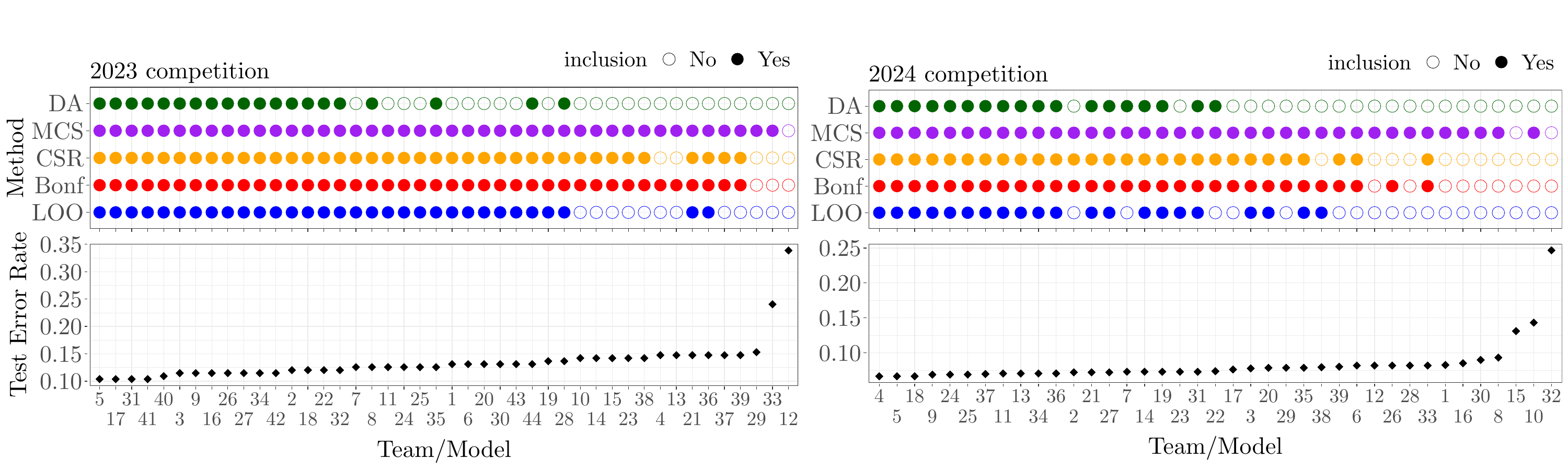}
	\caption{
		Comparison of the inclusion sets generated by the \texttt{DA-adj}$^{\times 50}$ method (DA) and other established methods for the 2023 (left) and 2024 (right) classification competitions. The methods compared are: MCS (\texttt{MCS}), CSR (\texttt{csranks}), Bonf (\texttt{Bonferroni}), and LOO (\texttt{LOO}). Each inclusion set is depicted as a colored interval. Our \texttt{DA-adj}$^{\times 50}$ method consistently produces smaller inclusion sets compared to other approaches, highlighting its superior efficiency in pinpointing the best-performing model.}
	\label{fig:real-data}
\end{figure}

\Cref{fig:real-data} illustrates a representative realization of the inclusion sets from each method. Notably, like other methods, the confidence set produced by \texttt{DA-adj}$^{\times 50}$ does not form a single interval. This discontinuity arises primarily because the significance of each test depends not only on differences in means but also intricately on their correlations with the minimum mean. Confidence sets for the worst-performing model, obtained through argmax inference, are presented separately in \Cref{Sec: real-data-argmax}.

\subsection{Simulations on DA-MCS} \label{Sec: Sim MCS}
In this subsection, we present simulation results for the DA-MCS method developed in \Cref{Sec: MCS}. The simulation settings closely follow those described in \Cref{Sec: power}, with the key difference being the specification of the mean vector:
\begin{align*}
	\bm{\mu} = (\underbrace{0,\ldots,0}_{|\Theta| \text{ entries}},\underbrace{\zeta,\ldots,\zeta}_{d - |\Theta| \text{ entries}})^\top \in \mathbb{R}^d,
\end{align*}
where the mean gap parameter $\zeta$ is set to $10\sqrt{\log|\Theta|/(2n)}$. Simulation results for other choices of $\zeta$ can be found in \Cref{Sec: additional-simulation-mcs}. We vary the size of the argmin set $\Theta$ over $|\Theta| \in \{2, 5, 10, 15, 20\}$ and consider three correlation levels $\rho \in \{0, 0.4, 0.8\}$. For each setting, we compare six methods described below in terms of their empirical uniform coverage rates at the nominal level $\alpha = 0.05$, computed as the proportion of simulations in which the true parameter set $\Theta$ is contained in the estimated set $\widehat{\Theta}$, i.e., $P(\Theta \subseteq \widehat{\Theta})$ based on $10,\!000$ repetitions. We also report the average length of the confidence sets.

In addition to the pointwise methods, namely \texttt{DA-plug} and \texttt{DA-adj} defined earlier, we consider four methods targeting uniform coverage.
\begin{itemize}
	\item \texttt{DA-MCS-plug}$^{1}$: The DA-MCS method with the plug-in selection rule $\widehat{s}_{\mathrm{plug}}$ and a one-step construction where the significance level is adjusted to $\alpha/d$ as detailed in \Cref{Sec: MCS}.
	\item \texttt{DA-MCS-adj}$^{1}$: The DA-MCS method with the noise-adjusted selection rule $\widehat{s}_{\mathrm{adj}}$ and a one-step construction where the significance level is adjusted to $\alpha/d$ as detailed in \Cref{Sec: MCS}.
	\item \texttt{DA-MCS-plug}$^{2}$: The DA-MCS method with the plug-in selection rule $\widehat{s}_{\mathrm{plug}}$ and a two-step construction where the significance level is adjusted to $\alpha/|\widehat{\Theta}^{(2)}|$ as detailed in \Cref{Sec: MCS}.
	\item \texttt{DA-MCS-adj}$^{2}$: The DA-MCS method with the noise-adjusted selection rule $\widehat{s}_{\mathrm{adj}}$ and a two-step construction where the significance level is adjusted to $\alpha/|\widehat{\Theta}^{(2)}|$ as detailed in \Cref{Sec: MCS}.
\end{itemize}

The simulation results are presented in \Cref{tab:uniform-coverage-10}. The first row reports the coverage rates of methods designed for pointwise coverage, while the bottom two rows report those of methods targeting uniform coverage. The results demonstrate that the DA-MCS method with the two-step construction consistently achieves superior uniform coverage compared to its one-step counterpart, with empirical coverage rates close to the nominal level $1-\alpha$ across all settings. In contrast, the DA methods tailored for pointwise coverage, namely \texttt{DA-plug} and \texttt{DA-adj}, exhibit substantially lower coverage, especially when the cardinality of the argmin set $|\Theta|$ is large. These findings highlight the importance of aligning the inferential method with the desired coverage objective: while pointwise methods offer greater power in terms of rejection rates, they may fail to ensure uniform coverage. Conversely, the DA-MCS methods provide valid uniform coverage, but can be overly conservative depending on the application.

\begin{table}[h]
    \centering
    \small
    \caption{Empirical coverage probabilities $P(\Theta \subseteq \widehat{\Theta})$ across varying cardinalities $|\Theta|$ and correlation levels $\rho$, evaluated for six different methods at the nominal level $1-\alpha = 0.95$. The mean gap $\zeta$ is set to $10\sqrt{\log|\Theta|/(2n)}$. Numbers in parentheses indicate the average length of the confidence sets. Under-coverage rates are shaded in progressively darker blue, over-coverage rates in progressively darker green, and rates close to the nominal level remain unshaded.}
    \label{tab:uniform-coverage-10}
    \renewcommand{\arraystretch}{1.0}
    \setlength{\tabcolsep}{2.7pt}

	\begin{tabular}{cc}

		\begin{minipage}{0.45\textwidth}
		\centering
		\texttt{DA-plug} (pointwise)\\[0.3em]
		\begin{tabular}{llll}\toprule
		  & $\rho=0$                     & $\rho=0.4$                   & $\rho=0.8$                   \\\midrule
		  $|\Theta|=2$  & \cellcolor{lightblue}0.897 {\tiny(2.51)} & \cellcolor{lightblue}0.898 {\tiny(2.49)} & \cellcolor{lightblue}0.896 {\tiny(2.46)} \\
		  $|\Theta|=5$  & \cellcolor{lightblue}0.804 {\tiny(4.75)} & \cellcolor{lightblue}0.804 {\tiny(4.75)} & \cellcolor{lightblue}0.812 {\tiny(4.74)} \\
		  $|\Theta|=10$ & \cellcolor{midblue}0.722 {\tiny(9.50)} & \cellcolor{midblue}0.723 {\tiny(9.50)} & \cellcolor{midblue}0.732 {\tiny(9.49)} \\
		  $|\Theta|=15$ & \cellcolor{strongblue}0.657 {\tiny(14.23)} & \cellcolor{strongblue}0.661 {\tiny(14.24)} & \cellcolor{strongblue}0.688 {\tiny(14.25)} \\
		  $|\Theta|=20$ & \cellcolor{strongblue}0.620 {\tiny(19.01)} & \cellcolor{strongblue}0.626 {\tiny(19.00)} & \cellcolor{strongblue}0.658 {\tiny(19.00)} \\\bottomrule
		\end{tabular}
		\end{minipage}
		&
		\begin{minipage}{0.45\textwidth}
		\centering
		\texttt{DA-adj} (pointwise)\\[0.3em]
		\begin{tabular}{llll}\toprule
		  & $\rho=0$                       & $\rho=0.4$                     & $\rho=0.8$                     \\\midrule
		  $|\Theta|=2$  & \cellcolor{vlightblue}0.903 {\tiny(2.53)} & \cellcolor{vlightblue}0.900 {\tiny(2.51)} & \cellcolor{lightblue}0.897 {\tiny(2.50)} \\
		  $|\Theta|=5$  & \cellcolor{lightblue}0.804 {\tiny(4.75)} & \cellcolor{lightblue}0.802 {\tiny(4.75)} & \cellcolor{lightblue}0.808 {\tiny(4.75)} \\
		  $|\Theta|=10$ & \cellcolor{midblue}0.717 {\tiny(9.51)} & \cellcolor{midblue}0.718 {\tiny(9.48)} & \cellcolor{midblue}0.705 {\tiny(9.48)} \\
		  $|\Theta|=15$ & \cellcolor{strongblue}0.653 {\tiny(14.25)} & \cellcolor{strongblue}0.636 {\tiny(14.25)} & \cellcolor{strongblue}0.631 {\tiny(14.23)} \\
		  $|\Theta|=20$ & \cellcolor{strongblue}0.613 {\tiny(19.02)} & \cellcolor{strongblue}0.597 {\tiny(19.01)} & \cellcolor{strongblue}0.585 {\tiny(18.96)} \\\bottomrule
		\end{tabular}
		\end{minipage} \\[5em]
		
		\begin{minipage}{0.45\textwidth}
		\centering
		\texttt{DA-MCS-plug}$^{1}$ (uniform)\\[0.3em]
		\begin{tabular}{llll}\toprule
		  & $\rho=0$                       & $\rho=0.4$                     & $\rho=0.8$                     \\\midrule
		  $|\Theta|=2$  & \cellcolor{midgreen}0.999 {\tiny(21.49)} & \cellcolor{midgreen}0.999 {\tiny(20.84)} & \cellcolor{midgreen}0.998 {\tiny(19.30)} \\
		  $|\Theta|=5$  & \cellcolor{midgreen}0.998 {\tiny(5.11)}  & \cellcolor{midgreen}0.996 {\tiny(5.10)}  & \cellcolor{midgreen}0.997 {\tiny(5.11)}  \\
		  $|\Theta|=10$ & \cellcolor{midgreen}0.995 {\tiny(9.99)} & \cellcolor{midgreen}0.996 {\tiny(9.99)} & \cellcolor{midgreen}0.995 {\tiny(9.99)} \\
		  $|\Theta|=15$ & \cellcolor{midgreen}0.993 {\tiny(14.99)} & \cellcolor{midgreen}0.994 {\tiny(14.99)} & \cellcolor{midgreen}0.995 {\tiny(14.99)} \\
		  $|\Theta|=20$ & \cellcolor{midgreen}0.992 {\tiny(19.99)} & \cellcolor{midgreen}0.994 {\tiny(19.99)} & \cellcolor{midgreen}0.995 {\tiny(19.99)} \\\bottomrule
		\end{tabular}
		\end{minipage}
		&
		\begin{minipage}{0.45\textwidth}
		\centering
		\texttt{DA-MCS-adj}$^{1}$ (uniform)\\[0.3em]
		\begin{tabular}{llll}\toprule
		  & $\rho=0$                       & $\rho=0.4$                     & $\rho=0.8$                     \\\midrule
		  $|\Theta|=2$  & \cellcolor{midgreen}0.998 {\tiny(21.03)} & \cellcolor{midgreen}0.999 {\tiny(20.68)} & \cellcolor{midgreen}0.999 {\tiny(19.31)} \\
		  $|\Theta|=5$  & \cellcolor{midgreen}0.998 {\tiny(5.10)}  & \cellcolor{midgreen}0.997 {\tiny(5.11)}  & \cellcolor{midgreen}0.998 {\tiny(5.09)}  \\
		  $|\Theta|=10$ & \cellcolor{midgreen}0.996 {\tiny(9.99)} & \cellcolor{midgreen}0.993 {\tiny(9.99)} & \cellcolor{midgreen}0.996 {\tiny(9.99)} \\
		  $|\Theta|=15$ & \cellcolor{midgreen}0.993 {\tiny(14.99)} & \cellcolor{midgreen}0.995 {\tiny(14.99)} & \cellcolor{midgreen}0.994 {\tiny(14.99)} \\
		  $|\Theta|=20$ & \cellcolor{midgreen}0.991 {\tiny(19.99)} & \cellcolor{midgreen}0.993 {\tiny(19.99)} & \cellcolor{midgreen}0.985 {\tiny(19.99)} \\\bottomrule
		\end{tabular}
		\end{minipage}\\[5em]
		
		\begin{minipage}{0.45\textwidth}
		\centering
		\texttt{DA-MCS-plug}$^{2}$ (uniform)\\[0.3em]
		\begin{tabular}{llll}\toprule
		  & $\rho=0$                     & $\rho=0.4$                     & $\rho=0.8$                     \\\midrule
		  $|\Theta|=2$  & \cellcolor{midgreen}0.986 {\tiny(8.15)}  & \cellcolor{midgreen}0.985 {\tiny(7.92)}  & \cellcolor{midgreen}0.983 {\tiny(7.26)}  \\
		  $|\Theta|=5$  & 0.956 {\tiny(4.95)}              & \cellcolor{lightgreen}0.962 {\tiny(4.95)} & \cellcolor{lightgreen}0.962 {\tiny(4.95)} \\
		  $|\Theta|=10$ & 0.958 {\tiny(9.94)}              & \cellcolor{lightgreen}0.961 {\tiny(9.99)} & 0.958 {\tiny(9.99)}              \\
		  $|\Theta|=15$ & 0.957 {\tiny(14.95)}             & \cellcolor{lightgreen}0.960 {\tiny(14.99)} & \cellcolor{lightgreen}0.967 {\tiny(14.94)} \\
		  $|\Theta|=20$ & \cellcolor{lightgreen}0.963 {\tiny(19.95)} & \cellcolor{lightgreen}0.963 {\tiny(19.94)} & \cellcolor{lightgreen}0.969 {\tiny(19.94)} \\\bottomrule
		\end{tabular}
		\end{minipage}
		&
		\begin{minipage}{0.45\textwidth}
		\centering
		\texttt{DA-MCS-adj}$^{2}$ (uniform)\\[0.3em]
		\begin{tabular}{llll}\toprule
		  & $\rho=0$                     & $\rho=0.4$                      & $\rho=0.8$                     \\\midrule
		  $|\Theta|=2$  & \cellcolor{midgreen}0.986 {\tiny(8.30)}  & \cellcolor{midgreen}0.985 {\tiny(8.08)}  & \cellcolor{midgreen}0.984 {\tiny(7.66)}  \\
		  $|\Theta|=5$  & 0.954 {\tiny(4.95)}              & 0.952 {\tiny(4.95)}               & 0.956 {\tiny(4.95)}              \\
		  $|\Theta|=10$ & 0.955 {\tiny(9.94)}              & \cellcolor{lightgreen}0.961 {\tiny(9.99)} & 0.957 {\tiny(9.94)}              \\
		  $|\Theta|=15$ & \cellcolor{lightgreen}0.962 {\tiny(14.94)} & 0.958 {\tiny(14.95)}              & \cellcolor{lightgreen}0.962 {\tiny(14.95)} \\
		  $|\Theta|=20$ & \cellcolor{lightgreen}0.963 {\tiny(19.95)} & \cellcolor{lightgreen}0.961 {\tiny(19.94)} & \cellcolor{lightgreen}0.963 {\tiny(19.94)} \\\bottomrule
		\end{tabular}
		\end{minipage}
		
		\end{tabular}
\end{table}

\section{Conclusion} \label{Sec: conclusion}
In this work, we proposed a DA method for the high-dimensional argmin inference problem that remains valid regardless of how the dimensionality scales with the sample size. We characterized the minimax separation rate for this problem and established its fundamental dependence on the cardinality of the confusion set. Furthermore, we showed that both the plug-in and noise-adjusted versions of our procedure adapt to the underlying confusion set and achieve minimax rate-optimal power. Our simulation study confirms the robustness of the proposed tests, which maintain the nominal level and exhibit strong power across a range of signal structures and correlation levels.

There are several promising avenues for future research. First, it would be valuable to extend our framework to general rank-$k$ inference problems, where the objective is to identify the index corresponding to the $k$-th smallest mean. Such an extension would broaden the applicability of our methodology and introduce new theoretical challenges. Second, it may be worthwhile to explore thresholding-based approaches for constructing $\bm{\gamma}_{\widehat{s}}$ in our test statistic. Specifically, rather than selecting a single index, one could include all indices whose means fall below a pre-specified threshold. This strategy may offer greater power, particularly in cases where multiple indices attain the minimum. Lastly, developing faster algorithms for the multiple-split procedure would also be a valuable direction for future work.

\paragraph{Acknowledgements}
We are grateful to the authors of~\cite{zhang2024winners} for kindly sharing the code and data used in the simulation study. 

\bibliographystyle{apalike}
\bibliography{reference}

\clearpage 
\appendix

\section{Proofs and technical lemmas} \label{App: Proof}
In this section, we collect the proofs of the main results and some technical lemmas. 

\subsection{Proof of \Cref{The: DA validity}}
This result is almost a direct consequence of the Berry--Esseen bound for Student's $t$-statistic~\citep{bentkus1996berry}, as similarly used in many past works on DA inference.
By the Berry--Esseen bound for Student's $t$-statistic~\citep[][Theorem 1.2]{bentkus1996berry}, we have that, conditional on $\widehat{s}$, which is independent of the first half of the data,
\begin{align*}
	\sup_{P \in \mathcal{P}_{0,r}} \sup_{t \in \mathbb{R}} \Bigg| P \Biggl(\frac{\sqrt{n}\bm{\gamma}_{\widehat{s}}^\top \bigl(\overline{\bX}^{(1)} -  \bm{\mu} \bigr) }{\sqrt{\bm{\gamma}_{\widehat{s}}^\top \widehat{\bm{\Sigma}}^{(1)} \bm{\gamma}_{\widehat{s}}}} \leq t \,\Bigg|\, \widehat{s} \Biggr) - \Phi(t) \Bigg| \leq  \min \Bigl\{1, C M_{\widehat{s}}\Bigr\} \leq \min \biggl\{1, C \!\!\max_{k \in [d] \setminus \{r\}} M_k\biggr\}.
\end{align*}
Now the result follows by taking the expectation over $\widehat{s}$ and noting that 
\begin{align*}
	& \sup_{P \in \mathcal{P}_{0,r}} \sup_{t \in \mathbb{R}} \Bigg| P \Biggl(\frac{\sqrt{n}\bm{\gamma}_{\widehat{s}}^\top \bigl(\overline{\bX}^{(1)} -  \bm{\mu} \bigr) }{\sqrt{\bm{\gamma}_{\widehat{s}}^\top \widehat{\bm{\Sigma}}^{(1)} \bm{\gamma}_{\widehat{s}}}} \leq t \Biggr) - \Phi(t) \Bigg| \\
	&  \leq \mE_{P}\Biggl[\sup_{P \in \mathcal{P}_{0,r}} \sup_{t \in \mathbb{R}} \Bigg| P \Biggl(\frac{\sqrt{n}\bm{\gamma}_{\widehat{s}}^\top \bigl(\overline{\bX}^{(1)} -  \bm{\mu} \bigr) }{\sqrt{\bm{\gamma}_{\widehat{s}}^\top \widehat{\bm{\Sigma}}^{(1)} \bm{\gamma}_{\widehat{s}}}} \leq t \,\Bigg|\, \widehat{s} \Biggr) - \Phi(t) \Bigg|\Biggr],
\end{align*}		
where the expectation outside is taken with respect to the randomness in $\widehat{s}$. This completes the proof of \Cref{The: DA validity}.

\paragraph{Remark.} Our proof of validity is straightforward and transparent, relying only on a conditional central limit theorem for student’s $t$-statistic. This simplicity may be viewed as an additional advantage of our approach. By contrast, existing validity proofs in the literature often require intricate arguments and heavy technical machinery, which can make them less accessible. For instance, \citet{zhang2024winners} establishes validity through a central limit theorem for cross-validation-type statistics, a setting that is substantially more challenging due to the dependence among the summands.

\subsection{Proof of \Cref{Thm: upper bound}}
We start focusing our analysis on the case $\widehat{s} = \widehat{s}_{\mathrm{plug}}$ in \Cref{App: plug} and then turn to the case $\widehat{s} = \widehat{s}{_\mathrm{adj}}$ in \Cref{App: adj}. 

\subsubsection{Proof for plug-in estimator $\widehat{s}_{\mathrm{plug}}$} \label{App: plug}
We now present the proof of \Cref{Thm: upper bound} by focusing first on the case where $\widehat{s} = \widehat{s}_{\mathrm{plug}}$. To simplify the notation and streamline the argument, we set $r=1$ without loss of generality and assume that $\mu_2 \leq \mu_3 \leq \ldots \leq \mu_d$ throughout the proof. For simplicity, we write $\mathbb{C}_1 = \mathbb{C}$. Given $\delta > 0$, which will be specified later, define the two events
\begin{align*}
	& \mathcal{E}_{1,\delta} \coloneqq \bigg\{ \bm{\gamma}_{\widehat{s}}^\top \widehat{\bm{\Sigma}}^{(1)} \bm{\gamma}_{\widehat{s}} \leq \frac{4\sigma^2}{\delta} \bigg\} \quad \text{and} \quad  \mathcal{E}_{2,\delta} \coloneqq \bigg\{ \bigg| \sqrt{n}\bm{\gamma}_{\widehat{s}}^\top \bigl(\overline{\bX}^{(1)} - \bm{\mu} \bigr) \bigg| \leq \sqrt{\frac{4\sigma^2}{\delta}} \bigg\}.
\end{align*}
Each of these events holds with probability at least $1 - \delta$, which can be verified by applying Markov's and Chebyshev's inequalities (conditional on $\hat s$) along with the inequality that $\mV(W_1-W_2) \leq 2 \mV(W_1) + 2\mV(W_2)$ for any random variables $W_1$ and $W_2$. 

Invoking the union bound, the type II error of the test under any distribution $P \in \mathcal{P}_{1,r}(\varepsilon;\tau)$ is bounded by
\begin{align*}
		P \Big( \sqrt{n}\bm{\gamma}_{\widehat{s}}^\top \overline{\bX}^{(1)} \leq z_{1-\alpha} \sqrt{\bm{\gamma}_{\widehat{s}}^\top \widehat{\bm{\Sigma}}^{(1)} \bm{\gamma}_{\widehat{s}}} \Big) &
	\leq \, P \Big( \sqrt{n}\bm{\gamma}_{\widehat{s}}^\top \overline{\bX}^{(1)} \leq z_{1-\alpha} \sqrt{4\sigma^2\delta^{-1}} \Big) + P(\mathcal{E}_{1,\delta}^c) \\[.5em]
	\leq \, & P \Bigl( \sqrt{n}(\mu_1 - \mu_{\widehat{s}}) \leq (z_{1-\alpha} + 1) \sqrt{4\sigma^2\delta^{-1}} \Bigr) + P(\mathcal{E}_{1,\delta}^c) + P(\mathcal{E}_{2,\delta}^c) \\[.5em]
	= \, & \underbrace{P \Bigl( \sqrt{n}(\mu_1 - \mu_{\widehat{s}}) \leq (z_{1-\alpha} + 1) \sqrt{4\sigma^2\delta^{-1}} \,\cap\, \widehat{s} \in  \mathbb{C} \Bigr)}_{\coloneqq \mathrm{(I)}} \\
	& + \underbrace{P \Bigl( \sqrt{n}(\mu_1 - \mu_{\widehat{s}}) \leq (z_{1-\alpha} + 1) \sqrt{4\sigma^2\delta^{-1}} \,\cap\, \widehat{s} \in \mathbb{C}^c \Bigr)}_{\coloneqq \mathrm{(II)}} + 2\delta. 
\end{align*}
It remains to show that each term vanishes under the condition of the theorem. 

\medskip

\noindent \textbf{Term (I):}
Starting with the first term $\mathrm{(I)}$, define the event $\mathcal{E}_{3,\delta}$ as
\begin{align*}
	\mathcal{E}_{3,\delta} \coloneqq \bigcap_{k \in \mathbb{C} \cup \{2\}}\biggl\{\big|\overline{X}^{(2)}_k - \mu_k \big| < \sqrt{\frac{2\sigma^2}{n}\log\biggl(\frac{2|\mathbb{C} \cup \{2\}|}{\delta}\biggr)} \biggr\}, 
\end{align*}
which holds with probability at least $1 - \delta$, as can be verified by using a standard sub-Gaussian tail bound \citep[e.g.,][Proposition 2.5]{wainwright2019high} and the union bound.

On the event $\mathcal{E}_{3,\delta} \cap \{\widehat{s} \in \mathbb{C}\}$, we have
\begin{align*}
	\mu_{\widehat{s}} \leq \, & \overline{X}^{(2)}_{\widehat{s}} + \sqrt{\frac{2\sigma^2}{n}\log\biggl(\frac{2|\mathbb{C}\cup \{2\}|}{\delta}\biggr)} \leq \overline{X}^{(2)}_{2} + \sqrt{\frac{2\sigma^2}{n}\log\biggl(\frac{2|\mathbb{C}\cup \{2\}|}{\delta}\biggr)} \\
	\leq \, & \mu_{2} + 2\sqrt{\frac{2\sigma^2}{n}\log\biggl(\frac{2|\mathbb{C}\cup \{2\}|}{\delta}\biggr)}.
\end{align*}  
Hence it holds that 
\begin{align*}
	\mathrm{(I)} \leq \, & P \Bigl( \sqrt{n}(\mu_1 - \mu_{\widehat{s}}) \leq (z_{1-\alpha} + 1) \sqrt{4\sigma^2\delta^{-1}} \,\cap\, \mathcal{E}_{3,\delta} \,\cap\, \{\widehat{s} \in \mathbb{C}\}\Bigr) + P(\mathcal{E}_{3,\delta}^c) \\[.5em]
	\leq \, & P\Biggl( \sqrt{n}(\mu_1 - \mu_{2}) \leq (z_{1-\alpha} + 1) \sqrt{4\sigma^2\delta^{-1}} + \sqrt{8\sigma^2\log\biggl(\frac{2|\mathbb{C}\cup \{2\}|}{\delta}\biggr)}\Biggr) + \delta.
\end{align*}

\smallskip

\noindent \textbf{Term (II):} Next for the second term, write $\mathbb{C}^c = \mathbb{C}_a^c \cup \mathbb{C}_b^c$ where
\begin{align*}
	& \mathbb{C}_a^c = \bigg\{k \in [d]\!\setminus\! (\{1\} \cup \Theta_{-1}): \frac{\mu_1 - \mu_2}{2} > \mu_k - \mu_2\bigg\} \ \text{and} \\
	& \mathbb{C}_b^c = \bigg\{k \in [d]\!\setminus\! (\{1\} \cup \Theta_{-1}) : \mu_k - \mu_2 > C_n \sqrt{\frac{\log (d)}{n}}\bigg\},
\end{align*}
so that we have
\begin{align*}
	\mathrm{(II)} \leq \, & P \Bigl( \sqrt{n}(\mu_1 - \mu_{\widehat{s}}) \leq (z_{1-\alpha} + 1) \sqrt{4\sigma^2\delta^{-1}} \,\cap\,\widehat{s} \in \mathbb{C}_a^c \Bigr) \\
	& + P \Bigl( \sqrt{n}(\mu_1 - \mu_{\widehat{s}}) \leq (z_{1-\alpha} + 1) \sqrt{4\sigma^2\delta^{-1}} \,\cap\, \widehat{s} \in \mathbb{C}_b^c \Bigr) \\[.5em]
	\leq \, & P \biggl( \frac{\sqrt{n}}{2}(\mu_1 - \mu_2) \leq (z_{1-\alpha} + 1) \sqrt{4\sigma^2\delta^{-1}} \biggr) + P(\widehat{s} \in \mathbb{C}_b^c).
\end{align*}
To deal with $P(\widehat{s} \in \mathbb{C}_b^c)$, define the event
\begin{align*}
	\mathcal{E}_{4,\delta} \coloneqq \bigcap_{k=2}^d\biggl\{\bigg|\overline{X}^{(2)}_k - \mu_k \bigg| < \sqrt{\frac{2\sigma^2}{n}\log\biggl(\frac{2d}{\delta}\biggr)} \biggr\}.
\end{align*}
Another application of the sub-Gaussian tail bound together with the union bound yields 
\begin{align*}
	P(\mathcal{E}_{4,\delta}^c) = P\Biggl( \bigcup_{k=2}^d\biggl\{\bigg|\overline{X}^{(2)}_k - \mu_k \bigg| \geq \sqrt{\frac{2\sigma^2}{n}\log\biggl(\frac{2d}{\delta}\biggr)} \biggr\}\Biggr) \leq \delta.
\end{align*}
From this, we obtain that 
\begin{align*}
	P(\widehat{s} \in \mathbb{C}_b^c) \leq \, & P \biggl(\mu_{\widehat{s}} - \mu_2 > C_n\sqrt{\frac{\log (d)}{n}} \,\cap\, \mathcal{E}_{4,\delta} \biggr) + P(\mathcal{E}_{4,\delta}^c) \\[.5em]	
	\leq \, & P \biggl(\mu_{\widehat{s}} - \mu_2 > C_n \sqrt{\frac{\log (d)}{n}} \,\cap\, \mathcal{E}_{4,\delta} \biggr) + \delta \\[.5em]
	\leq \, & P \Biggl(2\sqrt{\frac{2\sigma^2}{n}\log\biggl(\frac{2d}{\delta}\biggr)} > C_n \sqrt{\frac{\log (d)}{n}}\Biggr) + \delta,
\end{align*}
where the last inequality holds since under the event $\mathcal{E}_{4,\delta}$,
\begin{align*}
	\mu_{\widehat{s}} \leq \, & \overline{X}^{(2)}_{\widehat{s}} + \sqrt{\frac{2\sigma^2}{n}\log\biggl(\frac{2d}{\delta}\biggr)} \leq \overline{X}^{(2)}_{2} + \sqrt{\frac{2\sigma^2}{n}\log\biggl(\frac{2d}{\delta}\biggr)} \\
	\leq \,&  \mu_2 + 2\sqrt{\frac{2\sigma^2}{n}\log\biggl(\frac{2d}{\delta}\biggr)} .
\end{align*}
	
\smallskip

\noindent \textbf{Final Bound:} Putting things together, the type II error of the test is bounded above by
\begin{align*}
	& P \Big( \sqrt{n}\bm{\gamma}_{\widehat{s}}^\top \overline{\bX}^{(1)} \leq z_{1-\alpha} \sqrt{\bm{\gamma}_{\widehat{s}}^\top \widehat{\bm{\Sigma}}^{(1)} \bm{\gamma}_{\widehat{s}}} \Big) \leq \mathrm{(I)} + \mathrm{(II)} + 2\delta \\[.5em]
	\leq \, & P\Biggl( \sqrt{n}(\mu_1 - \mu_{2}) \leq (z_{1-\alpha} + 1) \sqrt{4\sigma^2\delta^{-1}} + \sqrt{8\sigma^2\log\biggl(\frac{2|\mathbb{C}\cup \{2\}|}{\delta}\biggr)}\Biggr) \\[.5em]
	& + P \biggl( \frac{\sqrt{n}}{2}(\mu_1 - \mu_2) \leq (z_{1-\alpha} + 1) \sqrt{4\sigma^2\delta^{-1}} \biggr) 
	+ P \Biggl(\sqrt{\frac{8\sigma^2}{n}\log\biggl(\frac{2d}{\delta}\biggr)} > C_n \sqrt{\frac{\log (d)}{n}}\Biggr)  + 4\delta.
\end{align*}
Recall that $\mu_1 - \mu_2 \geq C_n' \sqrt{n^{-1}(1 \vee \log|\mathbb{C}|)}$ for some positive sequence $C_n'$ diverging to infinity. Consequently, each of the terms above approaches zero uniformly over $P \in \mathcal{P}_{1,r}(\varepsilon;\tau)$ as $n \to \infty$, provided that $\sqrt{\delta} C_n' \to \infty$ and $C_n/\sqrt{\log(1/\delta)} \to \infty$. For instance, choosing $\delta = 1/2 \wedge (C_n'^{-1} \vee e^{-C_n})$ suffices to ensure these conditions. This completes the proof of \Cref{Thm: upper bound} with $\widehat{s}_{\mathrm{plug}}$.

\subsubsection{Proof for noise-adjusted estimator $\widehat{s}_{\mathrm{adj}}$} \label{App: adj}

We next prove \Cref{Thm: upper bound} by considering the DA argmin test using the noise-adjusted estimator $\widehat{s} = \widehat{s}_{\mathrm{adj}}$. The proof remains the same as that for the plug-in estimator $\widehat{s} = \widehat{s}_{\mathrm{plug}}$ until the point where we define the terms $\mathrm{(I)}$ and $\mathrm{(II)}$. It therefore suffices to show that both terms vanish under the conditions stated in the theorem. The main challenge lies in the fact that $\widehat{s}_{\mathrm{adj}}$ does not directly target $s = \sargmin_{2 \leq k \leq d} \mu_k$, as it incorporates variance estimators into the objective function. To address this, we carefully relate $\widehat{s}_{\mathrm{adj}}$ to $\widehat{s}_{\mathrm{plug}}$ and build on the earlier analysis for the plug-in estimator. Throughout the proof, we denote $\widehat{s} = \widehat{s}_{\mathrm{adj}}$ to simplify the notation. 

\medskip 
\noindent \textbf{Term (I):} We begin with the first term $\mathrm{(I)}$, which is recalled as
\begin{align*}
	\mathrm{(I)} = P \Bigl( \sqrt{n}(\mu_1 - \mu_{\widehat{s}}) \leq (z_{1-\alpha} + 1) \sqrt{4\sigma^2\delta^{-1}} \,\cap\, \widehat{s} \in  \mathbb{C} \Bigr).
\end{align*}
Define the event $\widetilde{\mathcal{E}}_{3,\delta}$ as
\begin{align*}
	\widetilde{\mathcal{E}}_{3,\delta} \coloneqq \bigcap_{k \in \mathbb{C} \cup \{2\}}\biggl\{\big| \overline{X}^{(2)}_k - \overline{X}^{(2)}_1 - \mu_k + \mu_1 \big| < \sqrt{\frac{8\sigma^2}{n}\log\biggl(\frac{2|\mathbb{C} \cup \{2\}|}{\delta}\biggr)} \biggr\}.
\end{align*}
Following the same argument as before, we can show that $\widetilde{\mathcal{E}}_{3,\delta}$ holds with probability at least $1 - \delta$, using the sub-Gaussian tail bound, the union bound, and the fact that the sum of two sub-Gaussian random variables with variance proxy $\sigma^2$ is also sub-Gaussian with variance proxy $4\sigma^2$.

For brevity, define $\Delta_{\delta, \mathbb{C}} \coloneqq \sqrt{8\sigma^2\log\bigl(2|\mathbb{C} \cup \{2\}|/\delta\bigr)}$. Under the event $\widetilde{\mathcal{E}}_{3,\delta} \cap \{\widehat{s} \in \mathbb{C}\}$, we then obtain the following inequalities:
\begin{align*}
	\mu_1 - \mu_{\widehat{s}} \geq \, & \overline{X}^{(2)}_1 - \overline{X}^{(2)}_{\widehat{s}} - \frac{\Delta_{\delta, \mathbb{C}}}{\sqrt{n}} \\[.5em]
	= \, & \frac{\overline{X}^{(2)}_1 - \overline{X}^{(2)}_{\widehat{s}}}{\sqrt{\bm{\gamma}_{\widehat{s}}^\top \widehat{\bm{\Sigma}}^{(2)} \bm{\gamma}_{\widehat{s}}}\vee \kappa} \cdot \sqrt{\bm{\gamma}_{\widehat{s}}^\top \widehat{\bm{\Sigma}}^{(2)} \bm{\gamma}_{\widehat{s}}}\vee \kappa - \frac{\Delta_{\delta, \mathbb{C}}}{\sqrt{n}} \\[.5em]
	\overset{(\star)}{\geq} \, & \frac{\overline{X}^{(2)}_1 - \overline{X}^{(2)}_{2}}{\sqrt{\bm{\gamma}_{2}^\top \widehat{\bm{\Sigma}}^{(2)} \bm{\gamma}_{2}}\vee \kappa} \cdot \sqrt{\bm{\gamma}_{\widehat{s}}^\top \widehat{\bm{\Sigma}}^{(2)} \bm{\gamma}_{\widehat{s}}}\vee \kappa - \frac{\Delta_{\delta, \mathbb{C}}}{\sqrt{n}} \\[.5em]
	\geq \, & \frac{\sqrt{\bm{\gamma}_{\widehat{s}}^\top \widehat{\bm{\Sigma}}^{(2)} \bm{\gamma}_{\widehat{s}}}\vee \kappa}{\sqrt{\bm{\gamma}_{2}^\top \widehat{\bm{\Sigma}}^{(2)} \bm{\gamma}_{2}}\vee \kappa} \Biggl( \mu_1 - \mu_2 - \frac{\Delta_{\delta, \mathbb{C}}}{\sqrt{n}} \Biggr) - \frac{\Delta_{\delta, \mathbb{C}}}{\sqrt{n}},
\end{align*}
where step $(\star)$ uses the definition of $\widehat{s}$. Hence, by replacing $\mu_1 - \mu_{\widehat{s}}$ in $\mathrm{(I)}$ with the established lower bound, it holds that 
\begin{align*}
	\mathrm{(I)} \leq \, & P\Biggl( \frac{\sqrt{\bm{\gamma}_{\widehat{s}}^\top \widehat{\bm{\Sigma}}^{(2)} \bm{\gamma}_{\widehat{s}}}\vee \kappa}{\sqrt{\bm{\gamma}_{2}^\top \widehat{\bm{\Sigma}}^{(2)} \bm{\gamma}_{2}}\vee \kappa} \Biggl\{ \sqrt{n}(\mu_1 - \mu_2) - \Delta_{\delta, \mathbb{C}} \Biggr\} - \Delta_{\delta, \mathbb{C}} \\[.5em] 
    & \quad \quad \quad \quad \quad \quad \leq (z_{1-\alpha} + 1) \sqrt{4\sigma^2\delta^{-1}} \,\cap\, \widetilde{\mathcal{E}}_{3,\delta} \,\cap\, \{\widehat{s} \in  \mathbb{C}\} \Biggr) +  P(\widetilde{\mathcal{E}}_{3,\delta}^c) \\[.5em]
    = \, & P \Biggl( \sqrt{n}(\mu_1 - \mu_2) \leq \Biggl\{1 + \frac{\sqrt{\bm{\gamma}_{2}^\top \widehat{\bm{\Sigma}}^{(2)} \bm{\gamma}_{2}}\vee \kappa}{\sqrt{\bm{\gamma}_{\widehat{s}}^\top \widehat{\bm{\Sigma}}^{(2)} \bm{\gamma}_{\widehat{s}}}\vee \kappa} \Biggr\} \Delta_{\delta,\mathbb{C}} \\[.5em]
	& \quad \quad +   \frac{\sqrt{\bm{\gamma}_{2}^\top \widehat{\bm{\Sigma}}^{(2)} \bm{\gamma}_{2}}\vee \kappa}{\sqrt{\bm{\gamma}_{\widehat{s}}^\top \widehat{\bm{\Sigma}}^{(2)} \bm{\gamma}_{\widehat{s}}}\vee \kappa} (z_{1-\alpha} + 1) \sqrt{4\sigma^2\delta^{-1}} \,\cap\, \widetilde{\mathcal{E}}_{3,\delta} \,\cap\, \{\widehat{s} \in \mathbb{C}\}\Biggr) + P(\widetilde{\mathcal{E}}_{3,\delta}^c) \\[.5em]
	\leq \, & P \Biggl( \sqrt{n}(\mu_1 - \mu_2) \leq \Bigl(2 + \kappa^{-1}\sqrt{\bm{\gamma}_{2}^\top \widehat{\bm{\Sigma}}^{(2)} \bm{\gamma}_{2}} \Bigr) \Delta_{\delta,\mathbb{C}} \\[.5em]
	& \quad \quad +   \Bigl(1 + \kappa^{-1}\sqrt{\bm{\gamma}_{2}^\top \widehat{\bm{\Sigma}}^{(2)} \bm{\gamma}_{2}} \Bigr)(z_{1-\alpha} + 1) \sqrt{4\sigma^2\delta^{-1}} \,\cap\, \widetilde{\mathcal{E}}_{3,\delta} \,\cap\, \{\widehat{s} \in \mathbb{C}\}\Biggr) + \delta,
\end{align*}
where the last inequlity uses $(p \vee r) / (q \vee r) \leq 1 + r^{-1}p$ for any $p,q \geq 0$ and $r >0$. Moreover, we define another event 
\begin{align*}
	\widetilde{\mathcal{E}}_{1,\delta} \coloneqq \bigg\{ \bm{\gamma}_{2}^\top \widehat{\bm{\Sigma}}^{(2)} \bm{\gamma}_{2} \leq \frac{4\sigma^2}{\delta} \bigg\},
\end{align*}
which holds with probability at least $1 - \delta$, similarly to $\mathcal{E}_{1,\delta}$. By incorportaing this event into the above inequality for $\mathrm{(I)}$ using the union bound, we have 
\begin{align*}
	\mathrm{(I)} \leq P \Bigl( \sqrt{n}(\mu_1 - \mu_2) \leq \bigl(2 + \kappa^{-1}\sqrt{4\sigma^2 \delta^{-1}}  \bigr)\Delta_{\delta,\mathbb{C}} + (z_{1-\alpha}+1)\bigl(\sqrt{4\sigma^2 \delta^{-1}} + 4\kappa^{-1}\sigma^2 \delta^{-1} \bigr)\Bigr) + 2\delta.
\end{align*}
The above upper bound vanishes under the condition on $\mu_1 - \mu_2 \geq C_n' \varepsilon^\star$, provided that $\delta$ decreases sufficiently slowly. For instance, one can take $\delta = 1/2 \wedge C_n'^{-1}$. Hence the term $\mathrm{(I)}$ vanishes under the conditions stated in the theorem.

\medskip 
\noindent \textbf{Term (II):} For the second term $\mathrm{(II)}$, it suffices to bound $P(\widehat{s} \in \mathbb{C}^c_b)$ as in the earlier analysis for the plug-in approach. Define the event
\begin{align*}
	\widetilde{\mathcal{E}}_{4,\delta} \coloneqq \bigcap_{k = 2}^d \biggl\{\big| \overline{X}^{(2)}_k - \overline{X}^{(2)}_1 - \mu_k + \mu_1 \big| < \sqrt{\frac{8\sigma^2}{n}\log\biggl(\frac{2d}{\delta}\biggr)} \biggr\},
\end{align*}
which satisfies $P(\widetilde{\mathcal{E}}_{4,\delta}^c) \leq \delta$, analogous to previous arguments. Then, it holds that
\begin{align*}
	P(\widehat{s} \in \mathbb{C}^c_b) \leq P \biggl(\mu_{\widehat{s}} - \mu_2 > C_n \sigma \sqrt{\frac{\log (d)}{n}} \,\cap\, \widetilde{\mathcal{E}}_{4,\delta} \biggr) + \delta.
\end{align*}
Unlike $\widehat{s}_{\mathrm{plug}}$, we cannot directly relate $\mu_{\widehat{s}}$ to $\mu_s$; so a more involved argument is required to formally show that the above upper bound vanishes. To this end, let $\Delta_{\delta,d} \coloneqq \sqrt{8\sigma^2\log(2d/\delta)}$ for brevity. Under the event $\widetilde{\mathcal{E}}_{4,\delta}$, we have 
\begin{align*}
	\mu_{\widehat{s}} - \mu_1 + \mu_1 - \mu_2 \leq \, & \overline{X}_{\widehat{s}}^{(2)} - \overline{X}_1^{(2)} + \mu_1 - \mu_2 + n^{-1/2}\Delta_{\delta,d} \\
	= \, & \sqrt{\bm{\gamma}_{\widehat{s}}^\top \widehat{\bm{\Sigma}}^{(2)} \bm{\gamma}_{\widehat{s}}}\vee \kappa \times \frac{\overline{X}_{\widehat{s}}^{(2)} - \overline{X}_1^{(2)}}{\sqrt{\bm{\gamma}_{\widehat{s}}^\top \widehat{\bm{\Sigma}}^{(2)} \bm{\gamma}_{\widehat{s}}}\vee \kappa} + \mu_1 - \mu_2 + n^{-1/2}\Delta_{\delta,d} \\
	\leq \, & \frac{\sqrt{\bm{\gamma}_{\widehat{s}}^\top \widehat{\bm{\Sigma}}^{(2)} \bm{\gamma}_{\widehat{s}}}\vee \kappa}{\sqrt{\bm{\gamma}_{\widehat{s}_{\mathrm{plug}}}^\top \widehat{\bm{\Sigma}}^{(2)} \bm{\gamma}_{\widehat{s}_{\mathrm{plug}}}}\vee \kappa} \times \bigl(\overline{X}_{\widehat{s}_{\mathrm{plug}}}^{(2)} - \overline{X}_1^{(2)}\bigr) + \mu_1 - \mu_2 + n^{-1/2}\Delta_{\delta,d},
\end{align*}
where the last inequality follows by definition of $\widehat s$. 
Now, again by the definition of $\widehat{s}_{\mathrm{plug}}$ and $\widehat{s} = \widehat{s}_{\mathrm{adj}}$, we make a key observation that
\begin{align*}
	\frac{\overline{X}_{\widehat{s}_{\mathrm{plug}}}^{(2)} - \overline{X}_1^{(2)}}{\sqrt{\bm{\gamma}_{\widehat{s}}^\top \widehat{\bm{\Sigma}}^{(2)} \bm{\gamma}_{\widehat{s}}}\vee \kappa} \overset{\mathrm{(i)}}{\leq} \frac{\overline{X}_{\widehat{s}}^{(2)} - \overline{X}_1^{(2)}}{\sqrt{\bm{\gamma}_{\widehat{s}}^\top \widehat{\bm{\Sigma}}^{(2)} \bm{\gamma}_{\widehat{s}}}\vee \kappa} \overset{\mathrm{(ii)}}{\leq} \frac{\overline{X}_{\widehat{s}_{\mathrm{plug}}}^{(2)} - \overline{X}_1^{(2)}}{\sqrt{\bm{\gamma}_{\widehat{s}_{\mathrm{plug}}}^\top \widehat{\bm{\Sigma}}^{(2)} \bm{\gamma}_{\widehat{s}_{\mathrm{plug}}}}\vee \kappa},
\end{align*}
where step $\mathrm{(i)}$ uses the definition of $\widehat{s}_{\mathrm{plug}}$ and step $\mathrm{(ii)}$ uses the definition of $\widehat{s}$. Combining the first and last expression, whenever the event $\widetilde{\mathcal{E}}_5 \coloneqq \bigl\{ \overline{X}_{\widehat{s}_{\mathrm{plug}}}^{(2)} - \overline{X}_1^{(2)} < 0 \bigr\}$ holds, it follows that
\begin{align*}
\frac{\sqrt{\bm{\gamma}_{\widehat{s}}^\top \widehat{\bm{\Sigma}}^{(2)} \bm{\gamma}_{\widehat{s}}}\vee \kappa}{\sqrt{\bm{\gamma}_{\widehat{s}_{\mathrm{plug}}}^\top \widehat{\bm{\Sigma}}^{(2)} \bm{\gamma}_{\widehat{s}_{\mathrm{plug}}}}\vee \kappa} \leq 1.
\end{align*}
Therefore, the probability can be bounded as follows:
\begin{align*}
	&P \biggl(\mu_{\widehat{s}} - \mu_2 > C_n \sigma \sqrt{\frac{\log (d)}{n}} \,\cap\, \widetilde{\mathcal{E}}_{4,\delta} \biggr) \\
	\leq \, & P \biggl( \bigl(\overline{X}_{\widehat{s}_{\mathrm{plug}}}^{(2)} - \overline{X}_1^{(2)}\bigr) + \mu_1 - \mu_2 + n^{-1/2}\Delta_{\delta,d}> C_n \sigma \sqrt{\frac{\log (d)}{n}}  \,\cap\, \widetilde{\mathcal{E}}_{4,\delta}  \,\cap\, \widetilde{\mathcal{E}}_{5} \biggr) + P\bigl(\widetilde{\mathcal{E}}_{5}^c\bigr) \\
	\leq \, & P \biggl( \mu_{\widehat{s}_{\mathrm{plug}}} - \mu_2 + 2n^{-1/2}\Delta_{\delta,d}> C_n \sigma \sqrt{\frac{\log (d)}{n}}  \,\cap\, \widetilde{\mathcal{E}}_{4,\delta}  \,\cap\, \widetilde{\mathcal{E}}_{5} \biggr) + P\bigl(\widetilde{\mathcal{E}}_{5}^c\bigr) \\
	\leq \, & P \biggl( \mu_{\widehat{s}_{\mathrm{plug}}} - \mu_2 + 2n^{-1/2}\Delta_{\delta,d}> C_n \sigma \sqrt{\frac{\log (d)}{n}}  \,\cap\, \widetilde{\mathcal{E}}_{4,\delta}  \,\cap\, \widetilde{\mathcal{E}}_{5} \biggr) + P\bigl(\widetilde{\mathcal{E}}_{5}^c\bigr).
\end{align*}
The first term above can be bounded by the same argument as in the proof of \Cref{Thm: upper bound} for the plug-in estimator. It remains to show that $P\bigl(\widetilde{\mathcal{E}}_{5}^c\bigr)$ vanishes. To this end, observe that 
\begin{align*}
	P\bigl(\widetilde{\mathcal{E}}_{5}^c\bigr) =\,& P\Bigl(\overline{X}_{\widehat{s}_{\mathrm{plug}}}^{(2)} - \overline{X}_1^{(2)} \geq 0 \Bigr) \leq P\Bigl(\overline{X}_{2}^{(2)} - \overline{X}_1^{(2)} \geq 0\Bigr) \\
	= \, & P\Bigl(\overline{X}_{2}^{(2)} - \overline{X}_1^{(2)} + \mu_1 - \mu_2 \geq \mu_1 - \mu_2 \Bigr) \\
	\leq \, & \frac{2\sigma^2}{n(\mu_1 - \mu_2)^2},
\end{align*}
where the first inequality uses the fact that $\overline{X}_{\widehat{s}_{\mathrm{plug}}}^{(2)}$ is the argmin index of the sample mean vectors $\overline{X}_2^{(2)}, \ldots, \overline{X}_d^{(2)}$ and the last inequality uses Chebyshev's inequality. Therefore, we have shown that the term $\mathrm{(II)}$ vanishes under the conditions stated in the theorem. Combining the bounds on terms $\mathrm{(I)}$ and $\mathrm{(II)}$ completes the proof of \Cref{Thm: upper bound} with $\widehat{s}_{\mathrm{adj}}$.

\subsection{Proof of \Cref{Cor: upper bound}} \label{App: cor upper bound}
By construction, the index $r$ is excluded from the confidence set if and only if the DA argmin test rejects the null hypothesis $H_0: r \in \Theta$. The result of \Cref{Cor: upper bound} then follows immediately from the power guarantee established in \Cref{Thm: upper bound}.

\subsection{Proof of \Cref{Thm: lower bound}} \label{App: lower bound}
We work with $n$ samples rather than $2n$ samples, which only affects a constant factor in the lower bound. Additionally, we explicitly indicate that the probability $P$ is taken over the i.i.d.~samples $\bX_1,\ldots,\bX_n$ by adding the superscript $\otimes n$ to $P$. As in the proof of \Cref{Thm: upper bound}, we set $r=1$ without loss of generality.

For $m \in \mathbb{Z}_{>0}$, the mean vector $\bm{\mu}^{(0)}$ consists of the first $m+1$ components set to zero, followed by the remaining $d-m-1$ components set to $b_n>0$, that is
\begin{align*}
	\bm{\mu}^{(0)} = (0,\underbrace{0,\ldots, 0}_{\text{$m$ entries}}, \underbrace{b_n,\ldots,b_n}_{\text{$d-m-1$ entries}})^\top \in \mathbb{R}^d.
\end{align*}
Here, $b_n$ is a positive sequence that varies with $n$ and will be specified later. Similarly, for each $i \in [m]$ and $\rho > 0$, the mean vector $\bm{\mu}^{(i)}$ is defined as
\begin{align*}
	\bm{\mu}^{(i)} = \bm{\mu}^{(0)} - \rho \cdot \boldsymbol{e}_{i+1} \in \mathbb{R}^d.
\end{align*}
In words, the mean vector $\bm{\mu}^{(i)}$ is obtained by decreasing the $(i+1)$-th component of $\bm{\mu}^{(0)}$ by $\rho$. For instance,
\begin{align*}
	\bm{\mu}^{(1)} = (0,-\rho,\underbrace{0,\ldots, 0}_{\text{$m-1$ entries}}, \underbrace{b_n,\ldots,b_n}_{\text{$d-m-1$ entries}})^\top \in \mathbb{R}^d.
\end{align*}

Let $P_i$ be the distribution of $N(\bm{\mu}^{(i)}, \sigma^2\bm{I}_d)$ for $i \in \{0,1,2,\ldots,m\}$, and let $P_i^{\otimes n}$ denote the $n$-fold product distribution of $P_i$. Define a mixture distribution of $P_1^{\otimes n},\ldots,P_m^{\otimes n}$ as
\begin{align*}
	P_{\mathrm{mix}}^{\otimes n} = \frac{1}{m} \sum_{i=1}^m P_i^{\otimes n}.
\end{align*}
Let $\phi(\bm{x}; \bm{\mu}, \sigma^2)$ be the density function of $N(\bm{\mu}, \sigma^2\bm{I}_d)$ evaluated at $\bm{x} \in \mathbb{R}^d$, 
and compute the chi-square divergence between $P_{\mathrm{mix}}^{\otimes n}$ and $P_0^{\otimes n}$ as
\begin{align*}
	\chi^2(P_{\mathrm{mix}}^{\otimes n} \| P_0^{\otimes n}) = \, & \mE_{P_0^{\otimes n}} \Biggl[ \biggl( \frac{dP_{\mathrm{mix}}^{\otimes n}}{dP_0^{\otimes n}}(\bX_1,\ldots,\bX_n) \biggr)^2\Biggr] - 1  \\[.5em]
	= \, & \mE_{P_0^{\otimes n}} \Biggl[ \Biggl( \frac{1}{m} \sum_{i=1}^m \prod_{j=1}^n \frac{\phi(\bX_j; \bm{\mu}^{(i)}, \sigma^2)}{\phi(\bX_j; \bm{\mu}^{(0)}, \sigma^2)} \Biggr)^2 \Biggr] - 1 \\[.5em]
	= \, & \frac{1}{m^2} \sum_{i=1}^m \sum_{j=1}^m \mE_{P_0^{\otimes n}} \Biggl[ \prod_{k=1}^n \frac{\phi(\bX_k; \bm{\mu}^{(i)}, \sigma^2)}{\phi(\bX_k; \bm{\mu}^{(0)}, \sigma^2)} \cdot \frac{\phi(\bX_k; \bm{\mu}^{(j)}, \sigma^2)}{\phi(\bX_k; \bm{\mu}^{(0)}, \sigma^2)} \Biggr] - 1 \\[.5em]
	= \, & \frac{1}{m^2} \sum_{i=1}^m \sum_{j=1}^m \Biggl( \mE_{P_0} \Biggl[ \frac{\phi(\bX; \bm{\mu}^{(i)}, \sigma^2)\phi(\bX; \bm{\mu}^{(j)}, \sigma^2)}{\phi(\bX; \bm{\mu}^{(0)}, \sigma^2)^2} \Biggr] \Biggr)^n - 1,
\end{align*}
where the last equality uses the fact that $\bX, \bX_1,\ldots, \bX_n$ are i.i.d.~samples from $P_0$. Focusing on the expectation inside, an explicit form is derived as
\begin{align*}
	\mE_{P_0} \Biggl[ \frac{\phi(\bX; \bm{\mu}^{(i)}, \sigma^2)\phi(\bX; \bm{\mu}^{(j)}, \sigma^2)}{\phi(\bX; \bm{\mu}^{(0)}, \sigma^2)^2} \Biggr] =\,& \exp\Bigl(\sigma^{-2} \bigl\langle \bm{\mu}^{(i)} - \bm{\mu}^{(0)}, \bm{\mu}^{(j)} - \bm{\mu}^{(0)} \bigr\rangle \Bigr) \\[.5em]
	=\,& \exp \Bigl( \sigma^{-2} \rho^2 \big\langle \boldsymbol{e}_{i+1}, \boldsymbol{e}_{j+1} \big\rangle \Bigr). 
\end{align*}
Returning to the chi-square divergence,
\begin{align*}
	\chi^2(P_{\mathrm{mix}}^{\otimes n} \| P_0^{\otimes n}) = \, & \frac{1}{m^2} \sum_{i=1}^m \sum_{j=1}^m \exp \Bigl( \sigma^{-2} \rho^2 \big\langle \boldsymbol{e}_{i+1}, \boldsymbol{e}_{j+1} \big\rangle \Bigr)^n - 1 \\[.5em]
	= \, & \frac{1}{m} \exp\bigl(n\sigma^{-2} \rho^2\bigr) - 1.
\end{align*}

We now set $\rho = \varepsilon$, $b_n > C_n \sigma \sqrt{n^{-1} \log(d)}$ and $m = \tau + 1$, which guarantees that each alternative distribution $P_i$ belongs to the class $\mathcal{P}_{1,r}(\varepsilon;\tau)$ as the mean difference satisfies $\mu_1 - \mu_i = \rho$ and $\mathbb{C}_1 = \{3, 4, \ldots, m+1\}$ yields cardinality $|\mathbb{C}_1| = m - 1 = \tau$.

With this setup, and denoting the total variation distance between $P$ and $Q$ as $\mathrm{TV}\bigl(P, Q\bigr)$, Ingster's $\chi^2$-method for minimax testing lower bounds~\citep{ingster1987minimax} yields that for sufficiently large $n$,
\begin{align*}
	\inf_{\psi \in \Psi_{\alpha}} \sup_{P \in \mathcal{P}_{1,r}(\varepsilon;\tau)} P^{\otimes n}(\psi = 0) \geq \, & \inf_{\psi \in \Psi_{\alpha}} P_{\mathrm{mix}}^{\otimes n}(\psi = 0) \geq  1 - \alpha - o(1) - \mathrm{TV}\bigl(P_{0}^{\otimes n}, P_{\mathrm{mix}}^{\otimes n}\bigr) \\[.5em]
	\geq \, & 1 - 2\alpha - \mathrm{TV}\bigl(P_{0}^{\otimes n}, P_{\mathrm{mix}}^{\otimes n}\bigr) \\[.5em]
	\geq \, &  1 - 2\alpha - \sqrt{\chi^2(P_{\mathrm{mix}}^{\otimes n} \| P_0^{\otimes n})},
\end{align*}
where the last inequality uses the inequality that $\mathrm{TV}(P,Q) \leq \sqrt{\chi^2(P\|Q)}$ for any two distributions $P$ and $Q$~\citep[][Section 2.4.1]{tsybakov2009}. Note that the little $o(1)$ term above is incorporated to account for the fact that $\psi$ is an asymptotically level-$\alpha$ test and $\alpha + o(1)$ is replaced by $2\alpha$ by taking $n$ sufficiently large.

Now, to ensure that the minimax type II error is at least $\beta$, we must have 
\begin{align*}
	1 - 2\alpha - \sqrt{\chi^2(P_{\mathrm{mix}}^{\otimes n} \| P_0^{\otimes n})} \geq \beta	\quad & \Longleftrightarrow \quad (1 - 2\alpha - \beta)^2 \geq \chi^2(P_{\mathrm{mix}}^{\otimes n} \| P_0^{\otimes n}) \\[.5em]
	&\Longleftrightarrow \quad \sqrt{\frac{\sigma^2}{n}\log \bigl(m(1-2\alpha-\beta)^2 + 1\bigr)} \geq \varepsilon.
\end{align*}
Moreover an algebraic argument shows that 
\begin{align*}
	\log \bigl(|\mathbb{C}_1|(1-2\alpha-\beta)^2 + (1-2\alpha-\beta)^2 + 1\bigr)
	\geq  \log\bigl(1 + (1-2\alpha-\beta)^2\bigr) \cdot \bigl(1 \vee \log |\mathbb{C}_1| \bigr).
\end{align*}
Hence a sufficient condition for the minimax type II error to be at least $\beta$ is
\begin{align*}
	\sqrt{\log\bigl(1 + (1-2\alpha-\beta)^2\bigr) \cdot \frac{\sigma^2}{n} \cdot \bigl(1 \vee \log |\mathbb{C}_1| \bigr)} \geq \varepsilon.
\end{align*}
Setting $c = \sqrt{\sigma^2\log\bigl(1 + (1-2\alpha-\beta)^2\bigr)}$ completes the proof of \Cref{Thm: lower bound}.

\subsection{Proof of \Cref{Cor: lower bound}} \label{App: lower bound cor}
Given $r \in [d]$, define 
\begin{align*}
	\mathcal{A}_{\alpha,r} \coloneqq \Bigl\{ \widehat{\Theta}: \liminf_{n \to \infty} \inf_{P \in \mathcal{P}_{0,r}}  P(r \in \widehat{\Theta}) \geq 1 - \alpha \Bigr\},
\end{align*}
which satisfies $\mathcal{A}_{\alpha} \subseteq \mathcal{A}_{\alpha,r}$. Now consider the test $\psi$ that rejects the null hypothesis $H_0: r \in \Theta$ if and only if $r \notin \widehat{\Theta}$. This establishes a one-to-one correspondence between $\mathcal{A}_{\alpha,r}$ and the set of asymptotic level-$\alpha$ tests $\Psi_{\alpha,r}$. Therefore it follows that
\begin{align*}
	\sup_{\widehat{\Theta} \in \mathcal{A}_{\alpha}} \inf_{P \in \mathcal{P}_{1,r}(\varepsilon;\tau)} P(r \notin \widehat{\Theta}) \leq \, & \sup_{\widehat{\Theta} \in \mathcal{A}_{\alpha,r}} \inf_{P \in \mathcal{P}_{1,r}(\varepsilon;\tau)} P(r \notin \widehat{\Theta})  \\
	= \, & 1 - \inf_{\widehat{\Theta} \in \mathcal{A}_{\alpha,r}} \sup_{P \in \mathcal{P}_{1,r}(\varepsilon;\tau)} P(r \in \widehat{\Theta}) \\
	= \, & 1 - \inf_{\psi \in \Psi_{\alpha,r}} \sup_{P \in \mathcal{P}_{1,r}(\varepsilon;\tau)} P(\psi = 0).
\end{align*}
Taking $\limsup_{n \to \infty}$ on both sides, the upper bound becomes $1-\beta$ by \Cref{Thm: lower bound}, which completes the proof.

\subsection{Proof of \Cref{Thm: robust}}
\label{subsec:proof-robust}

The proof of \Cref{Thm: robust} closely parallels that of \Cref{Thm: upper bound}, with the key distinction being the use of the MoM estimators in place of empirical means for estimating the argmin $s$. The core technical component in the proof of \Cref{Thm: upper bound} was a sub-Gaussian tail bound for the sample mean, which was used to establish high-probability bounds for the events $\mathcal{E}_{3,\delta}$, $\mathcal{E}_{4,\delta}$, $\widetilde{\mathcal{E}}_{3,\delta}$ and $\widetilde{\mathcal{E}}_{4,\delta}$. In the proof of \Cref{Thm: robust}, these events are defined analogously, with MoM estimators replacing sample means. Their associated probability bounds follow from the sub-Gaussian tail inequality for MoM estimators~\citep[e.g.,][Proposition 5]{hsu2016loss}, differing only in constant factors. In this setting, the parameter $\eta$ in the MoM framework serves the same role as $\delta$ in the proof of \Cref{Thm: upper bound}. The additional factor $e^{-n/18}$ in the definition of $\eta$ accounts for the constraint $V = 4.5\ceil{\log(1/\eta)} \leq n$, along with the condition $n \geq 18 \ceil{\log(1/\eta)}$. These choices follow the requirements in \citet[][Proposition 5]{hsu2016loss}. We omit further details, as the remainder of the argument proceeds almost identically to the proof of \Cref{Thm: upper bound} with only a minor modification.

\subsection{Proof of \Cref{Thm: uniform coverage}} \label{App: uniform coverage}

We start by observing that by the union bound,
\begin{align*}
	P(\Theta \subseteq \widehat{\Theta}_{\mathrm{DA}}^{\mathrm{uni}}) =   1 - P\bigl( \cup_{r \in \Theta} \bigl\{ r \notin \widehat{\Theta}_{\mathrm{DA}}^{\mathrm{uni}} \bigr\} \bigr) \geq  1 - \sum_{r \in \Theta} P(r \notin \widehat{\Theta}_{\mathrm{DA}}^{\mathrm{uni}}).
\end{align*}
By the (conditional) Berry--Esseen bound for the studentized mean~\citep[][Theorem 1.2]{bentkus1996berry}, we have
\begin{align*}
	P(r \notin \widehat{\Theta}_{\mathrm{DA}}^{\mathrm{uni}} \given D_2) \leq \frac{\alpha}{1 \vee |\widehat{\Theta}^{(2)}|} + \frac{C'}{\sqrt{n}},
\end{align*}
where $C'>0$ is a universal constant. Plugging this into the earlier expression and taking expectations with respect to $D_2$, we obtain
\begin{align*}
	P(\Theta \subseteq \widehat{\Theta}_{\mathrm{DA}}^{\mathrm{uni}}) \geq 1 - \mE_P\biggl[\frac{|\Theta|}{1 \vee |\widehat{\Theta}^{(2)}|}\biggr] \alpha - \frac{|\Theta|C'}{\sqrt{n}}.
\end{align*}
The last term is negligible under the assumption that $\sup_{P \in \mathcal{P}^{\leq 3}} |\Theta(P)| = o(\sqrt{n})$. To show that the first two terms are asymptotically lower bounded by $1-\alpha$, we only need to show that 
\begin{align*}
	\limsup_{n \to \infty}\sup_{P \in \mathcal{P}^{\leq 3}}\mE_P\biggl[\frac{|\Theta|}{1 \vee |\widehat{\Theta}^{(2)}|} - 1 \biggr] \leq 0. 
\end{align*}
Note that $|\widehat{\Theta}^{(2)}| = \sum_{i \in [d]} \mathds{1}(\widetilde{\psi}_i = 0)$, where $\widetilde{\psi}_i = \psi_i(D_2,n^{-1/2})$ is the DA argmin test applied to the second half of the data at level $n^{-1/2}$. By the definition of $\widehat{\Theta}^{(2)}$, we have
\begin{align*}
	\frac{|\Theta|}{1 \vee |\widehat{\Theta}^{(2)}|} \leq \frac{|\Theta|}{1 \vee \sum_{r \in \Theta}\mathds{1}(\widetilde{\psi}_r = 0)}.
\end{align*}
For a positive $\epsilon >0$ specified later, define the event
\begin{align*}
	\mathcal{B} \coloneqq \biggl\{|\Theta|^{-1}\sum_{r \in \Theta}\mathds{1}(\widetilde{\psi}_r = 0) \leq 1 - \epsilon \biggr\}.
\end{align*}
By Markov's inequality and the Berry--Esseen bound, we have 
\begin{align*}
	P(\mathcal{B}) \leq \frac{1}{\epsilon} -  \frac{1}{\epsilon |\Theta|} \sum_{r \in \Theta} P(\widetilde{\psi}_r = 0) \leq \frac{1}{\epsilon} -  \frac{1}{\epsilon |\Theta|} \sum_{r \in \Theta} \biggl(1 - \frac{1}{\sqrt{n}} - \frac{C'}{\sqrt{n}}\biggr) = \frac{1+C'}{\epsilon \sqrt{n}}.
\end{align*}
Using the preliminary results, we can bound the expectation as follows:
\begin{align*}
	\mE_P\biggl[\frac{|\Theta|}{1 \vee |\widehat{\Theta}^{(2)}|} - 1\biggr] \leq \, & \mE_P\biggl[\biggl\{\frac{|\Theta|}{1 \vee \sum_{r \in \Theta}\mathds{1}(\widetilde{\psi}_r = 0)} - 1\biggr\} \mathds{1}(\mathcal{B})\biggr] \\
		+& \mE_P\biggl[\biggl\{\frac{|\Theta|}{1 \vee \sum_{r \in \Theta}\mathds{1}(\widetilde{\psi}_r = 0)} - 1 \biggr\} \mathds{1}(\mathcal{B}^c)\biggr] \\
		\leq \, & |\Theta| P(\mathcal{B}) + \frac{\epsilon}{1-\epsilon} \leq |\Theta| \frac{1+C'}{\epsilon \sqrt{n}} + \frac{\epsilon}{1-\epsilon} \\
		\leq \, & C'' \frac{|\Theta|^{1/2}}{n^{1/4}},
\end{align*}
where the last inequality holds by taking $\epsilon = |\Theta|^{1/2}/n^{1/4}$ for sufficiently large $n$. The upper bound is negligible under the assumption that $\sup_{P \in \mathcal{P}^{\leq 3}} |\Theta(P)| = o(\sqrt{n})$. This completes the proof of \Cref{Thm: uniform coverage}.

\subsection{Proof of \Cref{Thm: CS for the smallest mean}} 
Fix any $P\in\mathcal{P}^{\le 3}$. By definition of the interval $\mathcal{C}_2$, we have
\begin{align*}
	P( \mu_{\star} \in \mathcal{C}_2) = \, &  P\biggl( \min_{k \in \widehat{\Theta}_{\mathrm{DA}}^{\mathrm{uni}}} \biggl\{ \overline{X}_k^{(1)} - z_{1 - \frac{\alpha}{2\widehat{d}}} \frac{\widehat{\sigma}_k^{(1)}}{\sqrt{n}} \biggr\} \leq \mu_\star \leq \min_{k \in \widehat{\Theta}_{\mathrm{DA}}^{\mathrm{uni}}} \biggl\{\overline{X}_k^{(1)} + z_{1 - \frac{\alpha}{2\widehat{d}}} \frac{\widehat{\sigma}_k^{(1)}}{\sqrt{n}} \biggr\} \biggr) \\
	\geq \, & P\biggl( \min_{k \in \widehat{\Theta}_{\mathrm{DA}}^{\mathrm{uni}}} \biggl\{ \overline{X}_k^{(1)} - z_{1 - \frac{\alpha}{2\widehat{d}}} \frac{\widehat{\sigma}_k^{(1)}}{\sqrt{n}} \biggr\} \leq \mu_\star \leq \min_{k \in \widehat{\Theta}_{\mathrm{DA}}^{\mathrm{uni}}} \biggl\{\overline{X}_k^{(1)} + z_{1 - \frac{\alpha}{2\widehat{d}}} \frac{\widehat{\sigma}_k^{(1)}}{\sqrt{n}} \biggr\}, \, \Theta \subseteq \widehat{\Theta}_{\mathrm{DA}}^{\mathrm{uni}} \biggr) \\
	= \, &  P\biggl( \min_{k \in \widehat{\Theta}_{\mathrm{DA}}^{\mathrm{uni}}} \biggl\{ \overline{X}_k^{(1)} - z_{1 - \frac{\alpha}{2\widehat{d}}} \frac{\widehat{\sigma}_k^{(1)}}{\sqrt{n}} \biggr\} \leq \min_{k \in \widehat{\Theta}_{\mathrm{DA}}^{\mathrm{uni}}} \mu_k \leq \min_{k \in \widehat{\Theta}_{\mathrm{DA}}^{\mathrm{uni}}} \biggl\{\overline{X}_k^{(1)} + z_{1 - \frac{\alpha}{2\widehat{d}}} \frac{\widehat{\sigma}_k^{(1)}}{\sqrt{n}} \biggr\}, \,  \Theta \subseteq \widehat{\Theta}_{\mathrm{DA}}^{\mathrm{uni}} \biggr) \\
	\geq \, & P\biggl( \forall k \in \widehat{\Theta}_{\mathrm{DA}}^{\mathrm{uni}}: \overline{X}_k^{(1)} - z_{1 - \frac{\alpha}{2\widehat{d}}} \frac{\widehat{\sigma}_k^{(1)}}{\sqrt{n}} \leq  \mu_k \leq \overline{X}_k^{(1)} + z_{1 - \frac{\alpha}{2\widehat{d}}} \frac{\widehat{\sigma}_k^{(1)}}{\sqrt{n}}, \,  \Theta \subseteq \widehat{\Theta}_{\mathrm{DA}}^{\mathrm{uni}} \biggr) \\
	\geq \, & 1 - P\left(\bigcup_{k \in \widehat{\Theta}_{\mathrm{DA}}^{\mathrm{uni}}}\Biggl\{\frac{\sqrt{n}|\overline{X}_k^{(1)}- \mu_k|}{\widehat{\sigma}_k^{(1)}} > z_{1 - \frac{\alpha}{2\widehat{d}}}\Biggr\}\right) - P\bigl(\Theta \not\subseteq \widehat{\Theta}_{\mathrm{DA}}^{\mathrm{uni}} \bigr) \\
	\geq \, & 1 - \mE_P \left[ \sum_{k \in \widehat{\Theta}_{\mathrm{DA}}^{\mathrm{uni}}} P\left(\frac{\sqrt{n}|\overline{X}_k^{(1)}- \mu_k|}{\widehat{\sigma}_k^{(1)}} > z_{1 - \frac{\alpha}{2\widehat{d}}}\, \Bigg| \, D_2 \right) \right] - P\bigl(\Theta \not\subseteq \widehat{\Theta}_{\mathrm{DA}}^{\mathrm{uni}} \bigr),
\end{align*}
where the second equality uses $\min_{k \in \widehat{\Theta}_{\mathrm{DA}}^{\mathrm{uni}}} \mu_k = \mu_\star$ whenever $\Theta \subseteq \widehat{\Theta}_{\mathrm{DA}}^{\mathrm{uni}}$, and the last inequality uses the (conditional) union bound.

\Cref{Thm: uniform coverage} guarantees that the last term $P(\Theta \not\subseteq \widehat{\Theta}_{\mathrm{DA}}^{\mathrm{uni}})$ tends to zero. For the remaining expectation, the (conditional) Berry--Esseen bound for the studentized mean~\citep[][Theorem 1.2]{bentkus1996berry} and the moment condition defining $\mathcal{P}^{\leq 3}$ give 
\begin{align*}
	P\left(\frac{\sqrt{n}|\overline{X}_k^{(1)}- \mu_k|}{\widehat{\sigma}_k^{(1)}} > z_{1 - \frac{\alpha}{2\widehat{d}}}\, \Bigg| \, D_2 \right) \leq \frac{\alpha}{\widehat{d}} + \frac{C'}{\sqrt{n}}
\end{align*}
for some universal constant $C'>0$. Consequently,
\begin{align*}
	\mE_P \left[ \sum_{k \in \widehat{\Theta}_{\mathrm{DA}}^{\mathrm{uni}}} P\left(\frac{\sqrt{n}|\overline{X}_k^{(1)}- \mu_k|}{\widehat{\sigma}_k^{(1)}} > z_{1 - \frac{\alpha}{2\widehat{d}}}\, \Bigg| \, D_2 \right)\right] \leq \alpha + \frac{C'\mE_P[\widehat{d}]}{\sqrt{n}} \leq \alpha + o(1),
\end{align*}
where the last inequality uses
$\sup_{P\in\mathcal{P}^{\le 3}}\mE_P[\widehat d]=o(n^{1/2})$.
Combining the bounds yields
\begin{align*}
	P(\mu_\star\in \mathcal{C}_2)\ge 1-\alpha+o(1),
\end{align*}
which completes the proof.

\subsection{Equivalence conditions for central limit theorem}
The following lemma establishes the equivalence between the truncated second moment condition in \eqref{Eq: moment condition} and Lindeberg's condition for the central limit theorem. 

\begin{lemma} \label{lem: UCLT}
	Let $\mathcal{P}$ be a class of distributions on $\mathbb{R}$ and assume that each $X_P \sim P \in \mathcal{P}$ has mean zero and variance $\sigma_P^2$. The following two conditions are equivalent:
	\begin{itemize}
		\item Condition (A). $\displaystyle \lim_{\lambda \to \infty} \sup_{P \in \mathcal{P}}\mE_P\biggl[\frac{X_P^2}{\sigma_P^2} \mathds{1}(|X_P| > \lambda \sigma_P) \biggr] = 0$;
		\item Condition (B). $\displaystyle \lim_{\lambda \to \infty}\sup_{P \in \mathcal{P}}\mE_P\biggl[ \frac{X_P^2}{\sigma_P^2} \biggl(1 \wedge \frac{|X_P|}{\lambda\sigma_P} \biggr)\biggr] = 0$.
	\end{itemize}
	\begin{proof}
		We first show that \emph{(A)} implies \emph{(B)}. To establish this implication, we begin by proving the following inequality, which holds for any $\epsilon, \lambda > 0$:
		\begin{align*}
			\mE_P\biggl[ \frac{X_P^2}{\sigma_P^2} \biggl(1 \wedge \frac{|X_P|}{\lambda\sigma_P} \biggr)\biggr] \leq \epsilon + \mE_P\biggl[\frac{X_P^2}{\sigma_P^2} \mathds{1}(|X_P| > \lambda \epsilon \sigma_P) \biggr].
		\end{align*}
		To this end, we first observe a basic inequality, which holds for any $y \geq 0$ and $\epsilon > 0$:
		\begin{align} \label{eq: min bound}
			1 \wedge y \leq \epsilon + \mathds{1}(y > \epsilon).
		\end{align}
		This inequality follows by noting that $1 \wedge y \leq \epsilon$ when $y \leq \epsilon$, and $1 \wedge y \leq 1 \leq \epsilon + 1$ when $y > \epsilon$. Applying the inequality~\eqref{eq: min bound} to the random variable $y = |X_P|/(\lambda\sigma_P)$, we obtain
		\begin{align*}
			\mE_P\biggl[ \frac{X_P^2}{\sigma_P^2} \biggl(1 \wedge \frac{|X_P|}{\lambda\sigma_P} \biggr)\biggr] \leq \,& \epsilon \mE_P\biggl[ \frac{X_P^2}{\sigma_P^2} \biggr] + \mE_P\biggl[\frac{X_P^2}{\sigma_P^2} \mathds{1}(|X_P| > \lambda \epsilon \sigma_P) \biggr] \\
			= \, & \epsilon + \mE_P\biggl[\frac{X_P^2}{\sigma_P^2} \mathds{1}(|X_P| > \lambda \epsilon \sigma_P) \biggr].
		\end{align*}
		We take the supremum over $P \in \mathcal{P}$ on both sides of the above inequality, which gives
		\begin{align*}
			\sup_{P \in \mathcal{P}}\mE_P\biggl[ \frac{X_P^2}{\sigma_P^2} \biggl(1 \wedge \frac{|X_P|}{\lambda\sigma_P} \biggr)\biggr] \leq \epsilon + \sup_{P \in \mathcal{P}}\mE_P\biggl[\frac{X_P^2}{\sigma_P^2} \mathds{1}(|X_P| > \lambda \epsilon \sigma_P) \biggr].
		\end{align*}
		Taking the limit $\lambda \to \infty$ followed by $\epsilon \to 0$ yields the conclusion that \emph{(A)} implies \emph{(B)}.

		To prove the converse, observe that for any $y \geq 0$, we have $\mathds{1}(y \geq 1) \leq 1 \wedge y$. Applying this inequality to $y = |X_P| / (\lambda\sigma_P)$ gives
		\begin{align*}
			\sup_{P \in \mathcal{P}}\mE\biggl[ \frac{X_P^2}{\sigma_P^2} \mathds{1}(|X_P| > \lambda \sigma_P) \biggr] \leq \sup_{P \in \mathcal{P}}\mE_P\biggl[ \frac{X_P^2}{\sigma_P^2} \biggl(1 \wedge \frac{|X_P|}{\lambda\sigma_P} \biggr)\biggr].
		\end{align*}
		The second implication then follows by taking the limit as $\lambda \to \infty$. This completes the proof of the lemma.
	\end{proof}	
\end{lemma}

\subsection{Details on the validity of the MCS procedure by \citet{hansen2011model}} \label{App: MCS validity}
In this subsection, we provide a detailed proof of \citet[][Theorem 1(i)]{hansen2011model}. The original argument is somewhat abstract, so we include a complete and explicit proof for clarity. We gratefully acknowledge that the proof described below is due to Jing Lei.

Using their notation, for each such fixed $\mathcal{M}$, \citet{hansen2011model} require that the associated test and elimination pair $(\delta_{\mathcal{M}},e_{\mathcal{M}})$ satisfy
(a) $\limsup_{n\to\infty} P\bigl(\delta_{\mathcal{M}}=1 \mid H_{0,\mathcal{M}}\bigr)\le\alpha$;
(b) $\lim_{n\to\infty} P\bigl(\delta_{\mathcal{M}}=1 \mid H_{A,\mathcal{M}}\bigr)=1$; and
(c) $\lim_{n\to\infty} P\bigl(e_{\mathcal{M}}\in\mathcal{M}^* \mid H_{A,\mathcal{M}}\bigr)=0$, where $H_{0,\mathcal M}$  is the null hypothesis that $\mathcal M$ is optimal, and $H_{A,\mathcal M}$ is its complement. Under these conditions, they show that their MCS, denoted as $\widehat{\mathcal{M}}_{1-\alpha}^\ast$, is asymptotically valid, i.e.,
\begin{align*}
	\liminf_{n \to \infty} P \bigl(\mathcal{M}^\ast \subset \widehat{\mathcal{M}}_{1-\alpha}^\ast \bigr) \geq 1 - \alpha \quad \Longleftrightarrow \quad \limsup_{n \to \infty} P \bigl(\mathcal{M}^\ast \not\subset \widehat{\mathcal{M}}_{1-\alpha}^\ast \bigr) \leq \alpha.
\end{align*}
To prove this, define the event $\mathcal{E}$ that a good model in $\mathcal{M}^\ast$ is eliminated before the model sequence reaches the optimal model subset $\mathcal{M}^\ast$. Then 
\begin{align*}
	P \bigl(\mathcal{M}^\ast \not\subset \widehat{\mathcal{M}}_{1-\alpha}^\ast \bigr) & = P \bigl(\mathcal{M}^\ast \not\subset \widehat{\mathcal{M}}_{1-\alpha}^\ast, \, \mathcal{E} \bigr) + P \bigl(\mathcal{M}^\ast \not\subset \widehat{\mathcal{M}}_{1-\alpha}^\ast, \, \mathcal{E}^c \bigr) \\[.5em]
	& \leq P\Biggl( \bigcup_{i^\ast \in \mathcal{M}^\ast} \bigcup_{\substack{\mathcal{M} \text{ such that} \\ H_{A, \mathcal{M}} \text{ holds} }} \Bigl\{ e_{\mathcal{M}} = i^\ast \Bigr\} \Biggr) + P(\delta_{\mathcal{M}^\ast} = 1 \given H_{0,\mathcal{M}^\ast}) \\[.5em]
	& \leq \sum_{i^\ast \in \mathcal{M}^\ast} \sum_{\substack{\mathcal{M} \text{ such that} \\ H_{A, \mathcal{M}} \text{ holds} }} P\bigl( e_{\mathcal{M}} = i^\ast \given H_{A,\mathcal{M}}\bigr) + P(\delta_{\mathcal{M}^\ast} = 1 \given H_{0,\mathcal{M}^\ast}),
\end{align*}
where the last inequality follows from the union bound. Under the conditions on $e_{\mathcal{M}}$ and $\delta_{\mathcal{M}}$, the last term can be bounded above as follows: 
\begin{align*}
	\sum_{i^\ast \in \mathcal{M}^\ast} \sum_{\substack{\mathcal{M} \text{ such that} \\ H_{A, \mathcal{M}} \text{ holds} }} P\bigl( e_{\mathcal{M}} = i^\ast \given H_{A,\mathcal{M}}\bigr) \leq |\mathcal{M}^\ast| \times  2^{|\mathcal{M}_0|} \times o(1) \quad \text{and} \quad P(\delta_{\mathcal{M}^\ast} = 1 \given H_{0,\mathcal{M}^\ast}) \leq \alpha + o(1).
\end{align*}
Therefore, as long as $|\mathcal{M}_0|$ remains constant, it holds that $\limsup_{n \to \infty} P \bigl(\mathcal{M}^\ast \not\subset \widehat{\mathcal{M}}_{1-\alpha}^\ast \bigr) \leq \alpha$ as desired. However, when the size of the model space $|\mathcal{M}_0|$ increases with $n$, stronger uniform conditions on $e_{\mathcal{M}}$ (over all subsets $\mathcal{M}$ along the model path) or specific convergence rates are needed to ensure validity.

\section{Additional simulation results} \label{Sec: additional-simulation}
This section presents additional simulation results for the robust DA argmin tests (\Cref{Sec: additional-robust}), for the DA‑MCS procedures (\Cref{Sec: additional-simulation-mcs}), for the smallest-mean confidence sets (\Cref{Sec: simulation-smallest-mean}), and for the argmax inference (\Cref{Sec: real-data-argmax}).

\subsection{Robust DA argmin tests} \label{Sec: additional-robust}
This subsection provides additional simulation results for the robust DA argmin tests introduced in \Cref{Sec: robust}. The simulation settings are similar to those in \Cref{Sec: power} and \Cref{Sec: power-unequal}, except that we increase the number of observations to $2n=3000$ and we generate the data from a heavy-tailed distribution—specifically a multivariate $t$-distribution with $3$ degrees of freedom, which has a finite second moment but an infinite third moment. For each sample, a standard normal vector $\bZ \sim N(\boldsymbol{0}, \bm{I}_d)$ is drawn and a chi-squared random variable $U \sim \chi^2_3$ is generated independently. The observed data is then generated as 
\begin{align*}
	\bX = \bm{\mu} + \frac{1}{\sqrt{U/3}} \bm{L} \bZ,
\end{align*}
where $\bm{L}$ is the lower triangular matrix from the Cholesky decomposition of $\bm{\Sigma}$ such that $\bm{\Sigma} = \bm{L} \bm{L}^\top$. The location parameter $\bm{\mu}$ and the covariance matrix $\bm{\Sigma}$ are the same as those used in \Cref{Sec: power} and \Cref{Sec: power-unequal}.

In addition to the MoM estimator, we also consider an alternative robust estimator, namely Catoni's M-estimator~\citep{catoni2012challenging}, described below. Like the MoM estimator, Catoni's estimator achieves sub-Gaussian concentration under the assumption of only finite variance. 

\paragraph{Catoni's M-estimator.} Suppose we observe $X_1,\ldots,X_n$ with finite variance $\sigma^2$. Catoni's estimator $\widehat{\theta}_{\tilde{\alpha}}$ \citep{catoni2012challenging} is defined as the solution of the equation
\begin{align*}
	\sum_{i=1}^n f\bigl(\tilde{\alpha}(X_i - \widehat{\theta}_{\tilde{\alpha}})\bigr) = 0,
\end{align*}
where the function $f$ is given by
\begin{align*}
	f(u) =
\begin{cases}
\log \bigl(1 + u + u^2/2\bigr), & \text{if } u \geq 0, \\
-\log \bigl(1 - u + u^2/2\bigr), & \text{if } u < 0,
\end{cases}	
\end{align*}
and the tuning parameter $\tilde{\alpha}$ is defined as
\begin{align*}
	\tilde{\alpha} = \sqrt{\frac{2\log(1/\delta)}{n(\sigma^2 + \eta^2)}} \quad \text{and} \quad  \eta = \sqrt{\frac{2\sigma^2 \log(1/\delta)}{n -2\log(1/\delta)}}.
\end{align*}
In our implementation, we replace the unknown variance $\sigma^2$ with the sample variance $\widehat{\sigma}^2$ and set the confidence level parameter $\delta = 0.05$.

There are four robust methods that we consider in this section, which are described as follows:
\begin{itemize}
	\item \texttt{DA-plug-mom}: The robust DA argmin test using the MoM plug-in selection method $\widetilde{s}_{\mathrm{plug}}$ with the number of partitions $V = \floor{\sqrt{n}}$.
	\item \texttt{DA-plug-catoni}: This variant is defined as \texttt{DA-plug-mom} but uses Catoni's M-estimator instead of the MoM estimator.
	\item \texttt{DA-adj-mom}: The robust DA argmin test using the MoM noise-adjusted selection method $\widetilde{s}_{\mathrm{adj}}$ with the number of partitions $V = \floor{\sqrt{n}}$.
	\item \texttt{DA-adj-catoni}: This variant is defined as \texttt{DA-adj-mom} but uses Catoni's M-estimator instead of the MoM estimator.
\end{itemize}

The results are summarized in \Cref{tab:rejection-rates-robust,tab:type1-robust} for the homoskedastic setting and in \Cref{tab:rejection-rates-hetero,tab:type1-hetero} for the heteroskedastic setting, based on $10,\!000$ repetitions. Overall, the robust DA argmin tests using Catoni’s estimator perform comparably to their non-robust counterparts. In contrast, the MoM-based versions tend to exhibit slightly lower power, likely due to the inefficiency introduced by sample splitting. We also find that varying the choices of $V$ and $\tilde{\alpha}$ has little impact on performance, while more extreme settings lead to worse results.

\begin{table}[h!]
	\centering
	\begin{minipage}{\textwidth}
		\small	
		\centering
		\caption{Empirical power at the significance level $\alpha = 0.05$ for different mean structures and correlation levels under equal variance and robust settings. The highest power in each scenario is highlighted in bold, and deeper color intensity indicates higher power.} 
		\label{tab:rejection-rates-robust}
		\renewcommand{\arraystretch}{0.95}
		\setlength{\tabcolsep}{5pt}
		\begin{tabular}{lccc ccc ccc}
		\toprule
		\multirow{2}{*}{\textbf{Method}} 
		& \multicolumn{3}{c}{$\bm{\mu}^{(a)}$ + equal variance} 
		& \multicolumn{3}{c}{$\bm{\mu}^{(b)}$ + equal variance} 
		& \multicolumn{3}{c}{$\bm{\mu}^{(c)}$ + equal variance} \\
		& $\rho=0$ & $\rho=0.4$ & $\rho=0.8$ 
		& $\rho=0$ & $\rho=0.4$ & $\rho=0.8$ 
		& $\rho=0$ & $\rho=0.4$ & $\rho=0.8$ \\
		\cmidrule(lr){2-4} \cmidrule(lr){5-7} \cmidrule(lr){8-10}
		\texttt{DA-plug} & \cellcolor{lightorange}0.235 & \cellcolor{lightorange}0.320 & \cellcolor{midorange}0.536 & \cellcolor{lightorange}0.379 & \cellcolor{midorange}0.444 & \cellcolor{strongorange}0.686 & \cellcolor{vlightorange}0.198 & \cellcolor{lightorange}0.241 & \cellcolor{midorange}0.424 \\
		\texttt{DA-plug-mom} & \cellcolor{vlightorange}0.181 & \cellcolor{lightorange}0.238 & \cellcolor{lightorange}0.366 & \cellcolor{lightorange}0.355 & \cellcolor{midorange}0.414 & \cellcolor{strongorange}0.620 & \cellcolor{lightorange}0.205 & \cellcolor{lightorange}0.240 & \cellcolor{midorange}0.427 \\
		\texttt{DA-plug-catoni} & \cellcolor{lightorange}0.238 & \cellcolor{lightorange}0.334 & \cellcolor{midorange}0.539 & \cellcolor{lightorange}0.376 & \cellcolor{midorange}0.440 & \cellcolor{strongorange}0.699 & \cellcolor{lightorange}0.202 & \cellcolor{lightorange}0.233 & \cellcolor{midorange}0.427 \\
		\texttt{DA-adj} & \cellcolor{lightorange}0.232 & \cellcolor{midorange}0.450 & \cellcolor{deeporange}0.935 & \cellcolor{lightorange}0.378 & \cellcolor{midorange}0.509 & \cellcolor{deeporange}\textbf{0.928} & \cellcolor{lightorange}\textbf{0.206} & \cellcolor{lightorange}0.243 & \cellcolor{midorange}\textbf{0.484} \\
		\texttt{DA-adj-mom} & \cellcolor{vlightorange}0.179 & \cellcolor{lightorange}0.350 & \cellcolor{deeporange}0.829 & \cellcolor{lightorange}0.360 & \cellcolor{midorange}0.480 & \cellcolor{deeporange}0.860 & \cellcolor{lightorange} 0.202 & \cellcolor{lightorange}\textbf{0.246} & \cellcolor{midorange}0.466 \\
		\texttt{DA-adj-catoni} & \cellcolor{lightorange}\textbf{0.243} & \cellcolor{midorange}\textbf{0.455} & \cellcolor{deeporange}\textbf{0.938} & \cellcolor{lightorange}\textbf{0.381} & \cellcolor{midorange}\textbf{0.514} & \cellcolor{deeporange}0.927 & \cellcolor{lightorange}\textbf{0.206} & \cellcolor{lightorange}0.242 & \cellcolor{midorange}0.473 \\
		\bottomrule
		\end{tabular}
	\end{minipage}

	\vspace{0.5em}

	\begin{minipage}{\textwidth}
		\small
		\centering
		\caption{Empirical type I error at the significance level $\alpha = 0.05$ for different mean structures and correlation levels under equal variance and robust settings. Green indicates under-rejection (conservative tests), and white indicates appropriate rejection rates (correct coverage).}
		\renewcommand{\arraystretch}{0.95}
		\setlength{\tabcolsep}{5pt}
		\begin{tabular}{lccc ccc ccc}
		\toprule
		\multirow{2}{*}{\textbf{Method}} 
		& \multicolumn{3}{c}{$\bm{\mu}^{(a,0)}$ + equal variance} 
		& \multicolumn{3}{c}{$\bm{\mu}^{(b,0)}$ + equal variance} 
		& \multicolumn{3}{c}{$\bm{\mu}^{(c,0)}$ + equal variance} \\
		& $\rho=0$ & $\rho=0.4$ & $\rho=0.8$ 
		& $\rho=0$ & $\rho=0.4$ & $\rho=0.8$ 
		& $\rho=0$ & $\rho=0.4$ & $\rho=0.8$ \\
		\cmidrule(lr){2-4} \cmidrule(lr){5-7} \cmidrule(lr){8-10}
		\texttt{DA-plug} & \cellcolor{lightgreen}0.021 & \cellcolor{lightgreen}0.023 & \cellcolor{lightgreen}0.028 & \cellcolor{lightgreen}0.030 & \cellcolor{lightgreen}0.029 & \cellcolor{lightgreen}0.026 & 0.052 & 0.049 & 0.053 \\
		\texttt{DA-plug-mom} & \cellcolor{midgreen}0.015 & \cellcolor{midgreen}0.014 & \cellcolor{midgreen}0.019 & \cellcolor{lightgreen}0.028 & \cellcolor{lightgreen}0.029 & \cellcolor{lightgreen}0.025 & 0.054 & 0.052 & 0.053 \\
		\texttt{DA-plug-catoni} & \cellcolor{lightgreen}0.022 & \cellcolor{lightgreen}0.022 & \cellcolor{lightgreen}0.024 & \cellcolor{lightgreen}0.028 & \cellcolor{lightgreen}0.027 & \cellcolor{lightgreen}0.030 & 0.053 & 0.045 & 0.051 \\
		\texttt{DA-adj} & \cellcolor{lightgreen}0.022 & \cellcolor{lightgreen}0.022 & \cellcolor{lightgreen}0.024 & \cellcolor{lightgreen}0.027 & \cellcolor{lightgreen}0.028 & \cellcolor{lightgreen}0.029 & 0.047 & 0.050 & 0.048 \\
		\texttt{DA-adj-mom} & \cellcolor{midgreen}0.016 & \cellcolor{midgreen}0.016 & \cellcolor{lightgreen}0.021 & \cellcolor{lightgreen}0.028 & \cellcolor{lightgreen}0.032 & \cellcolor{lightgreen}0.027 & 0.053 & 0.053 & 0.049 \\
		\texttt{DA-adj-catoni} & \cellcolor{lightgreen}0.023 & \cellcolor{lightgreen}0.022 & \cellcolor{lightgreen}0.028 & \cellcolor{lightgreen}0.028 & \cellcolor{lightgreen}0.029 & \cellcolor{lightgreen}0.028 & 0.050 & 0.049 & 0.050 \\
		\bottomrule
		\end{tabular}
		\label{tab:type1-robust}
	\end{minipage}
	
	\vspace{0.5em}

	\begin{minipage}{\textwidth}
		\small	
		\centering
		\caption{Empirical power at the significance level $\alpha = 0.05$ for different mean structures and correlation levels under unequal variance and robust settings. The highest power in each scenario is highlighted in bold, and deeper color intensity indicates higher power.} \label{tab:rejection-rates-hetero}
		\renewcommand{\arraystretch}{0.95}
		\setlength{\tabcolsep}{5pt}
		\begin{tabular}{lccccccccc}
		\toprule
		\multirow{2}{*}{\textbf{Method}} & \multicolumn{3}{c}{$\bm{\mu}^{(a)}$ + unequal variance} & \multicolumn{3}{c}{$\bm{\mu}^{(b)}$ + unequal variance} & \multicolumn{3}{c}{$\bm{\mu}^{(c)}$ + unequal variance} \\
		& $\rho=0$ & $\rho=0.4$ & $\rho=0.8$ & $\rho=0$ & $\rho=0.4$ & $\rho=0.8$ & $\rho=0$ & $\rho=0.4$ & $\rho=0.8$ \\
		\cmidrule(lr){2-4} \cmidrule(lr){5-7} \cmidrule(lr){8-10}
		\texttt{DA-plug} & \cellcolor{vlightorange}0.049 & \cellcolor{vlightorange}0.053 & \cellcolor{vlightorange}0.045 & \cellcolor{vlightorange}0.066 & \cellcolor{vlightorange}0.069 & \cellcolor{vlightorange}0.066 & \cellcolor{vlightorange}0.115 & \cellcolor{vlightorange}0.126 & \cellcolor{lightorange}0.206 \\
		\texttt{DA-plug-mom} & \cellcolor{vlightorange}0.048 & \cellcolor{vlightorange}0.049 & \cellcolor{vlightorange}0.052 & \cellcolor{vlightorange}0.063 & \cellcolor{vlightorange}0.062 & \cellcolor{vlightorange}0.059 & \cellcolor{vlightorange}0.113 & \cellcolor{vlightorange}0.133 & \cellcolor{lightorange}0.206 \\
		\texttt{DA-plug-catoni} & \cellcolor{vlightorange}0.052 & \cellcolor{vlightorange}0.050 & \cellcolor{vlightorange}0.051 & \cellcolor{vlightorange}0.066 & \cellcolor{vlightorange}0.069 & \cellcolor{vlightorange}0.066 & \cellcolor{vlightorange}0.106 & \cellcolor{vlightorange}0.129 & \cellcolor{lightorange}0.203 \\
		\texttt{DA-adj} & \cellcolor{vlightorange}0.140 & \cellcolor{lightorange}0.281 & \cellcolor{strongorange}0.838 & \cellcolor{lightorange}0.224 & \cellcolor{midorange}0.406 & \cellcolor{deeporange}0.904 & \cellcolor{vlightorange}\textbf{0.146} & \cellcolor{vlightorange}\textbf{0.196} & \cellcolor{midorange}\textbf{0.468} \\
		\texttt{DA-adj-mom} & \cellcolor{vlightorange}0.106 & \cellcolor{lightorange}0.223 & \cellcolor{strongorange}0.656 & \cellcolor{vlightorange}0.168 & \cellcolor{lightorange}0.318 & \cellcolor{strongorange}0.777 & \cellcolor{vlightorange}0.146 & \cellcolor{vlightorange}0.180 & \cellcolor{midorange}0.412 \\
		\texttt{DA-adj-catoni} & \cellcolor{vlightorange}\textbf{0.147} & \cellcolor{lightorange}\textbf{0.284} & \cellcolor{strongorange}\textbf{0.852} & \cellcolor{lightorange}\textbf{0.230} & \cellcolor{midorange}\textbf{0.410} & \cellcolor{deeporange}\textbf{0.906} & \cellcolor{vlightorange}\textbf{0.146} & \cellcolor{vlightorange}0.186 & \cellcolor{midorange}\textbf{0.468} \\
		\bottomrule
		\end{tabular}
	\end{minipage}

	\vspace{0.5em}

	\begin{minipage}{\textwidth}
		\small
		\centering
		\caption{Empirical type I error at the significance level $\alpha = 0.05$ for different mean structures and correlation levels under unequal variance and robust settings. Green indicates under-rejection (conservative tests), and white indicates appropriate rejection rates (correct coverage).}
		\renewcommand{\arraystretch}{0.95}
		\setlength{\tabcolsep}{5pt}
		\begin{tabular}{lccc ccc ccc}
		\toprule
		\multirow{2}{*}{\textbf{Method}} 
		& \multicolumn{3}{c}{$\bm{\mu}^{(a,0)}$ + equal variance} 
		& \multicolumn{3}{c}{$\bm{\mu}^{(b,0)}$ + equal variance} 
		& \multicolumn{3}{c}{$\bm{\mu}^{(c,0)}$ + equal variance} \\
		& $\rho=0$ & $\rho=0.4$ & $\rho=0.8$ 
		& $\rho=0$ & $\rho=0.4$ & $\rho=0.8$ 
		& $\rho=0$ & $\rho=0.4$ & $\rho=0.8$ \\
		\cmidrule(lr){2-4} \cmidrule(lr){5-7} \cmidrule(lr){8-10}
		\texttt{DA-plug} & \cellcolor{midgreen}0.016 & \cellcolor{midgreen}0.017 & \cellcolor{midgreen}0.014 & \cellcolor{lightgreen}0.021 & \cellcolor{lightgreen}0.023 & \cellcolor{lightgreen}0.023 & 0.048 & 0.049 & 0.051 \\
		\texttt{DA-plug-mom} & \cellcolor{midgreen}0.016 & \cellcolor{midgreen}0.016 & \cellcolor{midgreen}0.016 & \cellcolor{lightgreen}0.020 & \cellcolor{midgreen}0.019 & \cellcolor{lightgreen}0.022 & 0.050 & 0.057 & 0.048 \\
		\texttt{DA-plug-catoni} & \cellcolor{midgreen}0.018 & \cellcolor{midgreen}0.016 & \cellcolor{midgreen}0.013 & \cellcolor{midgreen}0.019 & \cellcolor{lightgreen}0.024 & \cellcolor{lightgreen}0.022 & 0.049 & 0.050 & 0.050 \\
		\texttt{DA-adj} & \cellcolor{midgreen}0.016 & \cellcolor{midgreen}0.018 & \cellcolor{midgreen}0.018 & \cellcolor{lightgreen}0.024 & \cellcolor{lightgreen}0.025 & \cellcolor{lightgreen}0.025 & 0.046 & 0.052 & 0.054 \\
		\texttt{DA-adj-mom} & \cellcolor{midgreen}0.015 & \cellcolor{midgreen}0.018 & \cellcolor{midgreen}0.019 & \cellcolor{lightgreen}0.024 & \cellcolor{lightgreen}0.024 & \cellcolor{lightgreen}0.024 & 0.050 & 0.051 & 0.048 \\
		\texttt{DA-adj-catoni} & \cellcolor{midgreen}0.015 & \cellcolor{midgreen}0.017 & \cellcolor{midgreen}0.018 & \cellcolor{lightgreen}0.025 & \cellcolor{lightgreen}0.024 & \cellcolor{lightgreen}0.027 & 0.047 & 0.049 & 0.050 \\
		\bottomrule
		\end{tabular}
		\label{tab:type1-hetero}
	\end{minipage}
\end{table}

\clearpage 

\subsection{Simulations for DA-MCS procedures} \label{Sec: additional-simulation-mcs}

In this subsection, we present additional simulation results for the DA-MCS procedures introduced in \Cref{Sec: MCS}. The simulation settings are identical to those in \Cref{Sec: Sim MCS} except that we change the mean gap parameter to $\zeta = 3\sqrt{\log |\Theta| / (2n)}$ and $\zeta = 1$. 

The simulation results for $\zeta = 3\sqrt{\log |\Theta| / (2n)}$ and $\zeta = 1$ are presented in \Cref{tab:uniform-coverage-3} and \Cref{tab:uniform-coverage-1}, respectively. These findings align closely with the observations in \Cref{Sec: Sim MCS}: the DA-MCS procedures consistently achieve nominal coverage across all scenarios, whereas the pointwise methods (\texttt{DA-plug} and \texttt{DA-adj}) exhibit notable under-coverage. Specifically, in the more challenging scenario where $\zeta = 3\sqrt{\log |\Theta| / (2n)}$ (which is smaller than 1 with $n = 500$), both the one-step and two-step uniform procedures show a tendency to over-cover the true parameter set. In contrast, when $\zeta = 1$, the two-step procedures not only produce coverage rates closer to the nominal level but also result in shorter average confidence set lengths compared to their one-step counterparts, highlighting the superior efficiency of the two-step procedures.

\begin{table}[h]
    \centering
    \small
    \caption{Empirical coverage probabilities $P(\Theta \subseteq \widehat{\Theta})$ across varying cardinalities $|\Theta|$ and correlation levels $\rho$, evaluated for six different methods at the nominal level $1-\alpha = 0.95$. The mean gap $\zeta$ is set to $3\sqrt{\log|\Theta|/(2n)}$. Numbers in parentheses indicate the average length of the confidence sets. Under-coverage rates are shaded in progressively darker blue, over-coverage rates in progressively darker green, and rates close to the nominal level remain unshaded.}
    \label{tab:uniform-coverage-3}
    \renewcommand{\arraystretch}{1.0}
    \setlength{\tabcolsep}{2.7pt}

	\begin{tabular}{cc}

			\begin{minipage}{0.45\textwidth}
			\centering
			\texttt{DA-plug} (pointwise)\\[0.3em]
			\begin{tabular}{llll}\toprule
			  & $\rho=0$ & $\rho=0.4$ & $\rho=0.8$\\\midrule
			  $|\Theta|=2$  & \cellcolor{lightgreen}0.971 {\tiny(82.31)} & \cellcolor{lightgreen}0.971 {\tiny(83.11)} & 0.959 {\tiny(81.72)}\\
			  $|\Theta|=5$  & \cellcolor{lightblue}0.822 {\tiny(45.50)} & \cellcolor{lightblue}0.829 {\tiny(46.42)} & \cellcolor{lightblue}0.841 {\tiny(48.33)}\\
			  $|\Theta|=10$ & \cellcolor{midblue}0.718 {\tiny(33.31)} & \cellcolor{midblue}0.720 {\tiny(33.57)} & \cellcolor{midblue}0.748 {\tiny(33.36)}\\
			  $|\Theta|=15$ & \cellcolor{strongblue}0.662 {\tiny(31.69)} & \cellcolor{strongblue}0.669 {\tiny(31.63)} & \cellcolor{strongblue}0.688 {\tiny(31.53)}\\
			  $|\Theta|=20$ & \cellcolor{strongblue}0.623 {\tiny(32.75)} & \cellcolor{strongblue}0.616 {\tiny(32.84)} & \cellcolor{strongblue}0.651 {\tiny(32.56)}\\\bottomrule
			\end{tabular}
			\end{minipage}
			&
			\begin{minipage}{0.45\textwidth}
			\centering
			\texttt{DA-adj} (pointwise)\\[0.3em]
			\begin{tabular}{llll}\toprule
			  & $\rho=0$ & $\rho=0.4$ & $\rho=0.8$\\\midrule
			  $|\Theta|=2$  & \cellcolor{lightgreen}0.973 {\tiny(82.62)} & \cellcolor{lightgreen}0.966 {\tiny(83.35)} & \cellcolor{lightgreen}0.960 {\tiny(82.09)}\\
			  $|\Theta|=5$  & \cellcolor{lightblue}0.817 {\tiny(45.24)} & \cellcolor{lightblue}0.823 {\tiny(46.62)} & \cellcolor{lightblue}0.830 {\tiny(50.56)}\\
			  $|\Theta|=10$ & \cellcolor{midblue}0.720 {\tiny(33.30)} & \cellcolor{midblue}0.705 {\tiny(33.57)} & \cellcolor{midblue}0.711 {\tiny(34.20)}\\
			  $|\Theta|=15$ & \cellcolor{strongblue}0.644 {\tiny(31.64)} & \cellcolor{strongblue}0.641 {\tiny(31.95)} & \cellcolor{strongblue}0.632 {\tiny(31.66)}\\
			  $|\Theta|=20$ & \cellcolor{strongblue}0.607 {\tiny(32.48)} & \cellcolor{strongblue}0.587 {\tiny(32.84)} & \cellcolor{strongblue}0.582 {\tiny(32.05)}\\\bottomrule
			\end{tabular}
			\end{minipage}\\[5em]
		
			\begin{minipage}{0.45\textwidth}
			\centering
			\texttt{DA-MCS-plug}$^{1}$ (uniform)\\[0.3em]
			\begin{tabular}{llll}\toprule
			  & $\rho=0$ & $\rho=0.4$ & $\rho=0.8$\\\midrule
			  $|\Theta|=2$  & \cellcolor{midgreen}1.000 {\tiny(99.05)} & \cellcolor{midgreen}1.000 {\tiny(99.05)} & \cellcolor{midgreen}1.000 {\tiny(98.92)}\\
			  $|\Theta|=5$  & \cellcolor{midgreen}0.997 {\tiny(92.70)} & \cellcolor{midgreen}0.998 {\tiny(92.63)} & \cellcolor{midgreen}0.998 {\tiny(92.15)}\\
			  $|\Theta|=10$ & \cellcolor{midgreen}0.997 {\tiny(85.85)} & \cellcolor{midgreen}0.996 {\tiny(85.97)} & \cellcolor{midgreen}0.996 {\tiny(85.11)}\\
			  $|\Theta|=15$ & \cellcolor{midgreen}0.993 {\tiny(82.26)} & \cellcolor{midgreen}0.993 {\tiny(82.23)} & \cellcolor{midgreen}0.995 {\tiny(81.96)}\\
			  $|\Theta|=20$ & \cellcolor{midgreen}0.991 {\tiny(80.30)} & \cellcolor{midgreen}0.992 {\tiny(80.49)} & \cellcolor{midgreen}0.993 {\tiny(79.68)}\\\bottomrule
			\end{tabular}
			\end{minipage}
			&
			\begin{minipage}{0.45\textwidth}
			\centering
			\texttt{DA-MCS-adj}$^{1}$ (uniform)\\[0.3em]
			\begin{tabular}{llll}\toprule
			  & $\rho=0$ & $\rho=0.4$ & $\rho=0.8$\\\midrule
			  $|\Theta|=2$  & \cellcolor{midgreen}1.000 {\tiny(99.08)} & \cellcolor{midgreen}1.000 {\tiny(99.11)} & \cellcolor{midgreen}1.000 {\tiny(98.77)}\\
			  $|\Theta|=5$  & \cellcolor{midgreen}0.998 {\tiny(92.49)} & \cellcolor{midgreen}0.998 {\tiny(92.58)} & \cellcolor{midgreen}0.997 {\tiny(91.94)}\\
			  $|\Theta|=10$ & \cellcolor{midgreen}0.996 {\tiny(85.67)} & \cellcolor{midgreen}0.996 {\tiny(86.04)} & \cellcolor{midgreen}0.995 {\tiny(84.60)}\\
			  $|\Theta|=15$ & \cellcolor{midgreen}0.992 {\tiny(82.45)} & \cellcolor{midgreen}0.993 {\tiny(82.63)} & \cellcolor{midgreen}0.993 {\tiny(80.96)}\\
			  $|\Theta|=20$ & \cellcolor{midgreen}0.991 {\tiny(80.33)} & \cellcolor{midgreen}0.991 {\tiny(80.50)} & \cellcolor{midgreen}0.993 {\tiny(78.41)}\\\bottomrule
			\end{tabular}
			\end{minipage}\\[5em]
		
			\begin{minipage}{0.45\textwidth}
			\centering
			\texttt{DA-MCS-plug}$^{2}$ (uniform)\\[0.3em]
			\begin{tabular}{llll}\toprule
			  & $\rho=0$ & $\rho=0.4$ & $\rho=0.8$\\\midrule
			  $|\Theta|=2$  & \cellcolor{midgreen}1.000 {\tiny(98.93)} & \cellcolor{midgreen}1.000 {\tiny(98.94)} & \cellcolor{midgreen}1.000 {\tiny(98.77)}\\
			  $|\Theta|=5$  & \cellcolor{midgreen}0.996 {\tiny(90.97)} & \cellcolor{midgreen}0.997 {\tiny(91.08)} & \cellcolor{midgreen}0.996 {\tiny(90.52)}\\
			  $|\Theta|=10$ & \cellcolor{midgreen}0.992 {\tiny(81.67)} & \cellcolor{midgreen}0.993 {\tiny(82.20)} & \cellcolor{midgreen}0.994 {\tiny(81.34)}\\
			  $|\Theta|=15$ & \cellcolor{midgreen}0.987 {\tiny(78.14)} & \cellcolor{midgreen}0.988 {\tiny(77.17)} & \cellcolor{midgreen}0.991 {\tiny(76.92)}\\
			  $|\Theta|=20$ & \cellcolor{midgreen}0.985 {\tiny(75.52)} & \cellcolor{midgreen}0.982 {\tiny(75.33)} & \cellcolor{midgreen}0.987 {\tiny(74.52)}\\\bottomrule
			\end{tabular}
			\end{minipage}
			&
			\begin{minipage}{0.45\textwidth}
			\centering
			\texttt{DA-MCS-adj}$^{2}$ (uniform)\\[0.3em]
			\begin{tabular}{llll}\toprule
			  & $\rho=0$ & $\rho=0.4$ & $\rho=0.8$\\\midrule
			  $|\Theta|=2$  & \cellcolor{midgreen}1.000 {\tiny(98.92)} & \cellcolor{midgreen}1.000 {\tiny(99.01)} & \cellcolor{midgreen}0.999 {\tiny(98.68)}\\
			  $|\Theta|=5$  & \cellcolor{midgreen}0.996 {\tiny(90.75)} & \cellcolor{midgreen}0.996 {\tiny(91.11)} & \cellcolor{midgreen}0.997 {\tiny(90.38)}\\
			  $|\Theta|=10$ & \cellcolor{midgreen}0.993 {\tiny(81.74)} & \cellcolor{midgreen}0.990 {\tiny(82.21)} & \cellcolor{midgreen}0.994 {\tiny(80.99)}\\
			  $|\Theta|=15$ & \cellcolor{midgreen}0.987 {\tiny(78.03)} & \cellcolor{midgreen}0.988 {\tiny(77.70)} & \cellcolor{midgreen}0.986 {\tiny(76.40)}\\
			  $|\Theta|=20$ & \cellcolor{midgreen}0.984 {\tiny(75.39)} & \cellcolor{midgreen}0.983 {\tiny(75.25)} & \cellcolor{midgreen}0.985 {\tiny(73.57)}\\\bottomrule
			\end{tabular}
			\end{minipage}
		\end{tabular}
\end{table}

\begin{table}[ht]
    \centering
    \small
    \caption{Empirical coverage probabilities $P(\Theta \subseteq \widehat{\Theta})$ across varying cardinalities $|\Theta|$ and correlation levels $\rho$, evaluated for six different methods at the nominal level $1-\alpha = 0.95$. The mean gap $\zeta$ is set to $1$. Numbers in parentheses indicate the average length of the confidence sets. Under-coverage rates are shaded in progressively darker blue, over-coverage rates in progressively darker green, and rates close to the nominal level remain unshaded.}
    \label{tab:uniform-coverage-1}
    \renewcommand{\arraystretch}{1.0}
	\setlength{\tabcolsep}{2.7pt}

	\begin{tabular}{cc}

		\begin{minipage}{0.45\textwidth}
		\centering
		\texttt{DA-plug} (pointwise)\\[0.3em]
		\begin{tabular}{llll}\toprule
		  & $\rho=0$                      & $\rho=0.4$                    & $\rho=0.8$                    \\\midrule
		  $|\Theta|=2$  & \cellcolor{vlightblue}0.905 {\tiny(1.90)} & \cellcolor{lightblue}0.894 {\tiny(1.90)} & \cellcolor{vlightblue} 0.900{\tiny(1.90)} \\
		  $|\Theta|=5$  & \cellcolor{lightblue}0.800 {\tiny(4.74)} & \cellcolor{lightblue}0.811 {\tiny(4.76)} & \cellcolor{lightblue}0.817 {\tiny(4.76)} \\
		  $|\Theta|=10$ & \cellcolor{midblue}0.720 {\tiny(9.48)} & \cellcolor{midblue}0.712 {\tiny(9.49)} & \cellcolor{midblue}0.739 {\tiny(9.50)} \\
		  $|\Theta|=15$ & \cellcolor{strongblue}0.658 {\tiny(14.26)} & \cellcolor{strongblue}0.665 {\tiny(14.23)} & \cellcolor{strongblue}0.688 {\tiny(14.24)} \\
		  $|\Theta|=20$ & \cellcolor{strongblue}0.616 {\tiny(19.03)} & \cellcolor{strongblue}0.615 {\tiny(18.98)} & \cellcolor{strongblue}0.648 {\tiny(19.01)} \\\bottomrule
		\end{tabular}
		\end{minipage}
		&
		\begin{minipage}{0.45\textwidth}
		\centering
		\texttt{DA-adj} (pointwise)\\[0.3em]
		\begin{tabular}{llll}\toprule
		  & $\rho=0$                       & $\rho=0.4$                     & $\rho=0.8$                     \\\midrule
		  $|\Theta|=2$  & \cellcolor{vlightblue}0.904 {\tiny(1.90)} & \cellcolor{vlightblue}0.901 {\tiny(1.90)} & \cellcolor{lightblue}0.898 {\tiny(1.90)} \\
		  $|\Theta|=5$  & \cellcolor{lightblue}0.800 {\tiny(4.75)} & \cellcolor{lightblue}0.808 {\tiny(4.76)} & \cellcolor{lightblue}0.811 {\tiny(4.74)} \\
		  $|\Theta|=10$ & \cellcolor{midblue}0.714 {\tiny(9.48)} & \cellcolor{midblue}0.712 {\tiny(9.49)} & \cellcolor{midblue}0.705 {\tiny(9.51)} \\
		  $|\Theta|=15$ & \cellcolor{strongblue}0.652 {\tiny(14.22)} & \cellcolor{strongblue}0.644 {\tiny(14.25)} & \cellcolor{strongblue}0.638 {\tiny(14.25)} \\
		  $|\Theta|=20$ & \cellcolor{strongblue}0.607 {\tiny(19.00)} & \cellcolor{strongblue}0.590 {\tiny(18.97)} & \cellcolor{strongblue}0.579 {\tiny(18.97)} \\\bottomrule
		\end{tabular}
		\end{minipage}\\[5em]
		
		\begin{minipage}{0.45\textwidth}
		\centering
		\texttt{DA-MCS-plug}$^{1}$ (uniform)\\[0.3em]
		\begin{tabular}{llll}\toprule
		  & $\rho=0$                        & $\rho=0.4$                      & $\rho=0.8$                      \\\midrule
		  $|\Theta|=2$  & \cellcolor{midgreen}0.999 {\tiny(2.00)} & \cellcolor{midgreen}0.999 {\tiny(2.00)} & \cellcolor{midgreen}0.999 {\tiny(2.00)} \\
		  $|\Theta|=5$  & \cellcolor{midgreen}0.998 {\tiny(5.00)} & \cellcolor{midgreen}0.998 {\tiny(5.00)} & \cellcolor{midgreen}0.997 {\tiny(5.00)} \\
		  $|\Theta|=10$ & \cellcolor{midgreen}0.994 {\tiny(9.99)} & \cellcolor{midgreen}0.996 {\tiny(10.00)} & \cellcolor{midgreen}0.997 {\tiny(10.00)} \\
		  $|\Theta|=15$ & \cellcolor{midgreen}0.994 {\tiny(14.99)} & \cellcolor{midgreen}0.995 {\tiny(14.99)} & \cellcolor{midgreen}0.996 {\tiny(14.99)} \\
		  $|\Theta|=20$ & \cellcolor{midgreen}0.993 {\tiny(19.99)} & \cellcolor{midgreen}0.991 {\tiny(19.99)} & \cellcolor{midgreen}0.993 {\tiny(19.99)} \\\bottomrule
		\end{tabular}
		\end{minipage}
		&
		\begin{minipage}{0.45\textwidth}
		\centering
		\texttt{DA-MCS-adj}$^{1}$ (uniform)\\[0.3em]
		\begin{tabular}{llll}\toprule
		  & $\rho=0$                        & $\rho=0.4$                      & $\rho=0.8$                      \\\midrule
		  $|\Theta|=2$  & \cellcolor{midgreen}0.999 {\tiny(2.00)} & \cellcolor{midgreen}0.999 {\tiny(2.00)} & \cellcolor{midgreen}0.998 {\tiny(2.00)} \\
		  $|\Theta|=5$  & \cellcolor{midgreen}0.997 {\tiny(5.00)} & \cellcolor{midgreen}0.997 {\tiny(5.00)} & \cellcolor{midgreen}0.997 {\tiny(5.00)} \\
		  $|\Theta|=10$ & \cellcolor{midgreen}0.994 {\tiny(9.99)} & \cellcolor{midgreen}0.996 {\tiny(9.99)} & \cellcolor{midgreen}0.996 {\tiny(9.99)} \\
		  $|\Theta|=15$ & \cellcolor{midgreen}0.993 {\tiny(14.99)} & \cellcolor{midgreen}0.994 {\tiny(14.99)} & \cellcolor{midgreen}0.994 {\tiny(14.99)} \\
		  $|\Theta|=20$ & \cellcolor{midgreen}0.992 {\tiny(19.99)} & \cellcolor{midgreen}0.991 {\tiny(19.99)} & \cellcolor{midgreen}0.992 {\tiny(19.99)} \\\bottomrule
		\end{tabular}
		\end{minipage}\\[5em]
		
		\begin{minipage}{0.45\textwidth}
		\centering
		\texttt{DA-MCS-plug}$^{2}$ (uniform)\\[0.3em]
		\begin{tabular}{llll}\toprule
		  & $\rho=0$                       & $\rho=0.4$                      & $\rho=0.8$                      \\\midrule
		  $|\Theta|=2$  & 0.947 {\tiny(1.94)}               & 0.946 {\tiny(1.95)}               & 0.942 {\tiny(1.94)}               \\
		  $|\Theta|=5$  & 0.954 {\tiny(4.95)}               & 0.951 {\tiny(4.95)}               & 0.956 {\tiny(4.95)}               \\
		  $|\Theta|=10$ & \cellcolor{lightgreen}0.960 {\tiny(9.94)} & \cellcolor{lightgreen}0.960 {\tiny(9.94)} & \cellcolor{lightgreen}0.964 {\tiny(9.95)} \\
		  $|\Theta|=15$ & \cellcolor{lightgreen}0.963 {\tiny(14.94)} & \cellcolor{lightgreen}0.961 {\tiny(14.94)} & \cellcolor{lightgreen}0.969 {\tiny(14.94)} \\
		  $|\Theta|=20$ & \cellcolor{lightgreen}0.963 {\tiny(19.94)} & \cellcolor{lightgreen}0.962 {\tiny(19.94)} & \cellcolor{lightgreen}0.971 {\tiny(19.95)} \\\bottomrule
		\end{tabular}
		\end{minipage}
		&
		\begin{minipage}{0.45\textwidth}
		\centering
		\texttt{DA-MCS-adj}$^{2}$ (uniform)\\[0.3em]
		\begin{tabular}{llll}\toprule
		  & $\rho=0$                       & $\rho=0.4$                      & $\rho=0.8$                      \\\midrule
		  $|\Theta|=2$  & 0.951 {\tiny(1.95)}               & 0.943 {\tiny(1.95)}               & 0.941 {\tiny(1.95)}               \\
		  $|\Theta|=5$  & 0.954 {\tiny(4.95)}               & 0.952 {\tiny(4.95)}               & 0.955 {\tiny(4.94)}               \\
		  $|\Theta|=10$ & 0.958 {\tiny(9.94)}               & 0.959 {\tiny(9.94)}               & 0.957 {\tiny(9.95)}               \\
		  $|\Theta|=15$ & \cellcolor{lightgreen}0.961 {\tiny(14.95)} & 0.958 {\tiny(14.94)}               & \cellcolor{lightgreen}0.963 {\tiny(14.94)} \\
		  $|\Theta|=20$ & \cellcolor{lightgreen}0.967 {\tiny(19.95)} & \cellcolor{lightgreen}0.963 {\tiny(19.94)} & \cellcolor{lightgreen}0.960 {\tiny(19.94)} \\\bottomrule
		\end{tabular}
		\end{minipage}
		
		\end{tabular}
\end{table}

\subsection{Simulations of smallest‐mean confidence sets} \label{Sec: simulation-smallest-mean}
We assess the performance of the smallest‐mean confidence sets introduced in Section~\ref{Sec: MCS} via simulation. The simulation settings are the same as in \Cref{Sec: Sim MCS}, with the exception that we set the mean gap $\zeta = 1$ and increase the dimension to $d=1000$. For the data-adaptive set $\mathcal{C}_2$, we set $\gamma_n = \alpha / \log(n)$ and employ the \texttt{DA-MCS-adj}$^2$ procedure to construct the screening set $\widehat{\Theta}_{\mathrm{DA}}^{\mathrm{uni}}$. \Cref{tab:min-mean} reports the average widths of $\mathcal{C}_1$ and $\mathcal{C}_2$ across varying $|\Theta|$ and correlation levels $\rho$. When $|\Theta|$ is small, the data‐adaptive interval $\mathcal{C}_2$ is narrower than the non‐adaptive $\mathcal{C}_1$. However, as $|\Theta|$ grows, this advantage diminishes and $\mathcal{C}_2$ becomes wider due to the efficiency loss from sample splitting. Specifically, the critical value $z_{1-\alpha/(2d)}/\sqrt{2n}$ for $\mathcal{C}_1$ exceeds $z_{1-\alpha/(2\widehat{d})}/\sqrt{n}$ for $\mathcal{C}_2$ only when $\widehat{d} \ll d$. Once $\widehat{d}$ approaches $d$, the sample‐splitting penalty dominates, which is confirmed by the simulation results. Although improving the efficiency of the data‐adaptive procedure represents an intriguing open problem, it falls beyond the scope of this paper and we leave it for future research.

\begin{table}[!h]
    \centering
    \small
    \caption{Comparison of the average widths of the smallest-mean confidence sets $\mathcal{C}_1$ and $\mathcal{C}_2$ (defined in \Cref{Sec: MCS}) shown for varying mean configurations and correlation levels $\rho$. The data-adaptive procedure $\mathcal{C}_2$ yields narrower intervals than the non-adaptive set $\mathcal{C}_1$ when $|\Theta|$ is small, whereas it becomes wider as $|\Theta|$ increases. Numbers in parentheses indicate the average value of $\widehat{d}=|\widehat{\Theta}_{\mathrm{DA}}^{\mathrm{uni}}|$.}
    \label{tab:min-mean}
    \renewcommand{\arraystretch}{1.0}
    \setlength{\tabcolsep}{3pt}
    \begin{tabular}{cc}
		\begin{minipage}{0.45\textwidth}
            \centering
            Average widths of $\mathcal{C}_1$ \\[0.3em]
            \begin{tabular}{lccc}\toprule
              & $\rho=0$    & $\rho=0.4$  & $\rho=0.8$  \\\midrule
              $|\Theta|=2$  & 0.181      & 0.181      & 0.181      \\
              $|\Theta|=5$  & 0.181      & 0.181      & 0.181      \\
              $|\Theta|=10$ & 0.181      & 0.181      & 0.181      \\
              $|\Theta|=15$ & 0.181      & 0.181      & 0.181      \\
              $|\Theta|=20$ & 0.181      & 0.181      & 0.181      \\\bottomrule
            \end{tabular}
        \end{minipage}
		&
        \begin{minipage}{0.5\textwidth}
            \centering
            Average widths of $\mathcal{C}_2$ \\[0.3em]
            \begin{tabular}{lccc}\toprule
              & $\rho=0$    & $\rho=0.4$  & $\rho=0.8$  \\\midrule
              $|\Theta|=2$  & 0.142 {\tiny (1.99)}     & 0.142 {\tiny (1.99)}     & 0.142 {\tiny (1.99)}     \\
              $|\Theta|=5$  & 0.163 {\tiny (4.99)}     & 0.163 {\tiny (4.99)}     & 0.163 {\tiny (4.99)}     \\
              $|\Theta|=10$ & 0.177 {\tiny (9.99)}     & 0.177 {\tiny (9.99)}     & 0.177 {\tiny (9.99)}     \\
              $|\Theta|=15$ & 0.186 {\tiny (14.99)}     & 0.185 {\tiny (14.99)}     & 0.186 {\tiny (14.99)}     \\
              $|\Theta|=20$ & 0.191 {\tiny (19.99)}     & 0.191 {\tiny (19.99)}     & 0.191 {\tiny (19.99)}     \\\bottomrule
            \end{tabular}
        \end{minipage}
    \end{tabular}
\end{table}

\subsection{Real world data example: argmax inference}
\label{Sec: real-data-argmax}

In this subsection, we continue our analysis of the classification competition datasets introduced in \Cref{Sec: real data}. Here, our focus shifts to the argmax inference problem, where the objective is to construct confidence sets identifying the worst-performing model, characterized by the highest classification loss.

For data preprocessing, we introduce Gaussian noise with mean zero and variance $10^{-60}$ to the classification losses to mitigate numerical instability. Additionally, to create a more challenging inference scenario, we exclude specific teams: teams numbered $12$ and $33$ from the 2023 dataset, and teams numbered $10, 15,$ and $32$ from the 2024 dataset.

\begin{figure}[h]
\centering
\includegraphics[width=1\textwidth]{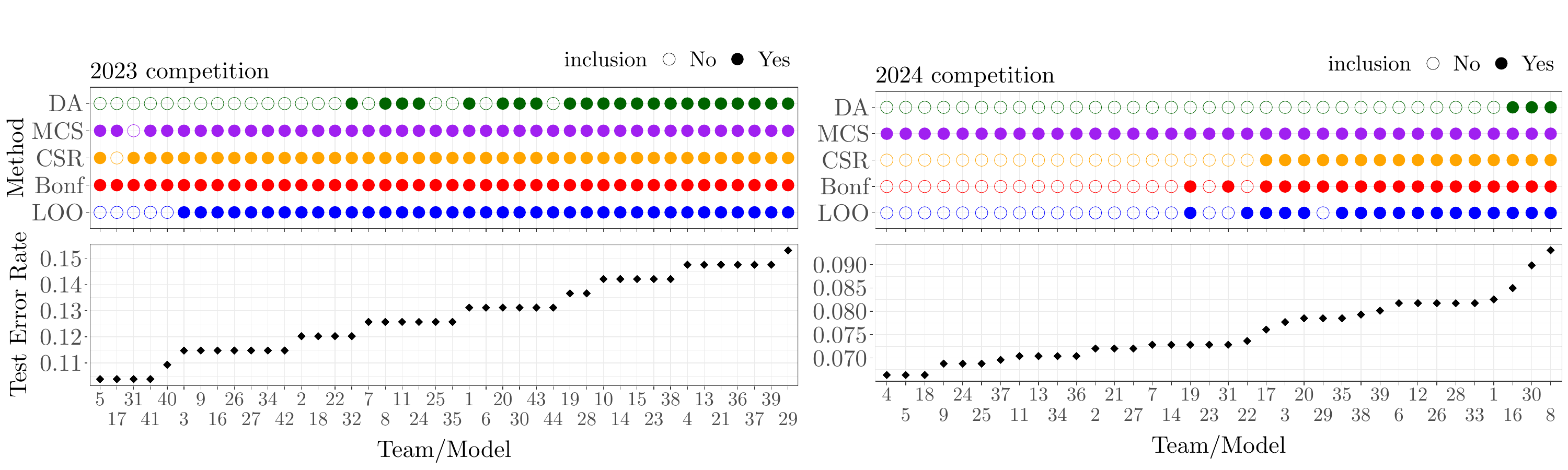}
\caption{
Comparison of inclusion sets generated by the proposed \texttt{DA-adj}$^{\times 50}$ method (DA) and other established techniques across the 2023 (left) and 2024 (right) classification competitions. The competing methods are MCS (\texttt{MCS}), CSR (\texttt{csranks}), Bonf (\texttt{Bonferroni}), and LOO (\texttt{LOO}). Each inclusion set is depicted as a colored interval. Our \texttt{DA-adj}$^{\times 50}$ method consistently produces smaller inclusion sets for the argmax, indicating its enhanced precision in identifying the model with the highest classification loss.
}
\label{fig:real-data-argmax}
\end{figure}

Similar to the visualizations provided in \Cref{fig:real-data}, \Cref{fig:real-data-argmax} illustrates representative inclusion sets generated by various argmax inference methods. The results clearly demonstrate that the \texttt{DA-adj}$^{\times 50}$ method consistently yields smaller inclusion sets compared to other established approaches. This observation aligns with the argmin inference results discussed in \Cref{Sec: real data} and further underscores the superior efficiency of the \texttt{DA-adj}$^{\times 50}$ approach in accurately pinpointing the best-performing models.

\end{document}